\documentclass[12pt]{amsart}
\usepackage{amssymb,amscd}
\usepackage[dvipdfmx]{graphicx,color}
\usepackage{graphicx}

\newtheorem*{Whitney towers}{Theorem~\ref{Whitney towers}}
\newtheorem*{h-towers}{Theorems ~\ref{half} \& \ref{$(n)$-solvable}}

\newtheorem*{surgery curves}{Theorem~\ref{surgery curves}}
\newtheorem*{cg=0}{Theorem~\ref{vanish}}

\newtheorem{thm}{Theorem}[section]

\newtheorem{pr}[thm]{Proposition} 
\newtheorem{lem}[thm]{Lemma}
\newtheorem{cor}[thm]{Corollary}

\newtheorem{cla}[thm]{Claim}

\theoremstyle{definition}
\newtheorem{defn}[thm]{Definition}

\newtheorem{prob}[thm]{Problem}
\newtheorem{note}[thm]{Note}

\numberwithin{equation}{section}
\numberwithin{figure}{section}

\newcommand{\x}{\times}
\newcommand{\np}{\newpage}
\newcommand{\Z}{\mathbb{Z}}
\newcommand{\N}{\mathbb{N}}

\newcommand{\Q}{\mathbb{Q}}

\def\yen{{\setbox0=\hbox{Y}Y\kern-.97\wd0\vbox{hrule height.lex width.98%
\wd0\kern.33ex\hrule height.lex width.98\wd0\kern.45ex}}}%%

\def\np{\newpage}

\begin{document}
\pagestyle{plain}
%{\color{cyan}{  
\title{
Local-move-identities for the $\Z[t,t^{-1}]$-Alexander polynomials of 2-links, 
the alinking number, 
and 
high dimensional analogues
}
\author{
  Eiji Ogasa
}

%{\color{cyan}{  }}

%\thanks{\hskip-4mm %E-mail: %kauffman@uic.edu\quad  
%ogasa@mail1.meijigakuin.ac.jp \newline
%Keywords: Virtual knots, Tori in $S^4$ }
% {\bf PACS nos.} 11-25w, 11-25Uv.
%\newline MSC2000 57N10, 57N13, 57N15. 
%\newline Differential Topology

\date{}

%{\color{cyan}{  }}

\begin{abstract}
A well-known identity 
$\hat\Delta_{L_+}-\hat\Delta_{L_-}=(t^\frac{1}{2}-t^\frac{-1}{2})\cdot\hat\Delta_{L_0}$ 
holds 
for three 1-links  $L_+$, $L_-$, and $L_0$ which satisfy a famous local-move-relation, 
where $\hat\Delta_{L}$ becomes the Alexander-Conway polynomial of %a 1-link 
$L$ 
if we let $z=t^\frac{1}{2}-t^\frac{-1}{2}$. 
We prove 
a new local-move-identity for the $\Z[t,t^{-1}]$-Alexander polynomials of 2-links, 
which is a 2-dimensional analogue of the 1-dimensional one.  
In the 1-dimensional link case there is a well-known relation between  
 the Alexander-Conway polynomial and the linking number.  
As its 2-dimensional analogue, 
we find a relation between the $\Z[t,t^{-1}]$-Alexander polynomials of 2-links and the alinking number of 2-links. 
We show high dimensional analogues of these results. 
Furthermore we prove 
that in the 2-dimensional case we cannot normalize the $\Z[t,t^{-1}]$-Alexander polynomials to be compatible with our identity  
but that in a high-dimensional case  
we can do that to be compatible with our new identity. 
\end{abstract} 
\maketitle

%\begin{figure}
%   \includegraphics[width=14cm]{ }
%\end{figure}

\section{Introduction and main results}\label{intro}
\noindent 
Suppose that three 1-dimensional links $L_+$, $L_-$, and $L_0$ $\subset S^4$   
 differ only in a 3-ball $B$ as shown below.  

\bigbreak
\hskip2cm%WinTpicVersion2.15
\unitlength 0.1in
\begin{picture}(46.22,14.00)(8.10,-22.10)
% CIRCLE 2 0 3 0
% 4 1356 1756 1716 2166 1716 2166 1696 2176
% 
\special{pn 8}%
\special{ar 1356 1356 546 546  0.8902751 6.2831853}%
\special{ar 1356 1356 546 546  0.0000000 0.8502422}%
% VECTOR 2 0 3 0
% 2 946 2116 1736 1386
% 
\special{pn 8}%
\special{pa 946 1716}%
\special{pa 1736 986}%
\special{fp}%
\special{sh 1}%
\special{pa 1736 986}%
\special{pa 1673 1017}%
\special{pa 1697 1022}%
\special{pa 1701 1046}%
\special{pa 1736 986}%
\special{fp}%
% VECTOR 2 0 3 0
% 2 1266 1676 976 1376
% 
\special{pn 8}%
\special{pa 1266 1276}%
\special{pa 976 976}%
\special{fp}%
\special{sh 1}%
\special{pa 976 976}%
\special{pa 1008 1038}%
\special{pa 1013 1014}%
\special{pa 1037 1010}%
\special{pa 976 976}%
\special{fp}%
% LINE 2 0 3 0
% 2 1416 1816 1726 2146
% 
\special{pn 8}%
\special{pa 1416 1416}%
\special{pa 1726 1746}%
\special{fp}%
% CIRCLE 2 0 3 0
% 4 3086 1796 3446 2206 3446 2206 3426 2216
% 
\special{pn 8}%
\special{ar 3086 1396 546 546  0.8902751 6.2831853}%
\special{ar 3086 1396 546 546  0.0000000 0.8502422}%
% VECTOR 2 0 3 0
% 2 2996 1716 2706 1416
% 
\special{pn 8}%
\special{pa 2996 1316}%
\special{pa 2706 1016}%
\special{fp}%
\special{sh 1}%
\special{pa 2706 1016}%
\special{pa 2738 1078}%
\special{pa 2743 1054}%
\special{pa 2767 1050}%
\special{pa 2706 1016}%
\special{fp}%
% LINE 2 0 3 0
% 2 3146 1856 3456 2186
% 
\special{pn 8}%
\special{pa 3146 1456}%
\special{pa 3456 1786}%
\special{fp}%
% LINE 2 0 3 0
% 4 3210 1920 2990 1720 3240 1950 3100 1820
% 
\special{pn 8}%
\special{pa 3210 1520}%
\special{pa 2990 1320}%
\special{fp}%
\special{pa 3240 1550}%
\special{pa 3100 1420}%
\special{fp}%
% VECTOR 2 0 3 0
% 2 3140 1730 3480 1410
% 
\special{pn 8}%
\special{pa 3140 1330}%
\special{pa 3480 1010}%
\special{fp}%
\special{sh 1}%
\special{pa 3480 1010}%
\special{pa 3418 1041}%
\special{pa 3441 1047}%
\special{pa 3445 1070}%
\special{pa 3480 1010}%
\special{fp}%
% LINE 2 0 3 0
% 2 3010 1870 2690 2170
% 
\special{pn 8}%
\special{pa 3010 1470}%
\special{pa 2690 1770}%
\special{fp}%
% CIRCLE 2 0 3 0
% 4 4886 1796 5246 2206 5246 2206 5226 2216
% 
\special{pn 8}%
\special{ar 4886 1396 546 546  0.8902751 6.2831853}%
\special{ar 4886 1396 546 546  0.0000000 0.8502422}%
% VECTOR 2 0 3 0
% 2 4570 2230 4540 1370
% 
\special{pn 8}%
\special{pa 4570 1830}%
\special{pa 4540 970}%
\special{fp}%
\special{sh 1}%
\special{pa 4540 970}%
\special{pa 4522 1037}%
\special{pa 4542 1023}%
\special{pa 4562 1036}%
\special{pa 4540 970}%
\special{fp}%
% VECTOR 2 0 3 0
% 2 5240 2220 5230 1370
% 
\special{pn 8}%
\special{pa 5240 1820}%
\special{pa 5230 970}%
\special{fp}%
\special{sh 1}%
\special{pa 5230 970}%
\special{pa 5211 1037}%
\special{pa 5231 1023}%
\special{pa 5251 1036}%
\special{pa 5230 970}%
\special{fp}%
% STR 2 0 3 0
% 3 1100 2640 1100 2740 2 0
% $L_+$
\put(11.0000,-23.4000){\makebox(0,0)[lb]{$L_+$}}%
% STR 2 0 3 0
% 3 2850 2640 2850 2740 2 0
% $L_{-}$
\put(28.5000,-23.4000){\makebox(0,0)[lb]{$L_{-}$}}%
% STR 2 0 3 0
% 3 4700 2680 4700 2780 2 0
% $L_0$
\put(47.0000,-23.8000){\makebox(0,0)[lb]{$L_0$}}%
\end{picture}%  
\bigbreak

\noindent It is very well-known that then we have the identity $(*)$ $$\hat\Delta_{L_+}-\hat\Delta_{L_-}=(t^\frac{1}{2}-t^\frac{-1}{2})\cdot\hat\Delta_{L_0},$$ 
where $\hat\Delta_{L}$ denotes 
the normalized Alexander polynomial of $L$. 
(If we let $z=t^\frac{1}{2}-t^\frac{-1}{2}$,  $\hat\Delta_{L}$ becomes the Alexander-Conway polynomial.) 
%(We mean that the normalized Alexander polynomial is the same as the Alexander-Conway polynomial. )
It is also very well-known that 
the Jones polynomial satisfies a similar local-move-identity 
and that 
there are several relations between local-moves on 1-links and their invariants.
(See \cite{Alexander, Cn, J, Kauffmanon, KauffmanS} etc.)

\bigbreak 
In \cite[Theorem 4.1]{Ogasa09} we showed a 2-link version of the identity $(*)$. 
We cite it after we state its corollary and an example. 
We strength 
it %\cite[Theorem 4.1]{Ogasa09} 
and obtain a main result of this paper, Theorem \ref{two}, cited in several paragraphs.   
%, which is a generalization of \cite[Theorem 4.1]{Ogasa09} cited in the few previous paragraphs.   

The corollary is as follows: 
Suppose that two 2-dimensional spherical knots, $K_+$ and $K_-$, $\subset S^4$ and 
a submanifold $K_0$$\subset S^4$  
differ only in a 4-ball $B$ trivially embedded in $S^4$ as shown in 
Figures \ref{rev}.1-\ref{rev}.3 in \S\ref{rev}.   
(This ordered set $(K_+, K_-,K_0)$ is called a (1,2)-pass-move-triple. 
An example %of $(K_+, K_-,K_0)$ 
is made from Figure \ref{intro}.1, \ref{intro}.3, and \ref{intro}.5,  %in this section, 
which are explained in the following paragraph. 
See \S\ref{rev} for the precise definition.)   
%of 
%%$(1,2)$-pass-moves and 
%$(1,2)$-pass-move-triples.)  
%An example of $(1,2)$-pass-move-triples is made from Figure \ref{intro}.1, \ref{intro}.3, and \ref{intro}.5 in this section as we explain in the following paragraph.)  
Then we have the following $(\#)$: 
there is a polynomial $\Delta_{\nu, K_*}(t) (*=+,-,0$ and  $\nu=1,2)$ 
which represents 
the $\Q[t,t^{-1}]$-$\nu$-Alexander polynomial for $K_*$,  and  
we have the identity \newline
$\Delta_{\nu, L_+}(t)-\Delta_{\nu, L_-}(t)=(t-1)\cdot\Delta_{\nu, L_0}(t).$

$K_+$ in Figure \ref{intro}.1,  $K_-$ in  Figure \ref{intro}.3, and $K_0$ in Figure \ref{intro}.5    
% make an example of $(1,2)$-pass-move-triple, and 
are constructed as follows: 
Embed $F=(S^1\x S^2)-$open$B^3$ in $S^4$.  
The 
boundary of $F$ in  $S^4$ is a $2$-knot.  
Let it be a trivial 2-knot $K_+$ as drawn in Figure \ref{intro}.1. 
Carry out a `local-move' on the $2$-knot $K_+$ 
in a 4-ball, which is denoted by a dotted circle  
in Figure \ref{intro}.2. This local-move is called the (1,2)-pass-move (see \S2 for the precise definition). 
Note that the above operation is done only in the 4-ball.  
The (1,2)-pass-move changes the trivial 2-knot in Figure \ref{intro}.1 (resp. \ref{intro}.2) 
into a 2-knot in Figure \ref{intro}.3. 
We can prove that the knot in Figure \ref{intro}.3 is nontrivial  
by using Seifert matrices and the Alexander polynomial. 
We use the fact that $S^1$ and $S^2$ can be `linked' in $S^4$. 
Note that $S^1$ and $S^2$ are included in $F$ as shown in Figure \ref{intro}.4.  
$K_0$ is drawn in Figure \ref{intro}.5.
If we give appropriate orientations, we can let 
$\Delta_{K_+}=t$,  
$\Delta_{K_-}=2t-1$, and $\Delta_{K_0}=-1$. 
Hence $\Delta_{K_+}-\Delta_{K_-}=(t-1)\cdot\Delta_{K_0}$ holds.

\cite[Theorem 4.1]{Ogasa09} is as follows: 
Let $L_+=(L_{+,1},...,L_{+,{m_+}})$ be a 2-dimensional closed oriented submanifold  $\subset S^4$. 
Let each $L_{+,i}$ be connected. 
Let $g_{+,i}$ be the genus of $L_{+,i}$. 
Let $m_+=1+\Sigma_1^{m_+} g_{+,i}$. 
Let $(L_+, L_-, L_0)$ be a $(1,2)$-pass-move-triple. 
%Suppose that $L_+, L_-,$ and $L_0$ differ only in a 4-ball $B$ trivially embedded in $S^4$ as shown in Figures \ref{rev}.1-\ref{rev}.3 in \S\ref{rev}.   
%%(This ordered set $(L_+, L_-, L_0)$ is called a (1,2)-pass-move-triple. See \S\ref{rev} for the precise definition of $(1,2)$-pass-moves and $(1,2)$-pass-move-triples.An example of $(1,2)$-pass-move-triples is made from Figure \ref{intro}.1, \ref{intro}.3, and \ref{intro}.5 
%%%in this section
%%as we explain in a few paragraphs.)  
%Then there is a polynomial $\Delta_{\nu, K_*}(t) (*=+,-,0)$ which represents the $\Q[t,t^{-1}]$-$\nu$-Alexander polynomial for $L_*$,  and  we have the following identity $(\nu=1,2)$  $\Delta_{\nu, L_+}(t)-\Delta_{\nu, L_-}(t)=(t-1)\cdot\Delta_{\nu, L_0}(t).$ 
%%where $\Delta_{K}$ represents the $\Q[t,t^{-1}]$-Alexander polynomial of $K$. 
Then we have the above $(\#)$, where we replace $K_*$ with $L_*$. 

In \cite[Proposition 4.3]{Ogasa09} we proved that we cannot normalize 
the Alexander polynomials  
%$\Delta_{K_+}$, $\Delta_{K_-}$, and $\Delta_{K_0},$ 
to be compatible with this local-move-identity.   

In \cite{KauffmanOgasa,  KauffmanOgasaII,  KauffmanOgasaB,   
Ogasa98n, Ogasa04, Ogasa07, Ogasa09, OgasaT3, OgasaIH},  
furthermore, 
we proved several relations between local-moves on $n$-knots and their invariants ($n\in\N$).

\bigbreak 
We state one of our main results, Theorem \ref{two}, 
%A main result of this paper is Theorem \ref{two}. 
%, which is a generalization of \cite[Theorem 4.1]{Ogasa09} cited in the few previous paragraphs.   
whose $\Q[t,t^{-1}]$-Alexander polynomial case is  \cite[Theorem 4.1]{Ogasa09}.    
The former is stronger than the latter because 
the former does not follow from the latter directly 
(see  Propositions \ref{Idaho} and \ref{Mississippi}).  
Furthermore we prove a high dimensional analogue of this main theorem 
(see Theorem \ref{Tokyo} quoted in this section.) 

%The $\Q[t,t^{-1}]$-Alexander polynomial case of Theorem \ref{two} is  \cite[Theorem 4.1]{Ogasa09}.    Theorem \ref{two} is stronger than \cite[Theorem 4.1]{Ogasa09} because Theorem \ref{two} does not follow from \cite[Theorem 4.1]{Ogasa09} directly (see  Propositions \ref{Idaho} and \ref{Mississippi}).  Furthermore we prove a high dimensional analogue of this theorem (see Theorem \ref{Tokyo} quoted in this section.) 

\bigbreak
\noindent   
{\bf Theorem \ref{two}.} {\it 
Let $L_+=(L_{+,1},...,L_{+,{m_+}})$ be a 2-dimensional closed oriented submanifold  $\subset S^4$. 
Let each $L_{+,i}$ be connected. 
Let $g_{+,i}$ be the genus of $L_{+,i}$. 
Let $m_+=1+\Sigma_1^{m_+} g_{+,i}$. 
Let $(L_+, L_-, L_0)$ be a $(1,2)$-pass-move-triple. 
Then there is a polynomial $\Delta_{\nu, K_*}(t) (*\newline=+,-,0$ and  $\nu=1,2)$ which represents 
the $\Z[t,t^{-1}]$-$\nu$-Alexander polynomial for $K_*$,  and  
we have the identity %$(\nu=1,2)$

$$\Delta_{\nu, L_+}(t)-\Delta_{\nu, L_-}(t)=(t-1)\cdot\Delta_{\nu, L_0}(t).$$ 
}
\bigbreak

$K_+$ in Figure \ref{intro}.1,   $K_-$ in  Figure \ref{intro}.3, and $K_0$ in Figure \ref{intro}.5 
make not only an example of \cite[Theorem 4.1]{Ogasa09} 
but also Theorem \ref{two}. 

%We show an example of the identity in Theorem \ref{two} here.It is also an example of the identity in  \cite[Theorem 4.1]{Ogasa09}.   $K_+$ in Figure \ref{intro}.1,   $K_-$ in  Figure \ref{intro}.3, and $K_0$ in Figure \ref{intro}.5  satisfy the identity $\Delta_{K_+}-\Delta_{K_-}=(t-1)\cdot\Delta_{K_0},$  where $\Delta_{K}$ represents the $\Z[t,t^{-1}]$-Alexander polynomial of 2-dimensional closed oriented submanifolds $\subset S^4$.   These $K_+, K_-, K_0$ make an example of $(1,2)$-pass-move-triple.  Embed $F=(S^1\x S^2)-$open$B^3$ in $S^4$.  The \newline boundary of $F$ in  $S^4$ is a $2$-knot.  Let it be a trivial 2-knot $K_+$ as drawn in Figure \ref{intro}.1. Carry out a `local-move' on the $2$-knot $K_+$ in a 4-ball, which is denoted by a dotted circle  in Figure \ref{intro}.2. This local-move is called the (1,2)-pass-move (see \S2 for the precise definition). Note that the above operation is done only in the 4-ball.  The (1,2)-pass-move changes the trivial 2-knot in Figure \ref{intro}.1 (resp. \ref{intro}.2) into a 2-knot in Figure \ref{intro}.3. We can prove that the knot in Figure \ref{intro}.3 is nontrivial by using Seifert matrices and the Alexander polynomial. We use the fact that $S^1$ and $S^2$ can be `linked' in $S^4$. Note that $S^1$ and $S^2$ are included in $F$ as shown in Figure \ref{intro}.4.  $K_0$ is drawn in Figure \ref{intro}.5.  If we give appropriate orientations, we can let $\Delta_{K_+}=t$,  $\Delta_{K_-}=2t-1$, and $\Delta_{K_0}=-1$. Hence $\Delta_{K_+}-\Delta_{K_-}=(t-1)\cdot\Delta_{K_0}$ holds. 

\begin{figure}
\hskip3cm\includegraphics[width=35mm]{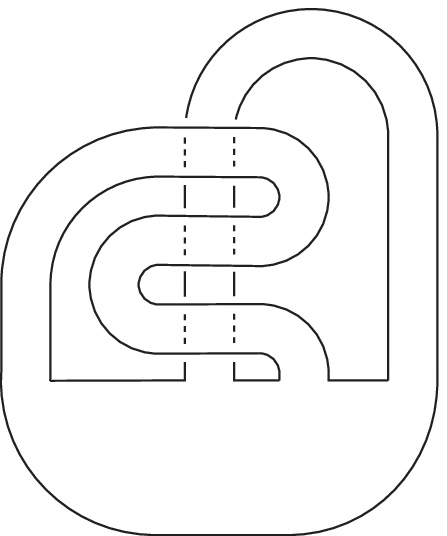}   

\hskip4cm   Figure \ref{intro}.1: {\bf A trivial $2$-knot $K_+$ and $F$. $\partial F=K_+.$} 
\end{figure}

\begin{figure}
\hskip4cm\includegraphics[width=35mm]{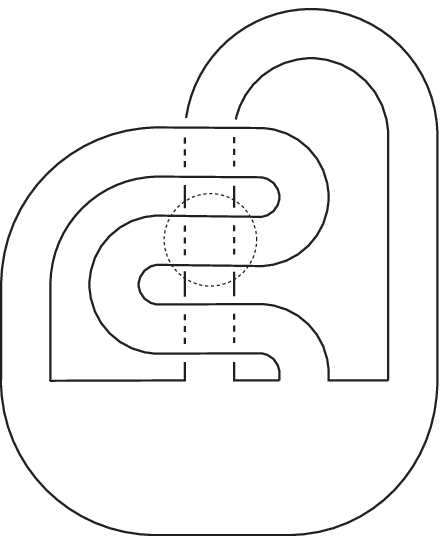}  %%zz-F3
\smallbreak
\hskip2cm   Figure \ref{intro}.2: {\bf A local-move will be carried out in the dotted  

\hskip20mm 4-ball. The resulting 2-knot is a nontrivial 2-knot $K_-$.}  
\end{figure}

\begin{figure}
\hskip4cm\includegraphics[width=35mm]{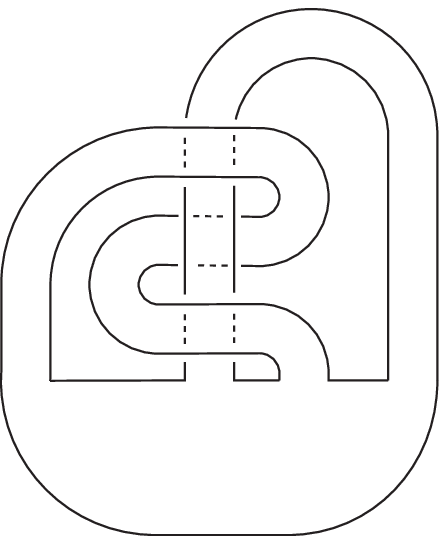}  
\smallbreak
\hskip3cm   Figure \ref{intro}.3: {\bf A nontrivial 2-knot $K_-$}

\end{figure}

\begin{figure}
\hskip4cm\includegraphics[width=40mm]{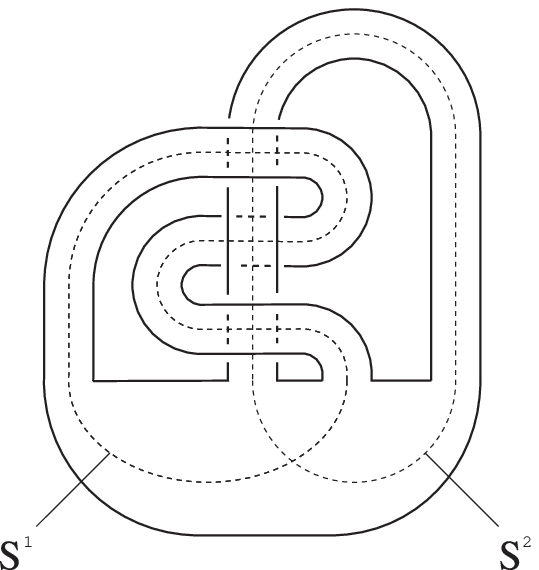}  
\smallbreak
\hskip2cm   Figure \ref{intro}.4: 
{\bf $S^1$ and $S^2$ in $F$ whose boundary is the $2$-knot} 
\end{figure}

\begin{figure}
\hskip4cm\includegraphics[width=35mm]{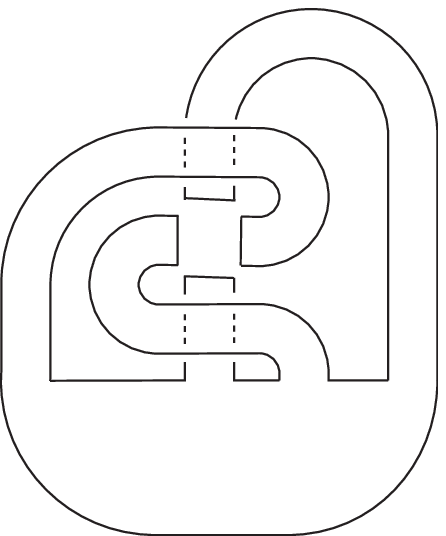}  
\smallbreak
\hskip3cm   Figure \ref{intro}.5: {\bf $K_0$ is a trivial 2-knot in this case.}
\end{figure}

\bigbreak 
Theorem \ref{Tokyo} is as follows. 
The terms and definitions needed for it  % Theorem \ref{Tokyo} 
are in the body of the paper. 
In the 2-dimensional case we cannot normalize 
the $\Z[t,t^{-1}]$-Alexander polynomial to be compatible with the local-move-identity  
%by  \cite[Proposition 4.3]{Ogasa09} 
as we explained a few paragraphs before,  
but in a case of $(4k+1)$-dimensional case 
we can define the `normalized' Alexander polynomial 
(see Definition \ref{4k+1} and Theorem \ref{norm})  
associated with a local-move defined in \S\ref{Osaka}. 

\bigbreak
\noindent
{\bf Theorem \ref{Tokyo}.} {\it 
Let $K_+$ be a $(4k+1)$-knot $\subset S^{4k+3}$. 
Let $(K_+, K_-, K_0)$ be a twist-move-triple.
Then 
$$\hat\Delta_{K_+}(t)-\hat\Delta_{K_-}(t)
=(t^\frac{1}{2}-t^\frac{-1}{2})\cdot\hat\Delta_{K_0}(t),$$
where $\hat\Delta_{K}(t)$ denotes the normalized Alexander polynomial of $K$. 
}

\bigbreak
In the 1-dimensional case we have the following fact (\cite{Hoste}): 
Let $L$ be a 2-component 1-link. 
Let $\hat\Delta_{L}(t)$ be the normalized Alexander polynomial of $L$. 
Then 
$\left.\displaystyle\frac{\hat\Delta_{K}(t)}{t^\frac{1}{2}-t^\frac{-1}{2}}\right|_{t=1}$ 
%$\Delta_{L}(1)$ 
is the linking number of $L$. 
Let $K_+,  K_-,$ and $K_0$ be as in the first paragraph of this section. 
Let $K_+$ be a 1-knot. 
Then $K_0$ is a 2-component 1-link and 
$\left.\displaystyle
\frac{ \hat\Delta_{K_+}(t)-\hat\Delta_{K_-}(t)}{(t^\frac{1}{2}-t^\frac{-1}{2})^2}
\right|_{t=1}$
%$\left.\displaystyle\frac{\Delta_{K_+}(t)-\Delta_{K_-}(t)}{t^\frac{1}{2}-t^\frac{-1}{2}}\right|_{t=1}$
is the linking number of $K_0$.

%In this paper we prove the 2-dimensional analogue of this result and the high-dimensional one. 
In \S\ref{m2}  we prove a 2-dimensional analogue of this result:  
we show a relation between 
the $\Z[t,t^{-1}]$-Alexander polynomial of 
2-dimensional closed oriented submanifolds $\subset S^4$ 
and the alinking number (Theorems \ref{alk} and \ref{Georgia}). 
In \S\ref{mhigh} furthermore we prove some high dimensional analogues 
(Theorem \ref{aletwi} and Corollary \ref{twi+-}). 
We cite the above theorems here. 
The terms and definitions needed for them are in the body of the paper. 

\bigbreak
\noindent
{\bf Theorem \ref{alk}.} {\it 
Let $L=(L_1, L_2)$ be an $(S^2, T^2)$-link $\subset S^4$. 
Let $\Delta_{1,L}(t)$ be a polynomial which represents 
the $\Z[t,t^{-1}]$-1-Alexander polynomial for $L$. 
Then 
$$
\begin{vmatrix}
\left.\displaystyle\frac{\Delta_{1,L}(t)}{(t-1)}\right|_{t=1}
\end{vmatrix}$$ 
is the pseudo-alinking number of $L$,  where $|\quad|$ denotes the absolute value.  
}
\bigbreak

\noindent
{\bf Theorem \ref{Georgia}.} {\it 
Let $L=(K_1,K_2)$ be a ribbon $(S^2, T^2)$-link. 
Then the alinking number of $L$ is 
$
\begin{vmatrix}
\left.\displaystyle\frac{\Delta_{1,L}(t)}{(t-1)}\right|_{t=1}
\end{vmatrix}$. 
}
\bigbreak

\noindent
{\bf Theorem \ref{aletwi}.} {\it 
Let $K$ be a $(4k+1)$-dimensional closed oriented subamanifold $\subset S^{4k+3}$ 
whose homotopy type is $S^{2k}\x S^{2k+1}$.  
Let $\hat\Delta_{K}(t)$ be the normalized Alexander polynomial of $K$.  
Then the pseudo-twinkling number of $K$ is 
$$\left.\displaystyle\frac{\hat\Delta_{K}(t)}{t^\frac{1}{2}-t^\frac{-1}{2}}\right|_{t=1}$$
}
\bigbreak

\noindent
{\bf Corollary \ref{twi+-}.} {\it 
Let $K_+$ be a $(4k+1)$-knot $\subset S^{4k+3}$. 
Let $(K_+, K_-, K_0)$ be a twist-move-triple.
Then the pseudo-twinkling number of $K_0$ is 
$$\left.\displaystyle
\frac{ \hat\Delta_{K_+}(t)-\hat\Delta_{K_-}(t)}{(t^\frac{1}{2}-t^\frac{-1}{2})^2}
\right|_{t=1}$$
}

\bigbreak

$$ \text{Table of  Contents} $$

\S\ref{intro} Introduction and main results

\S\ref{rev} Review of (1,2)-pass-moves on 2-knots 

\S\ref{revale}  
Review of the $\Q[t,t^{-1}]$-Alexander polynomial for 
(not necessarily connected) 

\hskip5mm
$n$-dimensional closed oriented submanifolds 
in $S^{n+2}$ $(n\geqq2)$

\S\ref{m2}  
Main theorems in the 2-dimensional case  

\S\ref{Osaka}
Review of twist-moves on high dimensional knots 

\S\ref{mhigh} Main theorems in the 4k+1 dimensional case 

\S\ref{proo2} Proof of results of \S\ref{m2} 

\S\ref{proohigh} Proof of results of \S\ref{mhigh}

\S\ref{problem} A problem

\bigbreak
\section{Review of (1,2)-pass-moves on 2-knots}\label{rev}
\noindent 
The local-move which associates the identity $(*)$ in the first paragraph in \S\ref{intro}  
is very easy as drawn there. 
In high dimensional case we must begin by explaining what kind of local-moves we use. 
We review the (1,2)-pass-move on  
2-dimensional closed oriented submanifolds $\subset S^4$, 
which are defined in  \cite{Ogasa04}.

We work in the smooth category. 
A (not necessarily connected) 2-dimensional smooth, closed oriented submanifold 
$L \subset S^4$ is called 
an  {\it $m$-component 2-$($dimensional$)$ $($spherical$)$ link}
if $L$ consists of $m$ connected components and 
each connected component is a 2-sphere. 
 If $L$ is 1-component 2-link, then $L$ is called a {\it $($spherical$)$ 2-knot}.  
We say that 
(not necessarily connected) 2-dimensional smooth, closed, oriented submanifolds
 $L_1$ and $L_2 \subset S^4$ are {\it equivalent} 
if there exists an orientation preserving diffeomorphism 
$f:$ $S^4$ $\rightarrow$ $S^4$  
such that $f(L_1)$=$L_2$  and 
that $f | _{L_1}:L_1\rightarrow L_2$ can be regarded 
as an order and orientation preserving diffeomorphism map.    

The (1,2)-pass-move is a local-move. 
Here, `local-move' means that 
when we make $K_+$ into $K_-$ (resp. $K_-$ into $K_0$, $K_0$ into $K_-$) 
vice versa in Definition \ref{12},   
we make a change only in $B$ and 
that we do not any requirement on diffeomorphism type or homeomorphism type 
of $K_+$ (resp. $K_-$, $K_0$) other than the change only in $B$

\begin{figure}
\input 4.1.tex
\bigbreak
\bigbreak
\bigbreak
\end{figure}

\begin{figure}
\input 4.2.tex
\bigbreak
\bigbreak
\end{figure}

\begin{figure}
%WinTpicVersion2.15
\unitlength 0.1in
\begin{picture}(37.87,32.10)(8.50,-32.70)
% ELLIPSE 2 0 3 0
% 4 2646 636 3024 804 3024 804 3024 804
% 
\special{pn 8}%
\special{ar 2646 236 378 168  0.0000000 6.2831853}%
% SPLINE 0 0 1 0
% 41 2538 3120 2538 3111 2542 3103 2548 3095 2557 3088 2567 3081 2580 3076 2592 3072 2608 3068 2623 3067 2639 3066 2655 3067 2670 3068 2684 3072 2699 3076 2710 3081 2720 3088 2729 3095 2736 3103 2738 3111 2740 3120 2738 3129 2736 3137 2729 3144 2720 3151 2710 3158 2699 3164 2684 3167 2670 3171 2655 3173 2639 3174 2623 3173 2608 3171 2592 3167 2580 3164 2567 3158 2557 3151 2548 3144 2542 3137 2538 3129 2538 3120
% 
\special{sh 0.300}%
\special{pn 20}%
\special{pa 2538 2720}%
\special{pa 2551 2692}%
\special{pa 2579 2676}%
\special{pa 2609 2668}%
\special{pa 2641 2666}%
\special{pa 2673 2669}%
\special{pa 2704 2678}%
\special{pa 2730 2696}%
\special{pa 2739 2725}%
\special{pa 2721 2750}%
\special{pa 2693 2765}%
\special{pa 2662 2772}%
\special{pa 2630 2774}%
\special{pa 2599 2769}%
\special{pa 2568 2759}%
\special{pa 2543 2739}%
\special{pa 2538 2720}%
\special{sp}%
% ELLIPSE 0 0 0 0
% 4 2639 629 2740 683 3307 1014 3307 1014
% 
\special{pn 20}%
\special{sh 0.600}%
\special{ar 2639 229 101 54  0.0000000 6.2831853}%
% ELLIPSE 2 0 3 0
% 4 2646 3174 3024 3342 3024 3342 3024 3342
% 
\special{pn 8}%
\special{ar 2646 2774 378 168  0.0000000 6.2831853}%
% LINE 2 0 3 0
% 2 3030 642 3030 3174
% 
\special{pn 8}%
\special{pa 3030 242}%
\special{pa 3030 2774}%
\special{fp}%
% LINE 2 0 3 0
% 2 2268 649 2268 3174
% 
\special{pn 8}%
\special{pa 2268 249}%
\special{pa 2268 2774}%
\special{fp}%
% ELLIPSE 2 0 3 0
% 4 1228 642 1606 811 1606 811 1606 811
% 
\special{pn 8}%
\special{ar 1228 242 378 169  0.0000000 6.2831853}%
% ELLIPSE 2 0 3 0
% 4 1228 3180 1606 3349 1606 3349 1606 3349
% 
\special{pn 8}%
\special{ar 1228 2780 378 169  0.0000000 6.2831853}%
% LINE 2 0 3 0
% 2 1613 649 1613 3180
% 
\special{pn 8}%
\special{pa 1613 249}%
\special{pa 1613 2780}%
\special{fp}%
% LINE 2 0 3 0
% 2 850 656 850 3180
% 
\special{pn 8}%
\special{pa 850 256}%
\special{pa 850 2780}%
\special{fp}%
% ELLIPSE 2 0 3 0
% 4 4252 629 4630 798 4630 798 4630 798
% 
\special{pn 8}%
\special{ar 4252 229 378 169  0.0000000 6.2831853}%
% ELLIPSE 2 0 3 0
% 4 4252 3167 4630 3336 4630 3336 4630 3336
% 
\special{pn 8}%
\special{ar 4252 2767 378 169  0.0000000 6.2831853}%
% LINE 2 0 3 0
% 2 4637 636 4637 3167
% 
\special{pn 8}%
\special{pa 4637 236}%
\special{pa 4637 2767}%
\special{fp}%
% LINE 2 0 3 0
% 2 3874 642 3874 3167
% 
\special{pn 8}%
\special{pa 3874 242}%
\special{pa 3874 2767}%
\special{fp}%
% STR 2 0 3 0
% 3 1086 3457 1086 3525 2 0
% t=-0.5
\put(10.8600,-31.2500){\makebox(0,0)[lb]{t=-0.5}}%
% STR 2 0 3 0
% 3 2436 3457 2436 3525 2 0
% t=0
\put(24.3600,-31.2500){\makebox(0,0)[lb]{t=0}}%
% STR 2 0 3 0
% 3 4056 3457 4056 3525 2 0
% t=0.5
\put(40.5600,-31.2500){\makebox(0,0)[lb]{t=0.5}}%
% ELLIPSE 2 0 3 0
% 4 2646 636 3024 804 3024 804 3024 804
% 
\special{pn 8}%
\special{ar 2646 236 378 168  0.0000000 6.2831853}%
% SPLINE 0 0 0 0
% 41 2538 3120 2538 3111 2542 3103 2548 3095 2557 3088 2567 3081 2580 3076 2592 3072 2608 3068 2623 3067 2639 3066 2655 3067 2670 3068 2684 3072 2699 3076 2710 3081 2720 3088 2729 3095 2736 3103 2738 3111 2740 3120 2738 3129 2736 3137 2729 3144 2720 3151 2710 3158 2699 3164 2684 3167 2670 3171 2655 3173 2639 3174 2623 3173 2608 3171 2592 3167 2580 3164 2567 3158 2557 3151 2548 3144 2542 3137 2538 3129 2538 3120
% 
\special{sh 0.600}%
\special{pn 20}%
\special{pa 2538 2720}%
\special{pa 2551 2692}%
\special{pa 2579 2676}%
\special{pa 2609 2668}%
\special{pa 2641 2666}%
\special{pa 2673 2669}%
\special{pa 2704 2678}%
\special{pa 2730 2696}%
\special{pa 2739 2725}%
\special{pa 2721 2750}%
\special{pa 2693 2765}%
\special{pa 2662 2772}%
\special{pa 2630 2774}%
\special{pa 2599 2769}%
\special{pa 2568 2759}%
\special{pa 2543 2739}%
\special{pa 2538 2720}%
\special{sp}%
% ELLIPSE 0 0 0 0
% 4 2639 629 2740 683 3307 1014 3307 1014
% 
\special{pn 20}%
\special{sh 0.600}%
\special{ar 2639 229 101 54  0.0000000 6.2831853}%
% ELLIPSE 2 0 3 0
% 4 2646 3174 3024 3342 3024 3342 3024 3342
% 
\special{pn 8}%
\special{ar 2646 2774 378 168  0.0000000 6.2831853}%
% LINE 2 0 3 0
% 2 3030 642 3030 3174
% 
\special{pn 8}%
\special{pa 3030 242}%
\special{pa 3030 2774}%
\special{fp}%
% LINE 2 0 3 0
% 2 2268 649 2268 3174
% 
\special{pn 8}%
\special{pa 2268 249}%
\special{pa 2268 2774}%
\special{fp}%
% ELLIPSE 2 0 3 0
% 4 1228 642 1606 811 1606 811 1606 811
% 
\special{pn 8}%
\special{ar 1228 242 378 169  0.0000000 6.2831853}%
% ELLIPSE 2 0 3 0
% 4 1228 3180 1606 3349 1606 3349 1606 3349
% 
\special{pn 8}%
\special{ar 1228 2780 378 169  0.0000000 6.2831853}%
% LINE 2 0 3 0
% 2 1613 649 1613 3180
% 
\special{pn 8}%
\special{pa 1613 249}%
\special{pa 1613 2780}%
\special{fp}%
% LINE 2 0 3 0
% 2 850 656 850 3180
% 
\special{pn 8}%
\special{pa 850 256}%
\special{pa 850 2780}%
\special{fp}%
% ELLIPSE 2 0 3 0
% 4 4252 629 4630 798 4630 798 4630 798
% 
\special{pn 8}%
\special{ar 4252 229 378 169  0.0000000 6.2831853}%
% ELLIPSE 2 0 3 0
% 4 4252 3167 4630 3336 4630 3336 4630 3336
% 
\special{pn 8}%
\special{ar 4252 2767 378 169  0.0000000 6.2831853}%
% LINE 2 0 3 0
% 2 4637 636 4637 3167
% 
\special{pn 8}%
\special{pa 4637 236}%
\special{pa 4637 2767}%
\special{fp}%
% LINE 2 0 3 0
% 2 3874 642 3874 3167
% 
\special{pn 8}%
\special{pa 3874 242}%
\special{pa 3874 2767}%
\special{fp}%
% ELLIPSE 0 0 3 0
% 4 2646 1945 3024 2114 3024 2114 3024 2114
% 
\special{pn 20}%
\special{ar 2646 1545 378 169  0.0000000 6.2831853}%
% BOX 2 5 2 0
% 2 2301 1837 2362 1898
% 
\special{pn 8}%
\special{sh 0}%
\special{pa 2301 1437}%
\special{pa 2362 1437}%
\special{pa 2362 1498}%
\special{pa 2301 1498}%
\special{pa 2301 1437}%
\special{ip}%
% BOX 2 5 2 0
% 2 2922 1824 2983 1884
% 
\special{pn 8}%
\special{sh 0}%
\special{pa 2922 1424}%
\special{pa 2983 1424}%
\special{pa 2983 1484}%
\special{pa 2922 1484}%
\special{pa 2922 1424}%
\special{ip}%
% BOX 2 5 2 0
% 2 2760 1770 2821 1830
% 
\special{pn 8}%
\special{sh 0}%
\special{pa 2760 1370}%
\special{pa 2821 1370}%
\special{pa 2821 1430}%
\special{pa 2760 1430}%
\special{pa 2760 1370}%
\special{ip}%
% BOX 2 5 2 0
% 2 2605 1749 2666 1810
% 
\special{pn 8}%
\special{sh 0}%
\special{pa 2605 1349}%
\special{pa 2666 1349}%
\special{pa 2666 1410}%
\special{pa 2605 1410}%
\special{pa 2605 1349}%
\special{ip}%
% BOX 2 5 2 0
% 2 2450 1770 2511 1830
% 
\special{pn 8}%
\special{sh 0}%
\special{pa 2450 1370}%
\special{pa 2511 1370}%
\special{pa 2511 1430}%
\special{pa 2450 1430}%
\special{pa 2450 1370}%
\special{ip}%
% ELLIPSE 0 0 3 0
% 4 2639 1756 3017 1925 3017 1925 3017 1925
% 
\special{pn 20}%
\special{ar 2639 1356 378 169  0.0000000 6.2831853}%
% LINE 0 0 3 0
% 32 2274 1803 2274 1972 2369 1878 2369 2046 2531 1925 2531 2087 2706 1932 2713 2100 2841 1884 2848 2080 2976 1851 2970 2013 2889 1635 2889 1689 2889 1756 2889 1803 2720 1608 2720 1655 2720 1709 2727 1763 2551 1601 2558 1655 2558 1709 2558 1783 2389 1628 2389 1702 2389 1756 2389 1810 2274 1810 2281 1965 2281 1803 2268 1965
% 
\special{pn 20}%
\special{pa 2274 1403}%
\special{pa 2274 1572}%
\special{fp}%
\special{pa 2369 1478}%
\special{pa 2369 1646}%
\special{fp}%
\special{pa 2531 1525}%
\special{pa 2531 1687}%
\special{fp}%
\special{pa 2706 1532}%
\special{pa 2713 1700}%
\special{fp}%
\special{pa 2841 1484}%
\special{pa 2848 1680}%
\special{fp}%
\special{pa 2976 1451}%
\special{pa 2970 1613}%
\special{fp}%
\special{pa 2889 1235}%
\special{pa 2889 1289}%
\special{fp}%
\special{pa 2889 1356}%
\special{pa 2889 1403}%
\special{fp}%
\special{pa 2720 1208}%
\special{pa 2720 1255}%
\special{fp}%
\special{pa 2720 1309}%
\special{pa 2727 1363}%
\special{fp}%
\special{pa 2551 1201}%
\special{pa 2558 1255}%
\special{fp}%
\special{pa 2558 1309}%
\special{pa 2558 1383}%
\special{fp}%
\special{pa 2389 1228}%
\special{pa 2389 1302}%
\special{fp}%
\special{pa 2389 1356}%
\special{pa 2389 1410}%
\special{fp}%
\special{pa 2274 1410}%
\special{pa 2281 1565}%
\special{fp}%
\special{pa 2281 1403}%
\special{pa 2268 1565}%
\special{fp}%
% STR 2 0 3 0
% 3 1190 3740 1190 3840 2 0
% Figure \ref{rev}.3: $L_0$ of a $(1,2)$-pass-move-triple
\put(11.9000,-34.4000){\makebox(0,0)[lb]{Figure \ref{rev}.3: $L_0$ of a $(1,2)$-pass-move-triple}}%
\end{picture}%
\bigbreak
\bigbreak
\end{figure}

\begin{defn}\label{12} 
Let $L_+$,  
$L_-$, and  
$L_0$ 
be (not necessarily connected) 2-dimensional closed oriented submanifolds $\subset S^4$. 
We say that ($L_+$, $L_-$, $L_0$) is a {\it $(1,2)$-pass-move-triple}   
if $L_+$, $L_-$, and $L_0$ differ only in a 4-ball $B$ trivially embedded in $S^4$ 
with the following properties: 
$B\cap L_+$ is drawn as in Figure \ref{rev}.1. 
$B\cap L_-$ is drawn as in Figure \ref{rev}.2.   
$B\cap L_0$ is drawn as in Figure \ref{rev}.3.     
Note that we do not assume how many connected components of $L_+$ intersect $B$. 
Furthermore we say that $L_+$ (resp. $L_-$) is obtained from $L_-$  (resp. $L_+$) 
by one {\it $(1,2)$-pass-move}. 
If $L$ is equivalent to $L'$ and if $L'$ is obtained from $L''$ by a sequence of 
(1,2)-pass-moves, we say that $L$ is {\it $(1,2)$-pass-move-equivalent} to $L''$.

We regard $B$ as 
(a close 2-disc $P$)$\times[0,1]\times\{t| -1\leqq t\leqq1\}$.
Let 
$B_t=$\newline
(the close 2-disc $P$)$\times[0,1]\times\{t \}$.  
Then $B=\cup B_t$. 
In Figures \ref{rev}.1, \ref{rev}.2 and \ref{rev}.3,  
we draw $B_{-0.5}, B_{0}, B_{0.5}$ $\subset B$. 
We draw $L_+$, $L_-$, and $L_0$ by the bold line. 
The fine line denotes $\partial B_t$.

$B\cap L_+$ (resp. $B\cap L_-$) is diffeomorphic to 
$D^2\amalg D^2\amalg (S^1\times [0,1])$, 
where $\amalg$ denotes the disjoint union. 
$B\cap L_+$ has the following properties:  
$B_t\cap L_+$ is empty for $-1\leqq t<0$ and $0.5<t\leqq1$.
$B_0\cap L_+$ is  
$D^2\times\{0.4\}\amalg D^2\times\{0.6\}\amalg(S^1\times [0,0.3])\amalg(S^1\times [0.7,1])$. 
$B_{0.5}\cap L_+$ is  $S^1\times [0.3,0.7]$. 
$B_t\cap L_+$ is diffeomorphic to $S^1\amalg S^1$ for $0<t<0.5$. 
(Here we draw $S^1\times [0,1]$ to have the corner 
in $B_0$ and in $B_{0.5}$. 
Strictly to say, $B\cap L_+$ in $B$ is a smooth embedding 
which is obtained by making the corner smooth naturally.)

$B\cap L_-$ has the following properties:  
$B_t\cap  L_-$ is empty for $-1\leqq t<-0.5$ and $0<t\leqq1$.
$B_0\cap L_-$ is   
$D^2\times\{0.4\}\amalg D^2\times\{0.6\}\amalg(S^1\times [0, 0.3])
\amalg(S^1\times [0.7, 1])$. 
$B_{-0.5}\cap  L_-$ is  $S^1\times [0.3, 0.7]$. 
$B_t\cap  L_-$ is diffeomorphic to $S^1\amalg S^1$ for $-0.5<t<0$. 
%(Here we draw $S^1\times [0,1]$ to have the corner in $B_0$ and in $B_{-0.5}$. Strictly to say, $B\cap L_-$ in $B$ is a smooth embedding which is obtained by making the corner smooth naturally.)

In Figure \ref{rev}.1 (resp. \ref{rev}.2) 
there are an oriented cylinder $S^1\times [0,1]$ 
and two oriented discs $D^2$. 
We do not make any assumption about 
the orientation of the cylinder.  
We suppose that 
each arrow $\overrightarrow{x}$, $\overrightarrow{y}$ 
in Figure \ref{rev}.1 (resp. \ref{rev}.2)  is a tangent vector of each disc at a point. 
(Note we use 
the same notations $\overrightarrow{x}$ (resp. $\overrightarrow{y}$) for different arrows.)
The orientation of each disc in 
Figure \ref{rev}.1 (resp. \ref{rev}.2) 
is determined by the each set $\{\overrightarrow{x},\overrightarrow{y}\}$. 
The orientation of $B\cap L_+$ (resp. $B\cap L_-$) 
coincides with that of the cylinder and that of the disc. 
We can suppose that there is a Seifert hypersurface $V$ such that 
$V\cap B$ is $P\x[0.3, 0.7]$.

$B\cap L_0$ is a disjoint union of two 2-discs and an annulus as drawn in Figure \ref{rev}.3.  
One of the 2-discs is in (the close 2-disc $P$)$\x\{0\}\x\{0\}$ and 
the other in  \newline
(the close 2-disc $P$)$\x\{1\}\x\{0\}$. 
The annulus is in $(\partial$(the close 2-disc $P$))$\x[0.4, 0.6]\x\{0\}$. 
\end{defn}

Recall that an example of (1,2)-pass-move-triples is drawn in \S1.

\begin{note}\label{hokkaido}
In the $(1,2)$-pass-move case we have the following examples: 
Let $V$ be a Seifert hypersurface for a 2-knot $K$. 
Suppose that $V$ is diffeomorphic to \newline
$((S_a^1\x S_a^2)\sharp(S_b^1\x S_b^2))-$open$B^3$. 
Take orientations of 
$S_i^1\x*$ and $*\x S_j^2$ ($i,j\in\{a,b\}$) 
so that the intersection product of 
$S_i^1\x*$ and $*\x S_j^2$ is $\delta_{ij}$. 
Suppose that the Seifert pairing of $S_a^1\x*$ and $*\x S_b^2$ is one. 
If we change the orientations of $S_b^1\x*$ and $*\x S_b^2$, 
then 
the intersection product of $S_i^1\x*$ and $*\x S_j^2$ does not change 
but the Seifert pairing of $S_a^1\x*$ and $*\x S_b^2$ cahnges $+1$ into $-1.$ 

It means the following:  
Suppose that 
we know that 
one of 2-dimensional closed oriented submanifolds, $K$ and $J$, $\subset S^4$ 
is $K_+$, and that the other $K_-$. % associated with a 4-ball $B$, 
Then we cannot distinguish $K_+$ from $K_-$ without the information 
of the orientation how $K_+$ and $K_-$ intersect $B$.  

On the other hand, in the case of the twist-move on the $(4k+1)$-dimensional submanifolds,  
we can distinguish $K_+$ from $K_-$. See Note \ref{aomori}.  
\end{note}

In \cite{Ogasa04} we introduce 
the ribbon-move for closed oriented 2-dimensional submanifolds $\subset S^4$. 
The ribbon-move is much connected with the (1,2)-pass-move  
(see  \cite[Proposition 4.2]{Ogasa04}). 
If we replace  `(1,2)-pass-move' (resp. `(1,2)-pass-move-triple') with 
 `ribbon-move' (resp. `ribbon-move-triple') in the theorems of this paper, 
similar theorems could hold. 
We draw (a part of a figure of) a ribbon-move-triple in Figure \ref{rev}.4-\ref{rev}.6.    
See \cite{Ogasa04} for the precise definition.  
%{\color{cyan}{ribbon-tripleなら(1,2)-pass-tripleかいな}}

\begin{figure}
%WinTpicVersion2.15
\unitlength 0.1in
\begin{picture}(40.72,34.25)(7.05,-34.85)
% ELLIPSE 2 0 3 0
% 4 2712 642 3104 817 3104 817 3104 817
% 
\special{pn 8}%
\special{ar 2712 242 392 175  0.0000000 6.2831853}%
% LINE 0 0 3 0
% 2 2817 3183 2817 2329
% 
\special{pn 20}%
\special{pa 2817 2783}%
\special{pa 2817 1929}%
\special{fp}%
% LINE 0 0 3 0
% 2 2593 3197 2593 2336
% 
\special{pn 20}%
\special{pa 2593 2797}%
\special{pa 2593 1936}%
\special{fp}%
% SPLINE 0 0 3 0
% 41 2600 3218 2601 3209 2605 3200 2611 3193 2620 3185 2631 3178 2643 3172 2657 3168 2673 3165 2689 3163 2705 3162 2721 3163 2737 3165 2753 3168 2767 3172 2779 3178 2790 3185 2799 3193 2805 3200 2809 3209 2810 3218 2809 3227 2805 3235 2799 3243 2790 3251 2779 3258 2767 3263 2753 3268 2737 3271 2721 3273 2705 3274 2689 3273 2673 3271 2657 3268 2643 3263 2631 3258 2620 3251 2611 3243 2605 3235 2601 3227 2600 3218
% 
\special{pn 20}%
\special{pa 2600 2818}%
\special{pa 2614 2790}%
\special{pa 2641 2773}%
\special{pa 2672 2765}%
\special{pa 2703 2762}%
\special{pa 2735 2765}%
\special{pa 2766 2772}%
\special{pa 2794 2788}%
\special{pa 2810 2815}%
\special{pa 2798 2844}%
\special{pa 2772 2861}%
\special{pa 2741 2870}%
\special{pa 2710 2874}%
\special{pa 2678 2872}%
\special{pa 2647 2864}%
\special{pa 2618 2850}%
\special{pa 2600 2824}%
\special{pa 2600 2818}%
\special{sp}%
% LINE 0 0 3 0
% 2 2593 670 2593 1524
% 
\special{pn 20}%
\special{pa 2593 270}%
\special{pa 2593 1124}%
\special{fp}%
% LINE 0 0 3 0
% 2 2817 656 2817 1517
% 
\special{pn 20}%
\special{pa 2817 256}%
\special{pa 2817 1117}%
\special{fp}%
% ELLIPSE 0 0 3 0
% 4 2705 635 2810 691 3398 1034 3398 1034
% 
\special{pn 20}%
\special{ar 2705 235 105 56  0.0000000 6.2831853}%
% ELLIPSE 2 0 3 0
% 4 2712 3274 3104 3449 3104 3449 3104 3449
% 
\special{pn 8}%
\special{ar 2712 2874 392 175  0.0000000 6.2831853}%
% LINE 2 0 3 0
% 2 3111 649 3111 3274
% 
\special{pn 8}%
\special{pa 3111 249}%
\special{pa 3111 2874}%
\special{fp}%
% LINE 2 0 3 0
% 2 2320 656 2320 3274
% 
\special{pn 8}%
\special{pa 2320 256}%
\special{pa 2320 2874}%
\special{fp}%
% ELLIPSE 2 0 3 0
% 4 1242 649 1634 824 1634 824 1634 824
% 
\special{pn 8}%
\special{ar 1242 249 392 175  0.0000000 6.2831853}%
% ELLIPSE 2 0 3 0
% 4 1242 3281 1634 3456 1634 3456 1634 3456
% 
\special{pn 8}%
\special{ar 1242 2881 392 175  0.0000000 6.2831853}%
% LINE 2 0 3 0
% 2 1641 656 1641 3281
% 
\special{pn 8}%
\special{pa 1641 256}%
\special{pa 1641 2881}%
\special{fp}%
% LINE 2 0 3 0
% 2 850 663 850 3281
% 
\special{pn 8}%
\special{pa 850 263}%
\special{pa 850 2881}%
\special{fp}%
% ELLIPSE 2 0 3 0
% 4 4378 635 4770 810 4770 810 4770 810
% 
\special{pn 8}%
\special{ar 4378 235 392 175  0.0000000 6.2831853}%
% ELLIPSE 2 0 3 0
% 4 4378 3267 4770 3442 4770 3442 4770 3442
% 
\special{pn 8}%
\special{ar 4378 2867 392 175  0.0000000 6.2831853}%
% LINE 2 0 3 0
% 2 4777 642 4777 3267
% 
\special{pn 8}%
\special{pa 4777 242}%
\special{pa 4777 2867}%
\special{fp}%
% LINE 2 0 3 0
% 2 3986 649 3986 3267
% 
\special{pn 8}%
\special{pa 3986 249}%
\special{pa 3986 2867}%
\special{fp}%
% ELLIPSE 0 0 3 0
% 4 2705 1531 2810 1587 3398 1930 3398 1930
% 
\special{pn 20}%
\special{ar 2705 1131 105 56  0.0000000 6.2831853}%
% ELLIPSE 0 0 3 0
% 4 2705 2315 2810 2371 3398 2714 3398 2714
% 
\special{pn 20}%
\special{ar 2705 1915 105 56  0.0000000 6.2831853}%
% ELLIPSE 0 0 3 0
% 4 2719 1867 3111 2042 3111 2042 3111 2042
% 
\special{pn 20}%
\special{ar 2719 1467 392 175  0.0000000 6.2831853}%
% LINE 2 2 3 0
% 2 2817 1538 4266 1538
% 
\special{pn 8}%
\special{pa 2817 1138}%
\special{pa 4266 1138}%
\special{dt 0.045}%
\special{pa 4266 1138}%
\special{pa 4265 1138}%
\special{dt 0.045}%
% LINE 2 2 3 0
% 2 2831 2322 4280 2322
% 
\special{pn 8}%
\special{pa 2831 1922}%
\special{pa 4280 1922}%
\special{dt 0.045}%
\special{pa 4280 1922}%
\special{pa 4279 1922}%
\special{dt 0.045}%
% ELLIPSE 0 0 3 0
% 4 4350 1538 4455 1594 5043 1937 5043 1937
% 
\special{pn 20}%
\special{ar 4350 1138 105 56  0.0000000 6.2831853}%
% ELLIPSE 0 0 3 0
% 4 4350 2329 4455 2385 5043 2728 5043 2728
% 
\special{pn 20}%
\special{ar 4350 1929 105 56  0.0000000 6.2831853}%
% LINE 0 0 3 0
% 2 4238 1545 4238 2308
% 
\special{pn 20}%
\special{pa 4238 1145}%
\special{pa 4238 1908}%
\special{fp}%
% LINE 0 0 3 0
% 2 4469 1566 4469 2329
% 
\special{pn 20}%
\special{pa 4469 1166}%
\special{pa 4469 1929}%
\special{fp}%
% STR 2 0 3 0
% 3 1095 3568 1095 3638 2 0
% t=-0.5
\put(10.9500,-32.3800){\makebox(0,0)[lb]{t=-0.5}}%
% STR 2 0 3 0
% 3 2495 3568 2495 3638 2 0
% t=0
\put(24.9500,-32.3800){\makebox(0,0)[lb]{t=0}}%
% STR 2 0 3 0
% 3 4175 3568 4175 3638 2 0
% t=0.5
\put(41.7500,-32.3800){\makebox(0,0)[lb]{t=0.5}}%
% ELLIPSE 2 0 3 0
% 4 2712 642 3104 817 3104 817 3104 817
% 
\special{pn 8}%
\special{ar 2712 242 392 175  0.0000000 6.2831853}%
% LINE 2 0 3 1
% 2 2320 656 2320 3274
% 
\special{pn 8}%
\special{pa 2320 256}%
\special{pa 2320 2874}%
\special{fp}%
% LINE 2 0 3 2
% 2 3111 649 3111 3274
% 
\special{pn 8}%
\special{pa 3111 249}%
\special{pa 3111 2874}%
\special{fp}%
% ELLIPSE 2 0 3 3
% 4 2712 3274 3104 3449 3104 3449 3104 3449
% 
\special{pn 8}%
\special{ar 2712 2874 392 175  0.0000000 6.2831853}%
% ELLIPSE 0 0 3 4
% 4 2705 635 2810 691 3398 1034 3398 1034
% 
\special{pn 20}%
\special{ar 2705 235 105 56  0.0000000 6.2831853}%
% LINE 0 0 3 5
% 2 2817 656 2817 1517
% 
\special{pn 20}%
\special{pa 2817 256}%
\special{pa 2817 1117}%
\special{fp}%
% LINE 0 0 3 6
% 2 2593 670 2593 1524
% 
\special{pn 20}%
\special{pa 2593 270}%
\special{pa 2593 1124}%
\special{fp}%
% SPLINE 0 0 3 7
% 41 2600 3218 2601 3209 2605 3200 2611 3193 2620 3185 2631 3178 2643 3172 2657 3168 2673 3165 2689 3163 2705 3162 2721 3163 2737 3165 2753 3168 2767 3172 2779 3178 2790 3185 2799 3193 2805 3200 2809 3209 2810 3218 2809 3227 2805 3235 2799 3243 2790 3251 2779 3258 2767 3263 2753 3268 2737 3271 2721 3273 2705 3274 2689 3273 2673 3271 2657 3268 2643 3263 2631 3258 2620 3251 2611 3243 2605 3235 2601 3227 2600 3218
% 
\special{pn 20}%
\special{pa 2600 2818}%
\special{pa 2614 2790}%
\special{pa 2641 2773}%
\special{pa 2672 2765}%
\special{pa 2703 2762}%
\special{pa 2735 2765}%
\special{pa 2766 2772}%
\special{pa 2794 2788}%
\special{pa 2810 2815}%
\special{pa 2798 2844}%
\special{pa 2772 2861}%
\special{pa 2741 2870}%
\special{pa 2710 2874}%
\special{pa 2678 2872}%
\special{pa 2647 2864}%
\special{pa 2618 2850}%
\special{pa 2600 2824}%
\special{pa 2600 2818}%
\special{sp}%
% LINE 0 0 3 8
% 2 2593 3197 2593 2336
% 
\special{pn 20}%
\special{pa 2593 2797}%
\special{pa 2593 1936}%
\special{fp}%
% LINE 0 0 3 9
% 2 2817 3183 2817 2329
% 
\special{pn 20}%
\special{pa 2817 2783}%
\special{pa 2817 1929}%
\special{fp}%
% ELLIPSE 2 0 3 0
% 4 1242 649 1634 824 1634 824 1634 824
% 
\special{pn 8}%
\special{ar 1242 249 392 175  0.0000000 6.2831853}%
% ELLIPSE 2 0 3 0
% 4 1242 3281 1634 3456 1634 3456 1634 3456
% 
\special{pn 8}%
\special{ar 1242 2881 392 175  0.0000000 6.2831853}%
% LINE 2 0 3 0
% 2 1641 656 1641 3281
% 
\special{pn 8}%
\special{pa 1641 256}%
\special{pa 1641 2881}%
\special{fp}%
% LINE 2 0 3 0
% 2 850 663 850 3281
% 
\special{pn 8}%
\special{pa 850 263}%
\special{pa 850 2881}%
\special{fp}%
% ELLIPSE 2 0 3 0
% 4 4378 635 4770 810 4770 810 4770 810
% 
\special{pn 8}%
\special{ar 4378 235 392 175  0.0000000 6.2831853}%
% ELLIPSE 2 0 3 0
% 4 4378 3267 4770 3442 4770 3442 4770 3442
% 
\special{pn 8}%
\special{ar 4378 2867 392 175  0.0000000 6.2831853}%
% LINE 2 0 3 0
% 2 4777 642 4777 3267
% 
\special{pn 8}%
\special{pa 4777 242}%
\special{pa 4777 2867}%
\special{fp}%
% LINE 2 0 3 0
% 2 3986 649 3986 3267
% 
\special{pn 8}%
\special{pa 3986 249}%
\special{pa 3986 2867}%
\special{fp}%
% STR 2 0 3 0
% 3 3000 3870 3000 3970 5 0
% Figure \ref{rev}.4: $L_+$ of a ribbon-move-triple. 
\put(30.0000,-35.7000){\makebox(0,0){Figure \ref{rev}.4: $L_+$ of a ribbon-move-triple. }}%
\end{picture}%
\bigbreak
\bigbreak
%\bigbreak
\end{figure}

\begin{figure}
%WinTpicVersion2.15
\unitlength 0.1in
\begin{picture}(39.27,33.05)(8.50,-33.65)
% ELLIPSE 2 0 3 0
% 4 2712 642 3104 817 3104 817 3104 817
% 
\special{pn 8}%
\special{ar 2712 242 392 175  0.0000000 6.2831853}%
% LINE 0 0 3 0
% 2 2817 3183 2817 2329
% 
\special{pn 20}%
\special{pa 2817 2783}%
\special{pa 2817 1929}%
\special{fp}%
% LINE 0 0 3 0
% 2 2593 3197 2593 2336
% 
\special{pn 20}%
\special{pa 2593 2797}%
\special{pa 2593 1936}%
\special{fp}%
% SPLINE 0 0 3 0
% 41 2600 3218 2601 3209 2605 3200 2611 3193 2620 3185 2631 3178 2643 3172 2657 3168 2673 3165 2689 3163 2705 3162 2721 3163 2737 3165 2753 3168 2767 3172 2779 3178 2790 3185 2799 3193 2805 3200 2809 3209 2810 3218 2809 3227 2805 3235 2799 3243 2790 3251 2779 3258 2767 3263 2753 3268 2737 3271 2721 3273 2705 3274 2689 3273 2673 3271 2657 3268 2643 3263 2631 3258 2620 3251 2611 3243 2605 3235 2601 3227 2600 3218
% 
\special{pn 20}%
\special{pa 2600 2818}%
\special{pa 2614 2790}%
\special{pa 2641 2773}%
\special{pa 2672 2765}%
\special{pa 2703 2762}%
\special{pa 2735 2765}%
\special{pa 2766 2772}%
\special{pa 2794 2788}%
\special{pa 2810 2815}%
\special{pa 2798 2844}%
\special{pa 2772 2861}%
\special{pa 2741 2870}%
\special{pa 2710 2874}%
\special{pa 2678 2872}%
\special{pa 2647 2864}%
\special{pa 2618 2850}%
\special{pa 2600 2824}%
\special{pa 2600 2818}%
\special{sp}%
% LINE 0 0 3 0
% 2 2593 670 2593 1524
% 
\special{pn 20}%
\special{pa 2593 270}%
\special{pa 2593 1124}%
\special{fp}%
% LINE 0 0 3 0
% 2 2817 656 2817 1517
% 
\special{pn 20}%
\special{pa 2817 256}%
\special{pa 2817 1117}%
\special{fp}%
% ELLIPSE 0 0 3 0
% 4 2705 635 2810 691 3398 1034 3398 1034
% 
\special{pn 20}%
\special{ar 2705 235 105 56  0.0000000 6.2831853}%
% ELLIPSE 2 0 3 0
% 4 2712 3274 3104 3449 3104 3449 3104 3449
% 
\special{pn 8}%
\special{ar 2712 2874 392 175  0.0000000 6.2831853}%
% LINE 2 0 3 0
% 2 3111 649 3111 3274
% 
\special{pn 8}%
\special{pa 3111 249}%
\special{pa 3111 2874}%
\special{fp}%
% LINE 2 0 3 0
% 2 2320 656 2320 3274
% 
\special{pn 8}%
\special{pa 2320 256}%
\special{pa 2320 2874}%
\special{fp}%
% ELLIPSE 2 0 3 0
% 4 1242 649 1634 824 1634 824 1634 824
% 
\special{pn 8}%
\special{ar 1242 249 392 175  0.0000000 6.2831853}%
% ELLIPSE 2 0 3 0
% 4 1242 3281 1634 3456 1634 3456 1634 3456
% 
\special{pn 8}%
\special{ar 1242 2881 392 175  0.0000000 6.2831853}%
% LINE 2 0 3 0
% 2 1641 656 1641 3281
% 
\special{pn 8}%
\special{pa 1641 256}%
\special{pa 1641 2881}%
\special{fp}%
% LINE 2 0 3 0
% 2 850 663 850 3281
% 
\special{pn 8}%
\special{pa 850 263}%
\special{pa 850 2881}%
\special{fp}%
% ELLIPSE 2 0 3 0
% 4 4378 635 4770 810 4770 810 4770 810
% 
\special{pn 8}%
\special{ar 4378 235 392 175  0.0000000 6.2831853}%
% ELLIPSE 2 0 3 0
% 4 4378 3267 4770 3442 4770 3442 4770 3442
% 
\special{pn 8}%
\special{ar 4378 2867 392 175  0.0000000 6.2831853}%
% LINE 2 0 3 0
% 2 4777 642 4777 3267
% 
\special{pn 8}%
\special{pa 4777 242}%
\special{pa 4777 2867}%
\special{fp}%
% LINE 2 0 3 0
% 2 3986 649 3986 3267
% 
\special{pn 8}%
\special{pa 3986 249}%
\special{pa 3986 2867}%
\special{fp}%
% ELLIPSE 0 0 3 0
% 4 2705 1531 2810 1587 3398 1930 3398 1930
% 
\special{pn 20}%
\special{ar 2705 1131 105 56  0.0000000 6.2831853}%
% ELLIPSE 0 0 3 0
% 4 2705 2315 2810 2371 3398 2714 3398 2714
% 
\special{pn 20}%
\special{ar 2705 1915 105 56  0.0000000 6.2831853}%
% ELLIPSE 0 0 3 0
% 4 2719 1867 3111 2042 3111 2042 3111 2042
% 
\special{pn 20}%
\special{ar 2719 1467 392 175  0.0000000 6.2831853}%
% STR 2 0 3 0
% 3 1095 3568 1095 3638 2 0
% t=-0.5
\put(10.9500,-32.3800){\makebox(0,0)[lb]{t=-0.5}}%
% STR 2 0 3 0
% 3 2495 3568 2495 3638 2 0
% t=0
\put(24.9500,-32.3800){\makebox(0,0)[lb]{t=0}}%
% STR 2 0 3 0
% 3 4175 3568 4175 3638 2 0
% t=0.5
\put(41.7500,-32.3800){\makebox(0,0)[lb]{t=0.5}}%
% ELLIPSE 2 0 3 0
% 4 2712 642 3104 817 3104 817 3104 817
% 
\special{pn 8}%
\special{ar 2712 242 392 175  0.0000000 6.2831853}%
% LINE 2 0 3 1
% 2 2320 656 2320 3274
% 
\special{pn 8}%
\special{pa 2320 256}%
\special{pa 2320 2874}%
\special{fp}%
% LINE 2 0 3 2
% 2 3111 649 3111 3274
% 
\special{pn 8}%
\special{pa 3111 249}%
\special{pa 3111 2874}%
\special{fp}%
% ELLIPSE 2 0 3 3
% 4 2712 3274 3104 3449 3104 3449 3104 3449
% 
\special{pn 8}%
\special{ar 2712 2874 392 175  0.0000000 6.2831853}%
% ELLIPSE 0 0 3 4
% 4 2705 635 2810 691 3398 1034 3398 1034
% 
\special{pn 20}%
\special{ar 2705 235 105 56  0.0000000 6.2831853}%
% LINE 0 0 3 5
% 2 2817 656 2817 1517
% 
\special{pn 20}%
\special{pa 2817 256}%
\special{pa 2817 1117}%
\special{fp}%
% LINE 0 0 3 6
% 2 2593 670 2593 1524
% 
\special{pn 20}%
\special{pa 2593 270}%
\special{pa 2593 1124}%
\special{fp}%
% SPLINE 0 0 3 7
% 41 2600 3218 2601 3209 2605 3200 2611 3193 2620 3185 2631 3178 2643 3172 2657 3168 2673 3165 2689 3163 2705 3162 2721 3163 2737 3165 2753 3168 2767 3172 2779 3178 2790 3185 2799 3193 2805 3200 2809 3209 2810 3218 2809 3227 2805 3235 2799 3243 2790 3251 2779 3258 2767 3263 2753 3268 2737 3271 2721 3273 2705 3274 2689 3273 2673 3271 2657 3268 2643 3263 2631 3258 2620 3251 2611 3243 2605 3235 2601 3227 2600 3218
% 
\special{pn 20}%
\special{pa 2600 2818}%
\special{pa 2614 2790}%
\special{pa 2641 2773}%
\special{pa 2672 2765}%
\special{pa 2703 2762}%
\special{pa 2735 2765}%
\special{pa 2766 2772}%
\special{pa 2794 2788}%
\special{pa 2810 2815}%
\special{pa 2798 2844}%
\special{pa 2772 2861}%
\special{pa 2741 2870}%
\special{pa 2710 2874}%
\special{pa 2678 2872}%
\special{pa 2647 2864}%
\special{pa 2618 2850}%
\special{pa 2600 2824}%
\special{pa 2600 2818}%
\special{sp}%
% LINE 0 0 3 8
% 2 2593 3197 2593 2336
% 
\special{pn 20}%
\special{pa 2593 2797}%
\special{pa 2593 1936}%
\special{fp}%
% LINE 0 0 3 9
% 2 2817 3183 2817 2329
% 
\special{pn 20}%
\special{pa 2817 2783}%
\special{pa 2817 1929}%
\special{fp}%
% ELLIPSE 2 0 3 0
% 4 1242 649 1634 824 1634 824 1634 824
% 
\special{pn 8}%
\special{ar 1242 249 392 175  0.0000000 6.2831853}%
% ELLIPSE 2 0 3 0
% 4 1242 3281 1634 3456 1634 3456 1634 3456
% 
\special{pn 8}%
\special{ar 1242 2881 392 175  0.0000000 6.2831853}%
% LINE 2 0 3 0
% 2 1641 656 1641 3281
% 
\special{pn 8}%
\special{pa 1641 256}%
\special{pa 1641 2881}%
\special{fp}%
% LINE 2 0 3 0
% 2 850 663 850 3281
% 
\special{pn 8}%
\special{pa 850 263}%
\special{pa 850 2881}%
\special{fp}%
% ELLIPSE 2 0 3 0
% 4 4378 635 4770 810 4770 810 4770 810
% 
\special{pn 8}%
\special{ar 4378 235 392 175  0.0000000 6.2831853}%
% ELLIPSE 2 0 3 0
% 4 4378 3267 4770 3442 4770 3442 4770 3442
% 
\special{pn 8}%
\special{ar 4378 2867 392 175  0.0000000 6.2831853}%
% LINE 2 0 3 0
% 2 4777 642 4777 3267
% 
\special{pn 8}%
\special{pa 4777 242}%
\special{pa 4777 2867}%
\special{fp}%
% LINE 2 0 3 0
% 2 3986 649 3986 3267
% 
\special{pn 8}%
\special{pa 3986 249}%
\special{pa 3986 2867}%
\special{fp}%
% LINE 2 2 3 0
% 2 2579 1531 1368 1531
% 
\special{pn 8}%
\special{pa 2579 1131}%
\special{pa 1368 1131}%
\special{dt 0.045}%
\special{pa 1368 1131}%
\special{pa 1369 1131}%
\special{dt 0.045}%
% LINE 2 2 3 0
% 2 2586 2336 1389 2336
% 
\special{pn 8}%
\special{pa 2586 1936}%
\special{pa 1389 1936}%
\special{dt 0.045}%
\special{pa 1389 1936}%
\special{pa 1390 1936}%
\special{dt 0.045}%
% ELLIPSE 0 0 3 0
% 4 1249 1531 1354 1587 1942 1930 1942 1930
% 
\special{pn 20}%
\special{ar 1249 1131 105 56  0.0000000 6.2831853}%
% ELLIPSE 0 0 3 1
% 4 1249 2322 1354 2378 1942 2721 1942 2721
% 
\special{pn 20}%
\special{ar 1249 1922 105 56  0.0000000 6.2831853}%
% LINE 0 0 3 2
% 2 1137 1538 1137 2301
% 
\special{pn 20}%
\special{pa 1137 1138}%
\special{pa 1137 1901}%
\special{fp}%
% LINE 0 0 3 3
% 2 1368 1559 1368 2322
% 
\special{pn 20}%
\special{pa 1368 1159}%
\special{pa 1368 1922}%
\special{fp}%
% STR 2 0 3 0
% 3 3160 3750 3160 3850 5 0
% Figure \ref{rev}.5: $L_-$ of a ribbon-move-triple. 
\put(31.6000,-34.5000){\makebox(0,0){Figure \ref{rev}.5: $L_-$ of a ribbon-move-triple. }}%
\end{picture}%
\bigbreak
\bigbreak
\end{figure}

\begin{figure}
%WinTpicVersion2.15
\unitlength 0.1in
\begin{picture}(39.52,33.55)(8.25,-34.15)
% ELLIPSE 2 0 3 0
% 4 2712 642 3104 817 3104 817 3104 817
% 
\special{pn 8}%
\special{ar 2712 242 392 175  0.0000000 6.2831853}%
% SPLINE 0 0 3 0
% 41 2600 3218 2601 3209 2605 3200 2611 3193 2620 3185 2631 3178 2643 3172 2657 3168 2673 3165 2689 3163 2705 3162 2721 3163 2737 3165 2753 3168 2767 3172 2779 3178 2790 3185 2799 3193 2805 3200 2809 3209 2810 3218 2809 3227 2805 3235 2799 3243 2790 3251 2779 3258 2767 3263 2753 3268 2737 3271 2721 3273 2705 3274 2689 3273 2673 3271 2657 3268 2643 3263 2631 3258 2620 3251 2611 3243 2605 3235 2601 3227 2600 3218
% 
\special{pn 20}%
\special{pa 2600 2818}%
\special{pa 2614 2790}%
\special{pa 2641 2773}%
\special{pa 2672 2765}%
\special{pa 2703 2762}%
\special{pa 2735 2765}%
\special{pa 2766 2772}%
\special{pa 2794 2788}%
\special{pa 2810 2815}%
\special{pa 2798 2844}%
\special{pa 2772 2861}%
\special{pa 2741 2870}%
\special{pa 2710 2874}%
\special{pa 2678 2872}%
\special{pa 2647 2864}%
\special{pa 2618 2850}%
\special{pa 2600 2824}%
\special{pa 2600 2818}%
\special{sp}%
% ELLIPSE 0 0 0 0
% 4 2705 635 2810 691 3398 1034 3398 1034
% 
\special{pn 20}%
\special{sh 0.600}%
\special{ar 2705 235 105 56  0.0000000 6.2831853}%
% ELLIPSE 0 0 3 0
% 4 2712 3274 3104 3449 3104 3449 3104 3449
% 
\special{pn 20}%
\special{ar 2712 2874 392 175  0.0000000 6.2831853}%
% LINE 2 0 3 0
% 2 3111 649 3111 3274
% 
\special{pn 8}%
\special{pa 3111 249}%
\special{pa 3111 2874}%
\special{fp}%
% LINE 2 0 3 0
% 2 2320 656 2320 3274
% 
\special{pn 8}%
\special{pa 2320 256}%
\special{pa 2320 2874}%
\special{fp}%
% ELLIPSE 2 0 3 0
% 4 1242 649 1634 824 1634 824 1634 824
% 
\special{pn 8}%
\special{ar 1242 249 392 175  0.0000000 6.2831853}%
% ELLIPSE 2 0 3 0
% 4 1242 3281 1634 3456 1634 3456 1634 3456
% 
\special{pn 8}%
\special{ar 1242 2881 392 175  0.0000000 6.2831853}%
% LINE 2 0 3 0
% 2 1641 656 1641 3281
% 
\special{pn 8}%
\special{pa 1641 256}%
\special{pa 1641 2881}%
\special{fp}%
% LINE 2 0 3 0
% 2 850 663 850 3281
% 
\special{pn 8}%
\special{pa 850 263}%
\special{pa 850 2881}%
\special{fp}%
% ELLIPSE 2 0 3 0
% 4 4378 635 4770 810 4770 810 4770 810
% 
\special{pn 8}%
\special{ar 4378 235 392 175  0.0000000 6.2831853}%
% ELLIPSE 2 0 3 0
% 4 4378 3267 4770 3442 4770 3442 4770 3442
% 
\special{pn 8}%
\special{ar 4378 2867 392 175  0.0000000 6.2831853}%
% LINE 2 0 3 0
% 2 4777 642 4777 3267
% 
\special{pn 8}%
\special{pa 4777 242}%
\special{pa 4777 2867}%
\special{fp}%
% LINE 2 0 3 0
% 2 3986 649 3986 3267
% 
\special{pn 8}%
\special{pa 3986 249}%
\special{pa 3986 2867}%
\special{fp}%
% ELLIPSE 0 0 3 0
% 4 2719 1867 3111 2042 3111 2042 3111 2042
% 
\special{pn 20}%
\special{ar 2719 1467 392 175  0.0000000 6.2831853}%
% STR 2 0 3 0
% 3 1095 3568 1095 3638 2 0
% t=-0.5
\put(10.9500,-32.3800){\makebox(0,0)[lb]{t=-0.5}}%
% STR 2 0 3 0
% 3 2495 3568 2495 3638 2 0
% t=0
\put(24.9500,-32.3800){\makebox(0,0)[lb]{t=0}}%
% STR 2 0 3 0
% 3 4175 3568 4175 3638 2 0
% t=0.5
\put(41.7500,-32.3800){\makebox(0,0)[lb]{t=0.5}}%
% ELLIPSE 2 0 3 0
% 4 2712 642 3104 817 3104 817 3104 817
% 
\special{pn 8}%
\special{ar 2712 242 392 175  0.0000000 6.2831853}%
% SPLINE 0 0 3 0
% 41 2600 3218 2601 3209 2605 3200 2611 3193 2620 3185 2631 3178 2643 3172 2657 3168 2673 3165 2689 3163 2705 3162 2721 3163 2737 3165 2753 3168 2767 3172 2779 3178 2790 3185 2799 3193 2805 3200 2809 3209 2810 3218 2809 3227 2805 3235 2799 3243 2790 3251 2779 3258 2767 3263 2753 3268 2737 3271 2721 3273 2705 3274 2689 3273 2673 3271 2657 3268 2643 3263 2631 3258 2620 3251 2611 3243 2605 3235 2601 3227 2600 3218
% 
\special{pn 20}%
\special{pa 2600 2818}%
\special{pa 2614 2790}%
\special{pa 2641 2773}%
\special{pa 2672 2765}%
\special{pa 2703 2762}%
\special{pa 2735 2765}%
\special{pa 2766 2772}%
\special{pa 2794 2788}%
\special{pa 2810 2815}%
\special{pa 2798 2844}%
\special{pa 2772 2861}%
\special{pa 2741 2870}%
\special{pa 2710 2874}%
\special{pa 2678 2872}%
\special{pa 2647 2864}%
\special{pa 2618 2850}%
\special{pa 2600 2824}%
\special{pa 2600 2818}%
\special{sp}%
% ELLIPSE 0 0 0 0
% 4 2705 635 2810 691 3398 1034 3398 1034
% 
\special{pn 20}%
\special{sh 0.600}%
\special{ar 2705 235 105 56  0.0000000 6.2831853}%
% ELLIPSE 2 0 3 0
% 4 2712 3274 3104 3449 3104 3449 3104 3449
% 
\special{pn 8}%
\special{ar 2712 2874 392 175  0.0000000 6.2831853}%
% LINE 2 0 3 0
% 2 3111 649 3111 3274
% 
\special{pn 8}%
\special{pa 3111 249}%
\special{pa 3111 2874}%
\special{fp}%
% LINE 2 0 3 0
% 2 2320 656 2320 3274
% 
\special{pn 8}%
\special{pa 2320 256}%
\special{pa 2320 2874}%
\special{fp}%
% ELLIPSE 2 0 3 0
% 4 1242 649 1634 824 1634 824 1634 824
% 
\special{pn 8}%
\special{ar 1242 249 392 175  0.0000000 6.2831853}%
% ELLIPSE 2 0 3 0
% 4 1242 3281 1634 3456 1634 3456 1634 3456
% 
\special{pn 8}%
\special{ar 1242 2881 392 175  0.0000000 6.2831853}%
% LINE 2 0 3 0
% 2 1641 656 1641 3281
% 
\special{pn 8}%
\special{pa 1641 256}%
\special{pa 1641 2881}%
\special{fp}%
% LINE 2 0 3 0
% 2 850 663 850 3281
% 
\special{pn 8}%
\special{pa 850 263}%
\special{pa 850 2881}%
\special{fp}%
% ELLIPSE 2 0 3 0
% 4 4378 635 4770 810 4770 810 4770 810
% 
\special{pn 8}%
\special{ar 4378 235 392 175  0.0000000 6.2831853}%
% ELLIPSE 2 0 3 0
% 4 4378 3267 4770 3442 4770 3442 4770 3442
% 
\special{pn 8}%
\special{ar 4378 2867 392 175  0.0000000 6.2831853}%
% LINE 2 0 3 0
% 2 4777 642 4777 3267
% 
\special{pn 8}%
\special{pa 4777 242}%
\special{pa 4777 2867}%
\special{fp}%
% LINE 2 0 3 0
% 2 3986 649 3986 3267
% 
\special{pn 8}%
\special{pa 3986 249}%
\special{pa 3986 2867}%
\special{fp}%
% LINE 0 0 3 0
% 2 2390 1972 2390 3358
% 
\special{pn 20}%
\special{pa 2390 1572}%
\special{pa 2390 2958}%
\special{fp}%
% LINE 0 0 3 0
% 2 2495 1755 2495 1888
% 
\special{pn 20}%
\special{pa 2495 1355}%
\special{pa 2495 1488}%
\special{fp}%
% LINE 0 0 3 0
% 2 2502 2336 2502 2224
% 
\special{pn 20}%
\special{pa 2502 1936}%
\special{pa 2502 1824}%
\special{fp}%
% LINE 0 0 3 0
% 2 2502 2581 2502 2462
% 
\special{pn 20}%
\special{pa 2502 2181}%
\special{pa 2502 2062}%
\special{fp}%
% LINE 0 0 3 0
% 2 2502 2812 2502 2742
% 
\special{pn 20}%
\special{pa 2502 2412}%
\special{pa 2502 2342}%
\special{fp}%
% LINE 0 0 3 0
% 2 2502 3134 2509 2931
% 
\special{pn 20}%
\special{pa 2502 2734}%
\special{pa 2509 2531}%
\special{fp}%
% LINE 0 0 3 0
% 2 2537 3141 2502 3134
% 
\special{pn 20}%
\special{pa 2537 2741}%
\special{pa 2502 2734}%
\special{fp}%
% LINE 0 0 3 0
% 2 2593 3183 2558 3169
% 
\special{pn 20}%
\special{pa 2593 2783}%
\special{pa 2558 2769}%
\special{fp}%
% LINE 0 0 3 0
% 2 3013 3309 3048 3344
% 
\special{pn 20}%
\special{pa 3013 2909}%
\special{pa 3048 2944}%
\special{fp}%
% LINE 0 0 3 0
% 2 2831 3218 3013 3309
% 
\special{pn 20}%
\special{pa 2831 2818}%
\special{pa 3013 2909}%
\special{fp}%
% LINE 0 0 3 0
% 2 3055 3351 3055 1937
% 
\special{pn 20}%
\special{pa 3055 2951}%
\special{pa 3055 1537}%
\special{fp}%
% LINE 0 0 3 0
% 2 2901 3393 2775 3267
% 
\special{pn 20}%
\special{pa 2901 2993}%
\special{pa 2775 2867}%
\special{fp}%
% LINE 0 0 3 0
% 2 2894 2028 2901 3393
% 
\special{pn 20}%
\special{pa 2894 1628}%
\special{pa 2901 2993}%
\special{fp}%
% LINE 0 0 3 0
% 2 2614 3421 2670 3281
% 
\special{pn 20}%
\special{pa 2614 3021}%
\special{pa 2670 2881}%
\special{fp}%
% LINE 0 0 3 0
% 2 2635 2042 2614 3421
% 
\special{pn 20}%
\special{pa 2635 1642}%
\special{pa 2614 3021}%
\special{fp}%
% LINE 0 0 3 0
% 2 2390 3344 2600 3260
% 
\special{pn 20}%
\special{pa 2390 2944}%
\special{pa 2600 2860}%
\special{fp}%
% LINE 0 0 3 0
% 2 2957 3162 2880 3169
% 
\special{pn 20}%
\special{pa 2957 2762}%
\special{pa 2880 2769}%
\special{fp}%
% LINE 0 0 3 0
% 2 2957 2973 2957 3162
% 
\special{pn 20}%
\special{pa 2957 2573}%
\special{pa 2957 2762}%
\special{fp}%
% LINE 0 0 3 0
% 2 2950 2742 2950 2889
% 
\special{pn 20}%
\special{pa 2950 2342}%
\special{pa 2950 2489}%
\special{fp}%
% LINE 0 0 3 0
% 2 2950 2357 2943 2560
% 
\special{pn 20}%
\special{pa 2950 1957}%
\special{pa 2943 2160}%
\special{fp}%
% LINE 0 0 3 0
% 2 2950 2077 2950 2259
% 
\special{pn 20}%
\special{pa 2950 1677}%
\special{pa 2950 1859}%
\special{fp}%
% LINE 0 0 3 0
% 2 2943 1741 2943 1874
% 
\special{pn 20}%
\special{pa 2943 1341}%
\special{pa 2943 1474}%
\special{fp}%
% LINE 0 0 3 0
% 2 2726 3113 2719 3148
% 
\special{pn 20}%
\special{pa 2726 2713}%
\special{pa 2719 2748}%
\special{fp}%
% LINE 0 0 3 0
% 2 2712 2889 2712 3092
% 
\special{pn 20}%
\special{pa 2712 2489}%
\special{pa 2712 2692}%
\special{fp}%
% LINE 0 0 3 0
% 2 2712 2616 2712 2735
% 
\special{pn 20}%
\special{pa 2712 2216}%
\special{pa 2712 2335}%
\special{fp}%
% LINE 0 0 3 0
% 2 2712 2364 2712 2518
% 
\special{pn 20}%
\special{pa 2712 1964}%
\special{pa 2712 2118}%
\special{fp}%
% LINE 0 0 3 0
% 2 2712 2091 2712 2266
% 
\special{pn 20}%
\special{pa 2712 1691}%
\special{pa 2712 1866}%
\special{fp}%
% LINE 0 0 3 0
% 2 2712 1944 2712 1986
% 
\special{pn 20}%
\special{pa 2712 1544}%
\special{pa 2712 1586}%
\special{fp}%
% LINE 0 0 3 0
% 2 2712 1720 2712 1846
% 
\special{pn 20}%
\special{pa 2712 1320}%
\special{pa 2712 1446}%
\special{fp}%
% STR 2 0 3 0
% 3 3120 3800 3120 3900 5 0
% Figure \ref{rev}.6: $L_0$ of a ribbon-move-triple. 
\put(31.2000,-35.0000){\makebox(0,0){Figure \ref{rev}.6: $L_0$ of a ribbon-move-triple. }}%
\end{picture}%
\bigbreak
\bigbreak
\end{figure}

%%%%%%%%%%%%%%%%%%%%%%%%%%%%%%%%%%%%%%%%
\bigbreak
\section{Review of the $\Q[t,t^{-1}]$-Alexander polynomial for 
(not necessarily connected) $n$-dimensional closed oriented submanifolds 
in $S^{n+2}$ $(n\geqq2)$}\label{revale}

\noindent 
In this section we review 
the $\Q[t,t^{-1}]$-%(resp. $\Z[t,t^{-1}]$-)
Alexander polynomial for 
(not necessarily connected) $n$-dimensional closed oriented submanifolds 
 $\subset S^{n+2}$ $(n\geqq2)$.   
In \S\ref{m2} we define the $\Z[t,t^{-1}]$-Alexander polynomial 
for 2-dimensional closed oriented submanifolds $\subset S^4$. 
In \S\ref{mhigh} we define 
the `normalized' Alexander polynomial 
for a kind of $(4k+1)$-submanifolds $\subset S^{4k+3}$.  
Of course these invariants are connected each other. 

%{\color{cyan}{
%other cases もwithout difficulty で論じられる、と一言入れとく
%これは§６のラストに入れた
%}}

\bigbreak  
%We work in the smooth category. 
Let $K=(K_1,...,K_\xi)$ be an $n$-dimensional closed oriented submanifold of 
$S^{n+2}$ ($n\in\N$). 
Let each $K_i$ be connected. 
If $K_i$ is PL homeomorphic to the standard sphere,  
$K_i$ is called an 
{\it $n$-dimensional $($spherical$)$ knot.}  
If each $K_i$ is an $n$-knot, 
$K$ is called an 
{\it $\xi$-component $n$-dimensional $($spherical$)$ link.}

It is known that 
the tubular neighborhood of $K$ is 
diffeomorphic to 
$K\x D^2$ (see \cite[pages 49 and 50]{Kirby}). 
Let $X=\overline{S^{n+2}-(K\x D^2)}$. 
By using the orientation of $S^{n+2}$ and that of $K$,  
we can determine 
an orientation of $\partial D^2$. Take a homomorphism
$\alpha:H_1(X;\Z)\rightarrow\Z$ to carry 
all $[\partial D^2]$ with the orientations to $+1$.
Take the infinite cyclic covering 
$\pi:$ $\widetilde X\rightarrow X$ 
associated with $\alpha$. 
$\widetilde X$ is called 
the {\it canonical cyclic covering space} of $K$. 
We can regard 
$H_p(\widetilde X;\Z)$ 
(resp. $H_p(\widetilde X;\Q)$)   
as a 
$\Z[t, t^{-1}]$-module 
(resp. $\Q[t, t^{-1}]$-module)   
by using 
the covering translation 
$\widetilde X\rightarrow\widetilde X$ 
defined by $\alpha$. 
It is called the {\it $\Z[t, t^{-1}]$-$p$-Alexander module}  
(resp. {\it $\Q[t, t^{-1}]$-$p$-Alexander module}).

\begin{defn}\label{North Dakota}
According to module theory, 
it holds that any $\Q[t, t^{-1}]$-module  is congruent to 
$$(\Q[t, t^{-1}]/{\lambda_1})\oplus\cdot\cdot\cdot\oplus
(\Q[t, t^{-1}]/{\lambda_l})
\oplus(\oplus^k\Q[t, t^{-1}]),$$ 
where we have the following:

\smallbreak
\noindent
(1) 
$\lambda_*\in$ $\Q[t, t^{-1}]$ is not zero,

\smallbreak
\noindent
(2)
$\lambda_*$ is not the $\Q[t, t^{-1}]$-balanced class of 1, 

\smallbreak
\noindent
(3) 
$k$ is the rank of the free part.     
%See e.g. \S7 of [Kw]. 
%$\Z[t, t^{-1}]$-modules do not split like this in general. Refer to page 211 of [Rolfsen].  

\smallbreak
Two polynomials, $f(t)$ and $g(t), \in\Q[t, t^{-1}]$ are said to be 
{\it $\Q[t, t^{-1}]$-balanced } 
%(written $f\doteq g$) 
if there is an integer $n$ 
and a nonzero rational number $r$ 
such that 
$f(t)=r\cdot t^n\cdot g(t)$.

\smallbreak
Let $H_p(\widetilde X;\Q)$ be 
as above. 
Then the {\it $\Bbb Q[t, t^{-1}]$-$p$-Alexander polynomial} is \newline
the $\Bbb Q[t, t^{-1}]$-balanced class of
$
\left\{
\begin{array}{ll}
\mbox{the product $\lambda_1\cdot...\cdot\lambda_l$} & 
\mbox{if $k=0$ and  $H_p(\widetilde X;\Q)\ncong0$, } \\
\mbox{0} & \mbox{if $k\neq0$, }  \\
\mbox{1} & \mbox{if $H_p(\widetilde X;\Q)\cong0$.}      
\end{array}
\right.
$ 
\end{defn}

In this paper manifolds (resp. submanifolds) include 
manifolds-with-boundary (resp. submanifolds-with-boundary).   
A {\it Seifert hypersurface} for 
an $n$-dimensional oriented closed submanifold  $K$ in $S^{n+2}$ 
is 
an $(n+1)$-dimensional oriented connected compact submanifold in $S^{n+2}$ 
whose boundary is $K$ ($n\in\N$).
Note that Seifert hypersurfaces exist by obstruction theory 
(see \cite[pages 49 and 50]{Kirby}). 
Note that there are two cases that 
$K$ is not connected and that $K$ is connected.

\bigbreak
Let $V$ be a Seifert hypersurface for the above $n$-submanifold $K$.    
Let $x_1,..., x_\mu$ be $p$-cycles in $V$ 
which compose a basis of $H_p(V;\Z)$/Tor. 
Let $y_1,..., y_\nu$ be $(n+1-p)$-cycles in $V$ 
which compose a basis of $H_{n+1-p}(V;\Z)$/Tor. 
Push $y_i$ into the positive (resp. negative) direction of the normal bundle of $V$. 
Call it $y_i^{+}$ (resp. $y_i^{-}$).  
%Push $y_i$ into the negative direction of the normal bundle of $V$. Call it $y_i^{-}$.  
A $(p,n+1-p)$-{\it positive Seifert matrix} for the above submanifold $K$ associated with $V$ represented by 
an ordered basis,  
$\{x_1,..., x_\mu\}$,   and 
an ordered basis,  
$\{y_1,..., y_\nu\}$,  
is a $(\mu\x\nu)$-matrix 
$$S=(s_{ij})=({\mathrm{lk}}(x_i, y_j^{+})).$$   
A $(p,n+1-p)$-{\it negative Seifert matrix} 
for the above submanifold $K$ associated with $V$ represented by 
an ordered basis,  
$\{x_1,..., x_\mu\}$,   and 
an ordered basis,  
$\{y_1,..., y_\nu\}$,  
is a matrix 
$$N=(n_{ij})=({\mathrm{lk}}(x_i, y_j^{-})).$$   
We have the following: 
Let $S$ and $N$ be as above. 
Then $S-N$ represents the map 
 $\{H_{p}(V;\Z)$/Tor\} $\x \{H_{n+1-p}(V;\Z)$/Tor\} $\rightarrow\Z$ 
which is defined by the intersection product. 
We call  $t\cdot S-N$ the $(p,n+1-p)$-{\it Alexander matrix}   
for $K$ associated with $V$ represented by 
an ordered basis,  
$\{x_1,..., x_\mu\}$,   and 
an ordered basis,  
$\{y_1,..., y_\nu\}$. 
$S$ and $N$ (resp. $S$ and $t\cdot S-N$,  $N$ and $t\cdot S-N$)  
are said to be {\it related}  
if 
$S$ and $N$ (resp. $S$ and $t\cdot S-N$,  $N$ and $t\cdot S-N$)  
are defined by using the same $V$, the same $\{x_1,..., x_\mu\}$, 
and the same $\{y_1,..., y_\nu\}$.   
%We call $(S,N)$ a related pair of $p$-Seifert matrices.
We sometimes abbreviate 
$(p,n+1-p)$-positive Seifert matrix  
(resp. $(p,n+1-p)$-negative Seifert matrix, $(p,n+1-p)$-Alexander matrix)   
to 
$p$-Seifert matrix  
(resp. $p$-negative Seifert matrix, $p$-Alexander matrix) 
when it is clear from the context.

\begin{pr}\label{square}  
Let $K$ be an $n$-dimensional oriented closed submanifold $\subset S^{n+2}$. 
Let 
$S_p$ 
$($resp. $N_p)$ 
be  
a $p$-positive 
$($resp. negative$)$  
Seifert matrix for $K$ associated with $V$ represented by 
an ordered basis,  
$\{x_1,..., x_\mu\}$,   and 
an ordered basis,  
$\{y_1,..., y_\nu\}$. 
Suppose 
$$\mu=\nu.$$
Suppose that the homomorphism on 
$H_{p-1}(\amalg_{-\infty}^{\infty}V\x[-1,1];\Q)
\to
H_{p-1}(\amalg_{-\infty}^{\infty}Y;\Q)$ 
defined by a $(p-1)$-Alexander matrix is injective. 
Then 
the $p$-$\Q[t, t^{-1}]$-Alexander polynomial is 
the $\Q[t, t^{-1}]$-balanced class of 
`the determinant of $p$-Alexander matrix'  
$${\rm{det}} (t\cdot S_p-N_p).$$ 
\end{pr}

\noindent
{\bf Note.}  
Of course  $\mu\neq\nu$ in general. 

\bigbreak
\noindent
{\bf Proof of Proposition \ref{square}.} 
Take the above $X=\overline{S^{n+2}-(K\x D^2)}$, $\widetilde X$, $V$.   
Let $V\x[-1,1]$ be the tubular neighborhood of $V$ in $X$. 
Let $Y=X-V$.  Consider the Mayer-Vietoris exact sequence: 
$$H_\natural(\amalg_{-\infty}^{\infty}V\x[-1,1];\Q)
\to
H_\natural(\amalg_{-\infty}^{\infty}Y;\Q)
\to  H_\natural(\widetilde X;\Q),$$  

\noindent 
where 
$\amalg_{-\infty}^{\infty}V\x[-1,1]$ is the lift of  $V\x[-1,1]$,  
and where  
$\amalg_{-\infty}^{\infty}Y$ is the lift of $Y$. 

This completes the proof of Proposition \ref{square}. \qed

\begin{pr}\label{tsuika}
Let $N_p$ be 
a $(p,n+1-p)$-{negative Seifert matrix} for $K$ associated with $V$ represented by 
an ordered basis,  
$\{x_1,..., x_\mu\}$,   and 
an ordered basis,  
$\{y_1,..., y_\nu\}$.  
Let $S_{n+1-p}$ be 
a $(n+1-p,p)$-{positive Seifert matrix} for $K$ associated with $V$ represented by 
an ordered basis,  
$\{y_1,..., y_\nu\}$,   and 
an ordered basis,  
$\{x_1,..., x_\mu\}$. Then we have 
$$N_p=(-1)^{p\cdot n+1}\cdot S_{n+1-p}.$$ 
\end{pr}

\noindent
{\bf Proof of Proposition \ref{tsuika}.}
By the definition of 
$x_i^{+}$ and  $y_i^{-}$, 
${\mathrm{lk}}(y_i, x_j^{+})$ 
$={\mathrm{lk}}(y_i^{-}, x_j)$. 
By \cite[page 541]{Levinepol}, 
${\mathrm{lk}}(y_i^{-}, x_j)$
$=(-1)^{p(n+1-p)+1}{\mathrm{lk}}(x_j, y_i^{-})$. 
Note that $p(1-p)$ is an even number.\qed

\bigbreak
Proposition \ref{tsuika} implies Proposition \ref{Mars}. 

\begin{pr}\label{Mars} 
Let $K$ be a $(2m+1)$-dimensional closed oriented submanifold  $\subset S^{2m+3}$. 
Let $S$ be an $(m+1,m+1)$-Seifert matrix. 
Then we have 
$$S=(-1)^{m}\cdot ^t \hskip-1mm S.$$ 
%It holds that 
%$A-(-1)^{m}\cdot ^t \hskip-1mm A$ represents the map 
% $\{H_{m+1}(V;\Z)$/Tor\} $\x \{H_{m+1}(V;\Z)$/Tor\} $\rightarrow\Z$,  
%which is defined by the intersection product. 
\end{pr}

Let $K$ be a $(4k+1)$-dimensional spherical knot 
$(k\in\N\cup\{0\})$. 
We regard naturally 
$(H_{2k+1}(V;\Z)/\mathrm{Tor}) \otimes {\Z_2}$ 
as a subgroup of $H_{2k+1}(V;\Z_2)$.  
Then we can take a basis  $\{x_1,..., x_\nu, y_1,..., y_\nu\}$ 
of 
$(H_{2k+1}(V;\Z)/\mathrm{Tor}) \otimes\Z_2$ 
such that 
$x_i\cdot x_j=0$, $y_i\cdot y_j=0$,  
$x_i\cdot y_j=\delta_{ij}$ 
for any pair $(i, j)$, 
where $\cdot$ denotes the $\Z_2$-intersection product. 
The {\it Arf  invariant} of $K$ is 
$$\Big(
\sum_{i=1}^\nu   {\mathrm{lk}}(x_i, x_i^{+})\cdot  {\mathrm{lk}}(y_i, y_i^{+})
\Big) 
\hskip1mm  \text{mod 2 }.$$ 

\bigbreak
Let $L=(L_1,..., L_\mu)$ be a $(4k+1)$-link 
$(k\in\N\cup\{0\}. \quad \mu\in\N-\{1\}.)$.  
We define the {Arf } invariant of $L$.  
There are two cases. 

\begin{itemize}
\item[(1)]
 Let $4k+1\geq5$. 
The Arf  invariant of $L$ is defined in the same manner as the knot case. 

\item[(2)]
 Let $4k+1=1$. 
See Appendix of \cite{Kirby} and \cite[Note right above Note 1.2.1]{Ogasa02}.   
\end{itemize}

\bigbreak
\section{Main theorems in the 2-dimensional case}\label{m2}
\noindent 
Two polynomials,  $f(t)$ and $g(t), \in\Z[t, t^{-1}]$ are said to be 
{\it $\Z[t, t^{-1}]$-balanced} 
%(written $f\doteq g$) 
if there is an integer $n$ 
such that 
$f(t)=\pm t^n\cdot g(t)$.

\begin{thm}\label{Z}
Let $L=(L_1,...,L_m)$ be a 2-dimensional closed oriented submanifold  $\subset S^4$. 
Let each $L_i$ be connected.   
Let $g_i$ be the genus of $L_i$.   
Let $m=1+\Sigma_1^m g_i$. 
Let $S_\nu(V)$ be a positive $\nu$-Seifert matrix associated with a 
Seifert hypersurface for $L$  $(\nu=1,2)$. 
Let $N_\nu(V)$ be its related negative $\nu$-Seifert matrix. 
Let $t\cdot S_\nu(V)-N_\nu(V)$ be their related $\nu$-Alexander matrix. 
Then 
the $\Z[t, t^{-1}]$-balanced class of %`the determinant of a $\nu$-Alexander matrix'  
$${\rm{det}} (t\cdot S_\nu(V)-N_\nu(V))$$ 
is a topological invariant of $L$. 
%, where $S_\nu$ is a positive $\nu$-Seifert matrix and $N_\nu$ is a related negative $\nu$-Seifert matrix. 
%%defined for the same Seifert hypersurface and basis of the homology group of the Seifert hypersurface. 
\end{thm}

\noindent
{\bf Note.}  (1)    
We have $b_0(L)=\frac{1}{2}b_1(L)+1$, where $b_j$ is the $j$-th betti number. 

\smallbreak\noindent(2) 
Since any Seifert hypersurface is connected by the definition, the $\Q[t,t^{-1}]$-balanced class of ${\rm{det}} (t\cdot S_\nu(L)-N_\nu(L))$ is determined by the $\Q[t,t^{-1}]$-module $H_1(\widetilde{X};\Q)$, where $\widetilde{X}$ is the infinite cyclic covering space of $L$. 

\smallbreak\noindent(3) 
Let $\Delta(t)$ be a polynomial which represents the $\nu$-$\Z[t,t^{-1}]$-Alexander polynomial of $L$.  
By the definition,  
the $\Q[t,t^{-1}]$-balanced class of $\Delta(t)$ 
is 
the $\nu$-$\Q[t,t^{-1}]$-Alexander polynomial associated with $L$.

\begin{defn}\label{Delaware}
 The $\Z[t,t^{-1}]$-balanced class defined in Theorem \ref{Z} is called the {\it $\nu$-$\Z[t,t^{-1}]$-Alexander polynomial} of $L$.  
\end{defn}

We call 
a 2-dimensional closed oriented submanifold  $L=(K_1, K_2)\subset S^4$ 
an {\it $(S^2, T^2)$-link} 
if $K_1$ (resp. $K_2$) is diffeomophic to $S^2$ (resp. $T^2$). 
Note that $b_0(L)=\frac{1}{2}b_1(L)+1$ holds. 

\begin{pr}\label{Idaho}
There are $(S^2, T^2)$-links, $A=(A_1,A_2)$ and $B=(B_1, B_2), \subset S^4$ 
such that 
their  $\nu$-$\Z[t,t^{-1}]$-Alexander polynomials are not equivalent  $(\nu=1,2)$  
and such that 
their  $\nu$-$\Q[t,t^{-1}]$-Alexander polynomials are equivalent. 
\end{pr}

\bigbreak  
%We prove Theorem \ref{two}, which is a generalization of Theorem 4.1 of  \cite{Ogasa09}.  The $\Q[t,t^{-1}]$-Alexander polynomial case of Theorem \ref{two} is Theorem 4.1 of  \cite{Ogasa09}.    Theorem \ref{two} is stronger than Theorem 4.1 of  \cite{Ogasa09} in the meaning of Propositions \ref{Idaho} and \ref{Mississippi}.

Recall the paragraph right before `Theorem \ref{two} cited in \S\ref{intro}'. 

\begin{thm}\label{two}
Let $L_+=(L_{+,1},...,L_{+,{m_+}})$ be a 2-dimensional closed oriented submanifold  $\subset S^4$. 
Let each $L_{+,i}$ be connected. 
Let $g_{+,i}$ be the genus of $L_{+,i}$. 
Let $m_+=1+\Sigma_1^{m_+} g_{+,i}$. 
Let $(L_+, L_-, L_0)$ be a $(1,2)$-pass-move-triple. 
Then there is a polynomial $\Delta_{\nu, K_*}(t) (*\newline=+,-,0$ and  $\nu=1,2)$ which represents 
the $\Z[t,t^{-1}]$-$\nu$-Alexander polynomial for $K_*$,  and  
we have the identity %$(\nu=1,2)$

$$\Delta_{\nu, L_+}(t)-\Delta_{\nu, L_-}(t)=(t-1)\cdot\Delta_{\nu, L_0}(t).$$ 
\end{thm} 

\noindent{\bf Note.} 
(1)   
If $m=1$, $K_+$ and $K_-$ are homeomorphic to $S^2$. 
Then $K_0$ is homeomorphic to $S^2$ or $S^2\amalg T^2$.
In each case Theorem \ref{two} holds. 
We do not need the condition on the homeomorphism type of $K_0$. 
%In general, if we have $H_0(K_0;\Bbb Z)\cong\Bbb Z^{\alpha'}$ and 
%$H_1(K_0;\Bbb Z)\cong\Bbb Z^{2\beta'}$, then $\alpha'=\beta'+1$. 

\smallbreak\noindent
(2) 
By  \cite[Proposition 4.3]{Ogasa09}, we cannot normalize 
the $\Z[t,t^{-1}]$-$\nu$-Alexander polynomial for $L_*$ 
to be compatible with the identities in Theorem \ref{two}. 
On the other hand, 
we can define the `normalized' Alexander polynomial 
in a case of the $(4k+1)$-dimensional case 
so that it is compatible with a local-move-identity. 
See Definition \ref{4k+1} and Theorem \ref{Tokyo} for detail.

\smallbreak\noindent
(3) 
If we remove the condition on the betti number, 
the identity does not hold in general by \cite[Proposition 4.2]{Ogasa09}.

\begin{pr}\label{Mississippi}
There are $(1,2)$-pass-move-triples 
$L=(L_+,L_-,L_0)$ and 
$L'=(L'_+,L'_-,L'_0)$ with the following properties:

\smallbreak\noindent$(1)$
 $b_0(L)=\frac{1}{2}b_1(L)+1$.
 $b_0(L')=\frac{1}{2}b_1(L')+1$.

\smallbreak\noindent$(2)$
$t-1$ $($resp. $t-1, 1)$   
represents  the $\Q[t, t^{-1}]$-Alexander polynomial of \newline
\quad\quad $L_+$ $($resp. $L_-, L_0)$. 

\smallbreak\noindent$(3)$ 
$t-1$ $($resp. $t-1, 1)$    
represents  the $\Q[t, t^{-1}]$-Alexander polynomial of \newline
\quad\quad $L'_+$ $($resp. $L'_-, L'_0)$. 

\smallbreak\noindent$(4)$
$4(t-1)$  $($resp. $3(t-1), 1)$  
represents  the $\Z[t, t^{-1}]$-Alexander polynomial of \newline
\quad\quad $L_+$ $($resp. $L_-, L_0)$. 

\smallbreak\noindent$(5)$  
$2(t-1)$ $($resp. $t-1, 1)$  
represents  the $\Z[t, t^{-1}]$-Alexander polynomial of \newline
\quad\quad $L'_+$ $($resp. $L'_-, L'_0)$. 
\end{pr}

\noindent{\bf Note.}  
Take two arbitrary different polynomials from 
$4(t-1)$, $3(t-1)$, $2(t-1)$, and $t-1$. 
Then they are not $\Z[t, t^{-1}]$-balanced but $\Q[t, t^{-1}]$-balanced.

\begin{pr}\label{Illinois}
Let $V$ be a Seifert hypersurface for an $(S^2, T^2)$-link $L=(K_1, K_2)$. 
Then we have the following: 

\smallbreak\noindent$(1)$
There is a basis,  $\{\tau_1,...\tau_n\}$,  of $H_2(V;\Z)$, 
where $n$ is an nonnegative integer.

\smallbreak\noindent$(2)$
There is a set $\{\sigma_1,...,\sigma_n\}\subset H_1(V;\Z)$ such that 
 $\{\pi(\sigma_1),...,\pi(\sigma_n)\}$ is a basis of $H_1(V;\Z)/{\rm Tor}$, 
where $\pi$ is the natural projection homomorphism $H_1(V;\Z)\to H_1(V;\Z)/{\rm Tor}$.  

\smallbreak\noindent$(3)$
The intersection product of $\sigma_i$ and $\tau_j$ in $V$ is 
$
\begin{cases}
0 & \text{if $i=1$.}\\
\delta_{ij} &\text{if $i\geq2$}\\
\end{cases}
$ 
\end{pr}

\begin{defn}\label{Arkansas}
Let  $L=\{K_1, K_2\}$ be an $(S^2, T^2)$-link and $V$ a Seifert hypersurface for $L$.   
Take sets,  $\{\sigma_1,...\sigma_n\}$ and $\{\tau_1,...,\tau_n\}$, 
 as in Proposition \ref{Illinois}.  
We define the {\it pseudo-alinking number} of $L$ to be 
the absolute value of the Seifert pairing of $\sigma_1$ and $\tau_1$. 
\end{defn}

\begin{thm}\label{alk}   
Let $L=(L_1, L_2)$ be an $(S^2, T^2)$-link $\subset S^4$. 
Let $\Delta_{1,L}(t)$ be a polynomial which represents 
the $\Z[t,t^{-1}]$-1-Alexander polynomial for $L$. 
Then 
$$
\begin{vmatrix}
\left.\displaystyle\frac{\Delta_{1,L}(t)}{(t-1)}\right|_{t=1}
\end{vmatrix}$$ 
is the pseudo-alinking number of $L$,  where $|\quad|$ denotes the absolute value.  
\end{thm}

\begin{defn}\label{Arizona} (\cite{Sato}.)
Let $L=(K_1,K_2)$ be an ordered closed oriented 2-dimensional submanifold $\subset S^4$. 
Let $K_1$ and $K_2$ be connected. 
Take any circle embedded in $K_i$. Give any orientation to the circle.  
Consider the linking number of the circle and $K_j$ ($i\neq j$). 
Make a set of all of the linking number. 
Then the set is regarded as $n\cdot\Z$ for a number $n\in\{0\}\cup\N$.  
Note that if $n=0$, then the set is $\{0\}$. 
We call this number $n$ the {\it alinking number} 
alk$(K_i\subset L, K_j\subset L)$  
of $K_i$ in $L$ around $K_j$ in $L$. 
Note that alk$(K_1\subset L, K_2\subset L)$ is not equal to 
alk$(K_2\subset L, K_1\subset L)$ in general.

Note that if $K_1$ is diffeomophic to $S^2$, 
alk$(K_1\subset L, K_2\subset L)$  
is zero. 
So, in this case, let the alinking number of $L$ mean 
alk$(K_2\subset L, K_1\subset L)$.    
\end{defn}

\begin{prob}\label{California}
Let $L=(K_1,K_2)$ be an $(S^2, T^2)$-link. 
Is the alinking number of $L$ different from the pseudo-alinking number of $L$ in general?  
\end{prob}

\noindent{\bf Note.} 
(1) Note \ref{Alaska} is related to Problem \ref{California}. 

\smallbreak\noindent(2)  
%\cite{Ogasasurfcobo} proved 
Proposition \ref{North Carolina} claims that 
the alinking number is a `surface-link cobordism' invariant. 
How about the pseudo-alinking number?  

\begin{thm}\label{Colorado}
Let $L=(K_1,K_2)$ be an $(S^2, T^2)$-link. 
Then the following three conditions are equivalent. 

\smallbreak\noindent$(1)$  
The alinking number of $L$ is zero. 

\smallbreak\noindent$(2)$    
The pseudo-alinking number of $L$ is zero. 

\smallbreak\noindent$(3)$  
$
\begin{vmatrix}
\left.\displaystyle\frac{\Delta_{1,L}(t)}{(t-1)}\right|_{t=1}
\end{vmatrix}$ is zero. 
\end{thm}

\begin{thm}\label{Connecticut}
Let $L=(K_1,K_2)$ be an $(S^2, T^2)$-link. 
Suppose that there is a Seifert hypersurfcace $V$ such that {\rm Tor}$H_1(V;\Z)\cong0$. 
Then the alinking number of $L$ is 
$
\begin{vmatrix}
\left.\displaystyle\frac{\Delta_{1,L}(t)}{(t-1)}\right|_{t=1}
\end{vmatrix}$. 
\end{thm}

%Ribbon の定義　 Ribbonならtor free のSeifertをはる　のを言う 俺の別のん、さらっと言って

Let $L=(K_1,K_2)$ be an $(S^2, T^2)$-link. 
We say that $L$ is {\it ribbon} if there is an immersion 
$f:B\amalg H\looparrowright S^4$ 
with the following properties, 
where $B$ is a 3-ball and $H$ is a genus one handle body: 
The self-intersection of $f$ consists of double points and 
is a disjoint union of 2-discs.  
Note that $f^{-1}(\text{each disc})$ is a disjoint union of two 2-discs. 
One of the two disc is included in the interior of $B\amalg H$. 
The intersection of $\partial(B\amalg H)$ and the other disc is the boundary of the other disc.

\begin{cor}\label{Georgia}
Let $L=(K_1,K_2)$ be a ribbon $(S^2, T^2)$-link. 
Then the alinking number of $L$ is 
$
\begin{vmatrix}
\left.\displaystyle\frac{\Delta_{1,L}(t)}{(t-1)}\right|_{t=1}
\end{vmatrix}$. 
\end{cor}

%\begin{pr}\label{Florida}  There is $(1,2)$-pass-move triples $(K_+, K_-, K_0)$ and $(J_+,J_-,J_0)$ with following properties: 

%\smallbreak\noindent$(1)$ $K_+$ (resp. $J_+$) is homeomorphic to $S^2$. 

%\smallbreak\noindent$(2)$ $K_0$ (resp. $J_0$) is a $(S^2,T^2)$-link. 

%\smallbreak\noindent$(3)$  The $\Z[t,t^{-1}]$-Alexander polynomial of $K_+ ($resp. $K_-, K_0)$ is the $\Z[t,t^{-1}]$-balanced class of   例 (resp. 例, $t).$  

%\smallbreak \noindent$(4)$  The $\Z[t,t^{-1}]$-Alexander polynomial of $J_+ ($resp. $J_-, J_0)$  is the $\Z[t,t^{-1}]$-balanced class of  例 $($resp. 例, $2t).$   \end{pr}

%\noindent{\bf Note.} The $\Q[t,t^{-1}]$-Alexander polynomial cannot distinguish $K_0$ from $J_0$.

If $(K_+, K_-, K_0)$ is a triple of 1-links as in \S1 and   
if $K_+$ and $K_-$ are 1-knots, 
then $K_0$ is always a 2-component 1-link. 
However the 2-dimensional case we have the following theorem. 

\begin{thm}\label{Aa}
There are $(1,2)$-pass-move-triples, 
$L=(L_+, L_-, L_0)$ and $L'=(L'_+, L'_-, L'_0)$,  
with the following properties: 

\smallbreak\noindent$(1)$
 $b_0(L)=\frac{1}{2}b_1(L)+1$.
 $b_0(L')=\frac{1}{2}b_1(L')+1$.

\smallbreak 
\noindent
$(2)$
$L_+$ and $L'_+$ are diffeomorphic to $S^2$. 
Hence 
$L_-$ and $L'_-$ are diffeomorphic to $S^2$. 

\smallbreak 
\noindent
$(3)$  
The $\Z[t,t^{-1}]$-1-Alexander polynomial of $L_+ ($resp. $L_-)$  
is equivalent to that of $L'_+ ($resp. 
\newline\hskip6mm
$L'_-)$.   

\smallbreak 
\noindent
$(4)$ 
 $L_0$ is diffeomorphic to $S^2\amalg T^2$. 
  $L'_0$ is diffeomorphic to $S^2$.   
Hence $L_0$ and $L'_0$ are not 
\newline\hskip6mm
diffeomorphic. 
\end{thm}

See  Corollary \ref{aletwi}.
In a $(4k+1)$-dimesional case 
we have similar situation to the 1-dimensional case, 
different from the 2-dimensional case.

\bigbreak
\section{Review of twist-moves on high dimensional knots}\label{Osaka}
\noindent 
In the following section (\S\ref{mhigh}) 
 we have high dimensional analogues of \S\ref{m2}. 
We prove a new local-move-identity for 
the `normalized' Alexander polynomial of 
a kind of $(4k+1)$-dimensional closed oriented submanifolds $\subset S^{4k+3}$ 
(Theorem \ref{Tokyo}). 
The local-move-identity is associated with the twist-move, which is reviewed in this section. 
We introduce the `pseudo-twinkling number'  
as an analogue of 
the pseudo-alinking number,  
the alinking number, and 
the linking number (Definition \ref{twinkling}). 
We show a relation between the `normalized' Alexander polynomial  
%of $n$-dimensional closed oriented submanifolds $\subset S^{n+2}$ 
and the pseudo-twinkling number (Theorem \ref{aletwi}).
The pseudo-twinkling number is an analogue of the pseudo-alinking number but 
a relation between  the  pseudo-twinkling number and the `normalized' Alexander polynomial 
in Corollary \ref{twi+-}  
is different from 
the relations between the pseudo-alinking number and the $\Z[t,t^{-1}]$-Alexander polynomial 
in \S\ref{m2}.

%\bigbreak
%Let $H_p(W;\Z)/{\rm Tor}\cong\Z^\mu$ (resp.  $H_{n+1-p}(W;\Z)/{\rm Tor}\cong\Z^\nu$). Let $x_1,...,x_\mu$ (resp. $y_1,...,y_\nu$)  be $p$-cycles (resp. $(n+1-p)$-cycles) which make a set of basis of $H_p(W;\Z)/{\rm Tor}\cong\Z^\mu$ (resp. $H_{n+1-p}(W;\Z)/{\rm Tor}\cong\Z^\nu$). Put $y_j$ in the positive (resp. negative) direction of the normal $\R$-bundle of $W$ in $S^{n+2}$, and call it $y_j^+$ (resp. $y_j^-$). The $\mu\x\nu$-matrix $$S(x_k, y_j)=({\rm lk}x_k, y_j^+) \quad{\rm (resp.} P(x_k, y_j)=({\rm lk}x_k, y_j^-))$$ is called a $(p,n+1-p)$-positive (resp. negative) Seifert matrix of $L$ associated with $W$ and basis $\{x_k\}$ and $\{x_j\}$ ($k=1,...,\nu,\quad j=1,...,\mu$).

\begin{figure}
\vskip-20mm

\includegraphics[width=9cm]{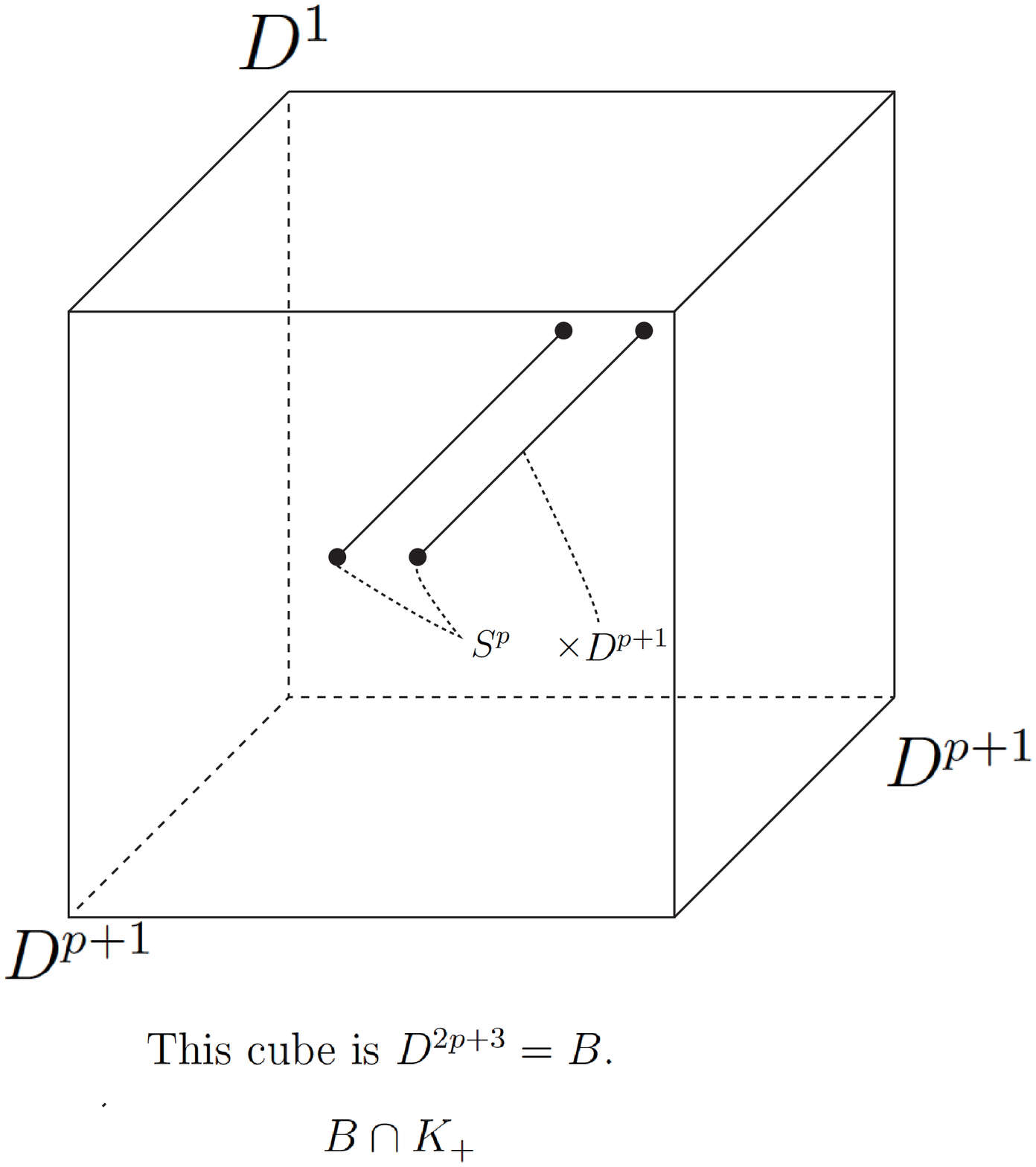}

{\bf Figure \ref{Osaka}.1.(1):  A twist-move-triple} 

\vskip7mm

\includegraphics[width=9cm]{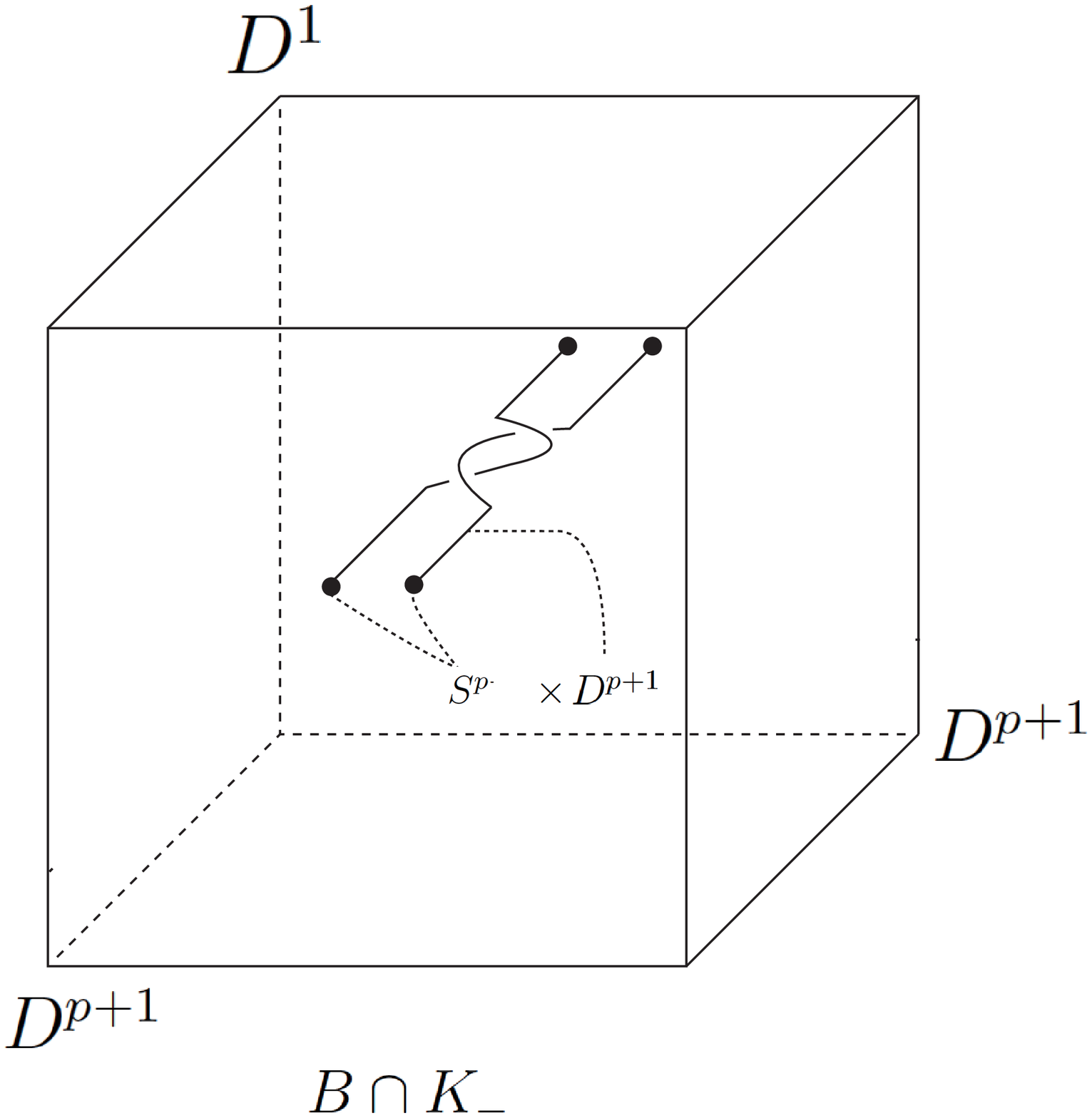}

\vskip-10mm

{\bf Figure \ref{Osaka}.1.(2):  A twist-move-triple} 

\end{figure}

\begin{figure}
\vskip-20mm

\includegraphics[width=10cm]{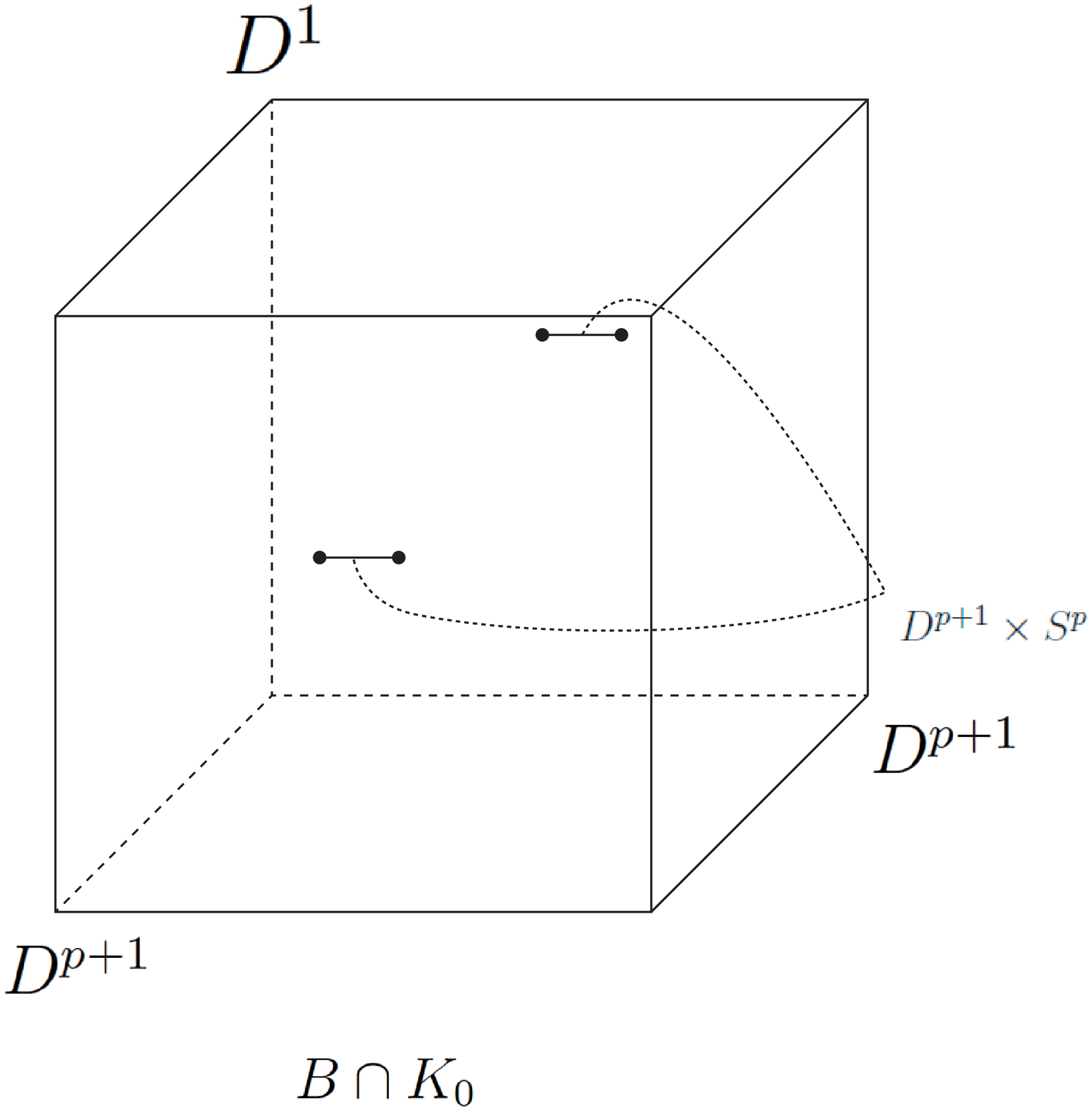}

Figure \ref{Osaka}.1.(3): {\bf A twist-move-triple} 
\end{figure}

\bigbreak
%We state the theorems after w
We review twist-moves on high dimensional knots in this section. 
(Note:  In \cite{Ogasa09} the twist-move is called the $XXII$-move.)
Figure \ref{Osaka}.1, which consists of the three figures (1), (2) and (3),  is a diagram of a twist-move-triple. 
Confirm the following:  if $p=0$, the twist-move is the crossing change on 1-links 
and 
%$K_+, K_-$ and $K_0$ are 
Figure \ref{Osaka}.1 is one drawn in the first paragraph in \S1.

\bigbreak
Let $K_+$,  $K_-$, $K_0$ be 
$(2p+1)$-dimensional closed oriented submanifolds $\subset S^{2p+3}$ 
($p\newline\in\N\cup\{0\}$). 
Let $B$ be a $(2p+3)$-ball trivially embedded in $S^{2p+3}$. 
Suppose that $K_+$ coincides with  $K_-$ (resp. $K_0$) in $\overline{S^{2p+3}-B}$. 
Take 
a single $(2p+2)$-dimensional $(p+1)$-handle $h_+$ (resp. $h_-$) embedded in $B$
such that (the handle)$\cap\partial B$ 
is the attaching part of the handle. 

\smallbreak\noindent
{\bf Note.} \cite{Haefligerunknot, HaefligerFr, Whitney, Whitneytrick, Wu} etc. 
imply that the core of $h_+$ (resp. $h_-$) is trivially embedded in $B$ under the above condition.     

\smallbreak
Suppose that $(h_+-$ its attaching part)$\cap(h_--$ its attaching part)$=\phi$.  
Suppose that their attaching parts coincide.  
Thus we can suppose that  
we regard $h_+\cup h_-$ as an oriented $(2p+2)$-submanifold $\subset S^{2p+3}$
if we give the opposite orientation to $h_-$.  
Then we can define a $(p+1)$-Seifert matrix for the $(2p+2)$-submanifold $h_+\cup h_-$.  
%by using a $(p+1)$-cycle.   
We can suppose that the $(p+1)$-Seifert matrix of 
$\partial(h_+\cup h_-)$ associated with $h_+\cup h_-$  
is a $1\x1$-matrix $(1)$.

\begin{note}\label{aomori} 
In the case of the twist-move on the $(4k+1)$-dimensional submanifolds 
we can distinguish $K_+$ from $K_-$ because 
the Seifert matrix is a $1\x1$-matrix $(1)$ 
even if we change the orientation of $h_+\cup h_-$. 
On the other hand, in the $(1,2)$-pass-move-triple case 
we cannot distinguish $K_+$ from $K_-$.  See Note \ref{hokkaido}. 
\end{note}

\noindent
{\bf Note.}   
Suppose that $p$ is an odd natural number, and let $p=2k+1$. 
The twist-move for $(4k+3)$-submanifolds $\subset S^{4k+5}$ 
($4k+3\in\N$, \quad $k\in\N\cup\{0\}$)   
has the following property: 
Suppose that $K_+$ is made into $K_-$  by the twist-move.   
Suppose that $K_+$ is PL homeomophic to the standard sphere.    
Then 
%$K_-$ is not PL homeomophic to the standard sphere in general.  Furthermore 
$H_*(K_-;\Z)$ is not congruent to $H_*(K_+;\Z)$ in general.   
Example: 
Make a Seifert hypersurface $V_*$  for a 3-knot $K_*$ ($*=+,-$) as follows.    
A framed link representation of $V_+$  is the Hopf link 
such that the framing of one component is zero 
and such that that of the other is two. 
A framed link representation of $V_-$  is the Hopf link 
such that the framing of each component is two.  
\bigbreak

Let $K_* (*=+,-)$ satisfy that  
$K_*\cap \mathrm{Int}B=$ 
$(\partial h_*-\partial B)$.  
Note the following. When we define $K_+$, $h_+$ exists in $B$ and $h_-$ does not exist in $B$. 
When we define $K_-$, $h_-$ exists in $B$ and $h_+$ does not exist in $B$. 
%
%
%Let $n=2p+1$. 
Let $P=K_+ \cap (S^{2p+3}-\mathrm{Int}B)$. 
Let $Q=h_+\cap \partial B$. 
Let  
$T=P\cup Q$. 
Then $T$ is 
an $(2p+1)$-dimensional oriented closed submanifold in $S^{2p+3}-\mathrm{Int}B$.
Let $K_0$ be 
%an $(2p+1)$-dimensional oriented closed submanifold 
$T$ in $S^{2p+3}$.
%Let $K_0$ be $\partial(V_+ -{\mathrm {Int}}B)$. 
Then we say that  an ordered set $(K_+$, $K_-$, $K_0)$ is related by a single {\it twist-move.}     
$(K_+$, $K_-$, $K_0)$ is called a {\it twist-move-triple}.
We say that 
$K_-$ (resp. $K_+$) 
is obtained from 
$K_+$ (resp. $K_-$) 
by a single 
 {\it negative-twist-move} 
(resp. {\it positive-twist-move})  
in $B$. 
%
%
%If $(K_+$, $K_-$, $K_0)$ is a twist-move-triple, then we also say that $(K_-$, $K_+$, $K_0)$ is a {\it twist-move-triple}.  If $K_+$ and $K_-$ differ by a single twist-move in $B$,  we also say that $K_-$ and $K_+$ differ by a single {\it twist-move} in $B$. 

See 
Figure \ref{Osaka}.2 
for a twist-move-triple of $(4k+1)$-knots. 

\begin{figure}
\includegraphics[width=14cm]{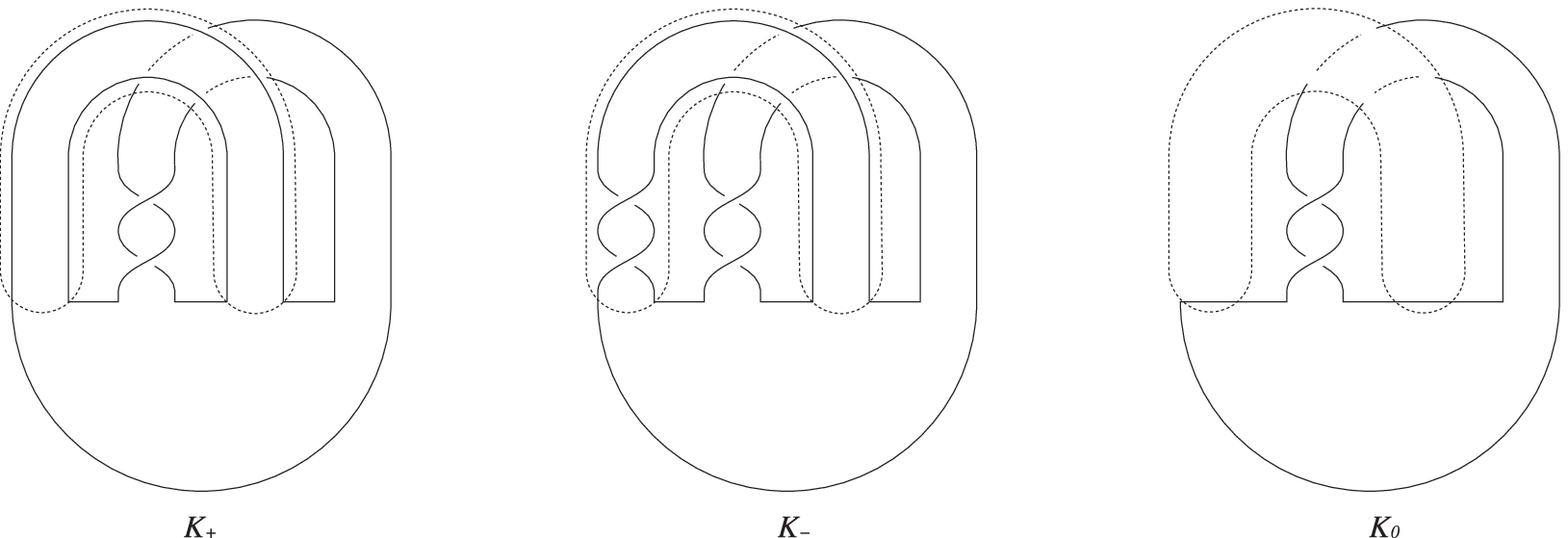}
\bigbreak
\hskip27mm Figure \ref{Osaka}.2: A twist-move-triple of $(4k+1)$-knots. 
\end{figure}

\bigbreak\noindent{\bf{Note.}} 
In the twist-move in the $(4k+1)$-dimensional case 
the homotopy type of $K_0$ is determined if $K_+$ is homotopy type equivalent to $S^{4k+1}$
by \cite{Adams, JW}.  
%{\color{cyan}{　←これほんまか？diffでなくてもよいか。考えよ  大丈夫。note.tex見よ}}
On the other hand, in the $(1,2)$-pass-move-triple case 
the homotopy type of $K_0$ is not determined even if $K_+$ is diffeomorphic to $S^2$. 
See Note to Theorem \ref{two}. 
\bigbreak

Let $(K_+, K_-, K_0)$ be related by a single twist-move in $B$.   
Then there is a Seifert hypersurface $V_*$ for $K_*$ ($*=+, -, 0$)  with the following properties. 

\smallbreak\noindent 
(1) $V_\sharp=V_0\cup h_\sharp$  ($\sharp=+,-$).
 $V_\sharp\cap B=h_\sharp$.  

\smallbreak
 \noindent 
(2) %\hskip5cm
$V_0\cap$ Int $B=\phi$. 
%\smallbreak
%
%\hskip5mm
$V_0\cap\partial B$ is the attaching part of $h_\sharp^p$. 

\noindent 
(The idea of the proof is the Thom-Pontrjagin construction.)

\bigbreak
 The ordered set ($V_+, V_-, V_0$) is called a {\it twist-move-triple of 
Seifert hypersurfaces} for  $(K_+, K_-, K_0)$. 
We say that 
$V_-$ (resp. $V_+$) is obtained from $V_+$ (resp. $V_-$) 
by a single 
 {\it negative-twist-move} 
(resp. {\it positive-twist-move})  
in $B$. 

See Figure \ref{Osaka}.3 for a twist-move-triple of Seifert hypersurfaces for $(4k+1)$-knots.

\begin{figure}
\includegraphics[width=14cm]{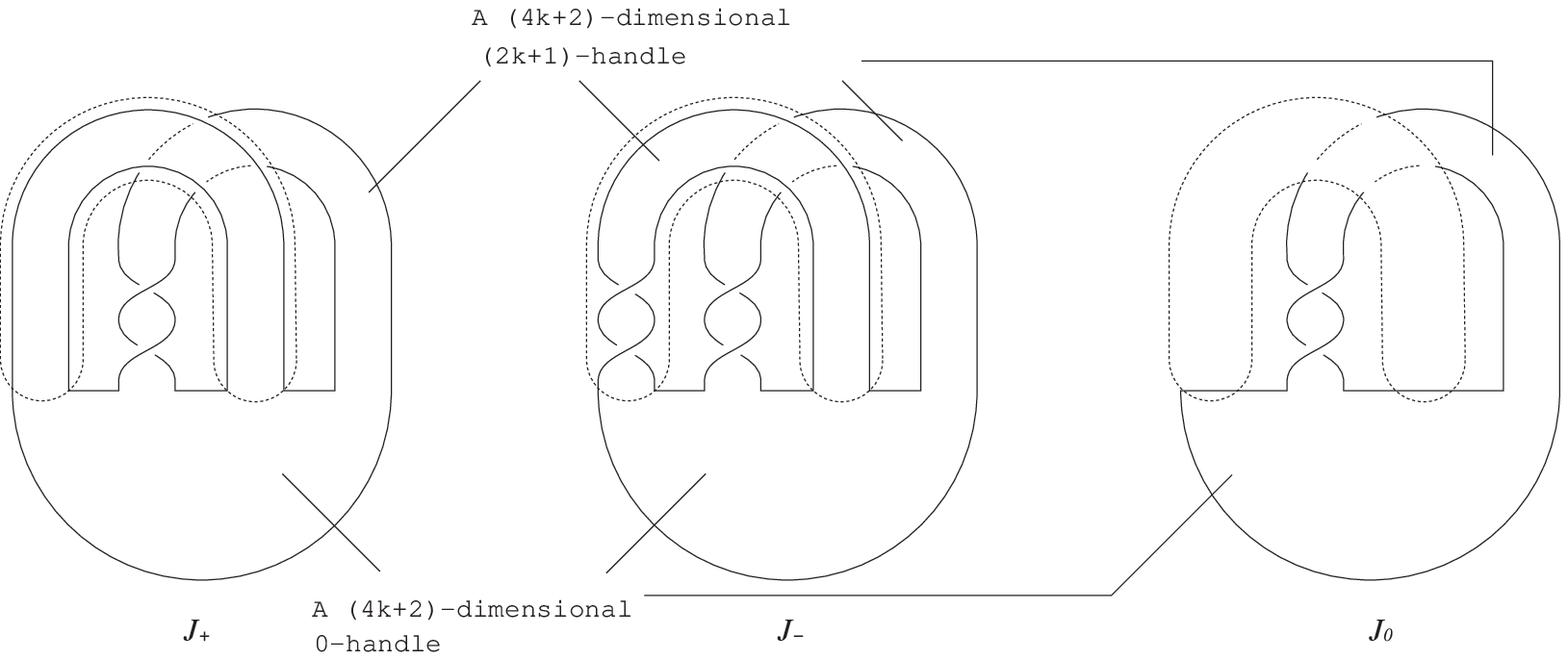}
\bigbreak
\hskip7mm Figure \ref{Osaka}.3: 
A twist-move-triple of Seifert hyopersurfaces for $(4k+1)$-knots 
\end{figure}

\bigbreak
In \cite{Ogasa98n, Ogasa09} 
we introduced the $(p,q)$-pass-move, which is a kind of local-moves.   
we found local-move-identities of the Alexander polynomial associated with 
the $(p,q)$-pass-move. 
We showed other relations between some invariants of knots and the $(p,q)$-pass-move. 
In \cite{KauffmanOgasa} we also proved such new results.

%%%%%%%%%%%%%%%%%%
%%%%%%%%%%%%%%%%%%
%%%%%%%%%%%%%%%%%%
\bigbreak
\section{Main theorems in the 4k+1 dimensional case}\label{mhigh}

%The $\Z[t, t^{-1}]$-balanced class of `the determinant of $(4k+1)$-Alexander matrix'  
%${\rm{det}} (t\cdot S_{4k+1}(L)-N_{4k+1}(L))$ is a topological invariant $(\nu=1,2)$ of $L$, 
%where $S_p(L)$ and $N_p(L)$ are related. 

\noindent 
We can define the `normalized' Alexander polynomial 
in a case of the $(4k+1)$-dimensional case 
so that it is compatible with a local-move-identity 
associated with the twist-move   
(see Definition \ref{4k+1} and Theorem \ref{Tokyo} for detail).   
On the other hand, 
in the 2-dimensional case
we cannot normalize the $\Z[t,t^{-1}]$-Alexander module  
so that it is compatible with the (1,2)-pass-move-identity 
(see \cite[Proposition 4.3]{Ogasa09}).

\begin{defn}\label{4k+1}
Let $k\in\{0\}\cup\N$. 
Let $K$ be a $(4k+1)$-dimensional closed oriented subamanifold $\subset S^{4k+3}$ 
whose homotopy type is $S^{4k+1}$. 
Let $V$ be a Seifert hypersurface for $K$. 
Let  $S_{2k+1}(V)$ be a $(2k+1)$-Seifert matrix and 
$N_{2k+1}(V)$ its related $(2k+1)$-negative Seifert matrix associated with 
a Seifert hypersurface $V$ for $K$. 
Call
$$\hat\Delta_{K}(t)={\rm{det}}(t^\frac{1}{2}\cdot S_{2k+1}(V)-t^\frac{-1}{2}\cdot N_{2k+1}(V))$$ 
the {\it normalized Alexander polynomial} for $K$.

Let $k\in\{0\}\cup\N$. 
Let $K$ be a $(4k+1)$-dimensional closed oriented subamanifold $\subset S^{4k+3}$ 
whose homotopy type is $S^{2k+1}\x S^{2k}$.  
Let  $S_{2k+1}(V)$ and $N_{2k+1}(V)$ be defined 
in the same manner as in the previous paragraph.  
%Let  $S_{2k+1}(V)$ be a $(2k+1)$-Seifert matrix and $N_{2k+1}(V)$ a related $(2k+1)$-negative Seifert matrix associated with Seifert hypersurface $V$ for $K$. 
Define the {\it normalized Alexander polynomial} $\hat\Delta_{K}(t)$ for $K$ 
to be 

$
\begin{cases}
{\rm{det}}(t^\frac{1}{2}\cdot S_{2k+1}(V)-t^\frac{-1}{2}\cdot N_{2k+1}(V))
&\\
\text{\hskip18mm 
if a $2k$-Alexander matrix associated with $V$ induces an injective map on} 
&\\\
\text{\hskip18mm $H_{2k}(\amalg_{-\infty}^{\infty}V\x[-1,1];\Q)
\to
H_{2k}(\amalg_{-\infty}^{\infty}Y;\Q)$,} & \\  
%the $(2p+1)$-$\Q[t,t^{-1}]$-Alexander polynomial is not the $\Q[t,t^{-1}]$-balanced class of zero, 
0& 
\text{\hskip-135mm else.
%if no $2k$-Alexander matrix associated with $V$ induces an injective map.
}   \\
%the $(2p+1)$-$\Q[t,t^{-1}]$-Alexander polynomial is the  $\Q[t,t^{-1}]$-balanced class of zero. 
\end{cases}
$
\end{defn}   

%{\color{cyan}{injectioveのところ、$\Z$ か $\Q$かをはっきりさせるほうがよい。たとえば、rankとかnulityとか使うのも手か}}

\noindent
{\bf Note.}   
(1) 
Recall that any $2k$-Alexander matrix associated with $V$
induces a homomorphism 
$H_{2k}(\amalg_{-\infty}^{\infty}V\x[-1,1];\Q)
\to
H_{2k}(\amalg_{-\infty}^{\infty}Y;\Q)$ 
as in Proof of Proposition \ref{square}. 

\smallbreak\noindent 
(2)By the definition of Alexander matrices we have the following: 
If a $2k$-Alexander matrix associated with $V$ induces (resp. does not induce) an injective map,  
then 
any (resp. no) $2k$-Alexander matrix associated with $V$ induces an injective map.

%分類かぶってるところしょうめいせなあかんのんやな

\begin{thm}\label{norm} 
The normalized Alexander polynomial $\hat\Delta_{K}(t)$ 
does not depend on the choice of $V$,  
and hence is a topological invariant.  
\end{thm}

\begin{thm}\label{Tokyo}
Let $K_+$ be a $(4k+1)$-knot $\subset S^{4k+3}$. 
Let $(K_+, K_-, K_0)$ be a twist-move-triple.
Then 
$$\hat\Delta_{K_+}(t)-\hat\Delta_{K_-}(t)
=(t^\frac{1}{2}-t^\frac{-1}{2})\cdot\hat\Delta_{K_0}(t),$$  
where $\hat\Delta_{K}(t)$ denotes the normalized Alexander polynomial of $K$. 
\end{thm}

See Figure \ref{Osaka}.2 
for an example of a twist-move-triple of $(4k+1)$-knots 
which satisfy the identity in Theorem \ref{Tokyo}.  
There,  we regard 
$S(V_*), N(V_*),$ and $\hat\Delta_{K_*}(t)$ as follows. 

 $
S(V_+)=
\begin{pmatrix}
0&-1\\
0&-1
\end{pmatrix}, 
N(V_+)= 
\begin{pmatrix}
0&0\\
-1&-1
\end{pmatrix}, 
\hat\Delta_{K_+}(t)=1, 
$

$
S(V_-)=
\begin{pmatrix}
-1&-1\\
0&-1
\end{pmatrix},  
N(V_-)= 
\begin{pmatrix}
-1&0\\
-1&-1
\end{pmatrix},  
\hat\Delta_{K_-}(t)=t+\frac{1}{t}-1,
$

$
S(V_0)=(-1),
N(V_0)=(-1),$  and 
$\hat\Delta_{K_+}(t)=-t^\frac{1}{2}+t^\frac{-1}{2}. $

\bigbreak

We say that $x\in H_i(X;\Z)$ is {\it order finite} (resp. {\it order infinite}) 
if $x\in {\rm Tor}H_i(X;\Z)$ (resp. $\notin {\rm Tor}H_i(X;\Z)$).  
Suppose that $x\in H_i(X;\Z)$ is nonzero and order finite. 
Let $p$ be the minimum number of $\{n\in\N|nx=0\}$.  
Then we say that $x$ is {\it order $p$}.  
We say that $x$ is {\it order zero} 
if $x=0\in H_i(X;\Z)$. 

\begin{defn}\label{divisible}   
%Let $x\in H_i(X;\Z).$    
We say that $x\in H_i(X;\Z)$ is {\it divisible} 
if $x$ is order infinite   %$x\notin {\rm Tor}H_i(X;\Z)$ 
and 
if there is $y\in H_i(X;\Z)$ such that 
$x=ny$ for an integer $n$ with the condition $|n|>1$. 
We suppose that  $x$ is order infinite  
when we say that $x\in H_i(X;\Z)$ is divisible (resp. non-divisible). 
%$x\notin {\rm Tor}H_i(X;\Z)$. 
If $y\in H_i(X;\Z)$ is order infinite, 
there is a non-divisible  $i$-cycle  
$z\in H_i(X;\Z)$    %, and $\notin {\rm Tor}H_i(X;\Z)$ 
such that there is an integer $m$ with the condition $y=mz$ ($m$ may be $\pm1$). 
Call $z$ a {\it non-divisible $i$-cycle associated with $y$.}  
\end{defn}

\begin{defn}\label{twinkling}
Let $k\in\{0\}\cup\N.$ 
Let $K$ be a $(4k+1)$-dimensional closed oriented subamanifold $\subset S^{4k+3}$ 
whose homotopy type is $S^{2k}\x S^{2k+1}$.  
Let $V$ be a Seifert hypersurface for $K$. 
We define the {\it pseudo-twinkling number} of $K$ to be 

$
\begin{cases}
\text{s($\tau, \tau)$}&
\text{
if there is a non-divisible $(2k+1)$-cycle $\tau\subset V$ such that 
}\\

 & \text{
for any $(2k+1)$-cycle $\alpha\subset V$  
the intersection product $\tau\cdot\alpha$ in $V$ is zero,  
}\\

\text{0}& \text{ else,  
%if there is not such a cycle,
}\\
\end{cases}
$

\noindent
where s$(\alpha, \beta)$ denotes the Seifert paring of $(2k+1)$-cycles $\alpha$ and $\beta$. 
Note that if $k=0$, the twinkling number is the linking number.  
\end{defn}

\noindent{\bf Note.}  
We would define the `twinkling number' to be 
s$(\gamma, \gamma)$, where  $\gamma$ is a generator of $H_{2k+1}(S^{2k}\x S^{2k+1})$. 
So we call the above one the pseudo-twinkling number by an analogy of 
the relation between the alinking number and the pseudo-alinking number 
although we do not discuss  the twinkling number so much in this paper. 
The author does not know whether the twinkling number and the pseudo-twinkling number are non-equivalent in general. 
He could prove that 
if there is a Seifert hypersurface $V$ such that ${\rm Tor}H_*(V;\Z)\cong0$, 
they are equivalent.  
Note \ref{Hawaii} is related to this question. 
He thinks that we have results which are analogues of Theorem \ref{Connecticut} and Corollary \ref{Georgia}. 
He could prove that the pseudo-twinkling number is `submanifold-cobordism' invariant, where  
submanifold-cobordism is defined in a similar fashion to that of knot cobordism by using (the submanifold)$\x[0,1]$. 
(See the definition right before Proposition \ref{North Carolina}  
for an example of submanifold-cobordism.)
He does not think that the twinkling number is `submanifold-cobordism' invariant. 

%０になるのは一致かな  Tor free なら一致かな　
%simple なら一致
%を書くか？
%
%bPとの関係も書くか？

\begin{pr}\label{hongo}
The pseudo-twinkling number of $K$ does not depend on the choice of $V$ and that of $\tau$, 
and hence is a topological invariant.   

\end{pr}

\begin{thm}\label{aletwi}
Let $K$ be a $(4k+1)$-dimensional closed oriented subamanifold $\subset S^{4k+3}$ 
whose homotopy type is $S^{2k}\x S^{2k+1}$.  
Let $\hat\Delta_{K}(t)$ be the normalized Alexander polynomial of $K$.  
Then the pseudo-twinkling number of $K$ is 
$$\left.\displaystyle\frac{\hat\Delta_{K}(t)}{t^\frac{1}{2}-t^\frac{-1}{2}}\right|_{t=1}$$
\end{thm}

Proposition \ref{hongo} and Theorem \ref{aletwi} imply the following.

\begin{cor}\label{twi+-}
Let $K_+$ be a $(4k+1)$-knot $\subset S^{4k+3}$. 
Let $(K_+, K_-, K_0)$ be a twist-move-triple.
Then the pseudo-twinkling number of $K_0$ is 
$$\left.\displaystyle
\frac{ \hat\Delta_{K_+}(t)-\hat\Delta_{K_-}(t)}{(t^\frac{1}{2}-t^\frac{-1}{2})^2}
\right|_{t=1}$$
\end{cor}

\noindent{\bf Note.} 
(1)
Compare `Theorem \ref{aletwi} and Corollary \ref{twi+-}' 
with 
`Theorem \ref{alk} and Theorem \ref{Aa}'. 

\smallbreak\noindent(2) 
There is a relation among $\hat\Delta_{K}(t)$, the bP-subgroup, and the inertia group 
by way of  \cite[Theorem 3.4]{Ogasa09} and Corollary \ref{twi+-}. %Theorem \ref{alk}.   
%We can say that this relation is an analogue of Theorem \ref{Colorado}. 

\bigbreak
Other results in \cite{KauffmanOgasa, Ogasa09} written in the $\Q[t,t^{-1}]$-term 
could be generalized into $\Z[t,t^{-1}]$-term in some fashion without difficulty 
although we must take care of \cite[\S10]{KauffmanOgasa}.

\bigbreak
\section{Proof of results in \S\ref{m2}}\label{proo2}

\begin{defn}\label{subsur}
Let $X$ be an $x$-dimensional submanifold 
of an $m$-dimesional manifold $M$ ($x,m\in\N, x<m$). 
Suppose that we can embed $X\x [0,1]$ in $M$ so that $X\x\{0\}=X$. 
Suppose that 
an $(x+1)$-dimensional handle $h^p$ is embedded in $M$ and  
is attached to $X\x [0,1]$ ($p\in\N\cup\{0\}, 0\leqq p\leqq x$).  
Suppose that the attaching part of $h^p$ is embedded in $X\x\{1\}$. 
See Figure \ref{proo2}.1. 
Suppose that $h^p\cap(X\x[0,1])$ is only the attaching part of $h^p$.  
Let 
$X'$=
$\overline{\partial(h^p\cup(X\x[0,1]))-(X\x\{0\})}$.
Note that there are two cases, $\partial X=\phi$ and $\partial X\neq\phi$. 
%Note that $X$ is the bottom and $X'$ the top of this handle decomposition. 
Then we say that 
 $X'$ is obtained from $X$ 
by {\it the surgery by using the embedded handle $h^p$}. 
We do not say that we use $X\x[0,1]$ if there is no danger of confusion.  

\noindent
{\bf Note.}  
Of course we can define `embedded surgery' even if we cannot embed $X\x[0,1]$ in $M$. 
However we do not need the case in this paper.

\begin{figure}
%WinTpicVersion2.15
\unitlength 0.1in
\begin{picture}(29.70,17.59)(16.60,-24.83)
% POLYGON 0 0 3 0
% 6 1660 2234 4390 2227 4390 2605 1667 2605 1667 2605 1660 2234
% 
\special{pn 20}%
\special{pa 1660 1834}%
\special{pa 4390 1827}%
\special{pa 4390 2205}%
\special{pa 1667 2205}%
\special{pa 1667 2205}%
\special{pa 1660 1834}%
\special{fp}%
% SPLINE 2 0 3 0
% 23 2241 2227 2248 2150 2262 2017 2332 1709 2353 1660 2500 1359 2535 1268 2598 1247 2703 1198 2822 1170 2920 1149 2962 1142 3018 1128 3144 1135 3291 1191 3410 1317 3466 1492 3543 1786 3564 1961 3578 2101 3578 2150 3578 2227 3585 2227
% 
\special{pn 8}%
\special{pa 2241 1827}%
\special{pa 2244 1795}%
\special{pa 2247 1763}%
\special{pa 2250 1731}%
\special{pa 2253 1700}%
\special{pa 2256 1668}%
\special{pa 2260 1636}%
\special{pa 2264 1604}%
\special{pa 2268 1572}%
\special{pa 2273 1541}%
\special{pa 2278 1509}%
\special{pa 2284 1477}%
\special{pa 2291 1446}%
\special{pa 2298 1415}%
\special{pa 2307 1384}%
\special{pa 2316 1353}%
\special{pa 2327 1323}%
\special{pa 2338 1293}%
\special{pa 2351 1264}%
\special{pa 2365 1235}%
\special{pa 2381 1208}%
\special{pa 2396 1181}%
\special{pa 2413 1154}%
\special{pa 2429 1127}%
\special{pa 2444 1100}%
\special{pa 2459 1071}%
\special{pa 2473 1042}%
\special{pa 2485 1012}%
\special{pa 2495 979}%
\special{pa 2503 945}%
\special{pa 2511 911}%
\special{pa 2523 882}%
\special{pa 2544 863}%
\special{pa 2574 853}%
\special{pa 2607 844}%
\special{pa 2637 831}%
\special{pa 2665 817}%
\special{pa 2693 802}%
\special{pa 2723 791}%
\special{pa 2753 783}%
\special{pa 2785 777}%
\special{pa 2817 771}%
\special{pa 2848 764}%
\special{pa 2879 757}%
\special{pa 2911 751}%
\special{pa 2942 746}%
\special{pa 2974 739}%
\special{pa 3005 731}%
\special{pa 3036 725}%
\special{pa 3068 724}%
\special{pa 3100 727}%
\special{pa 3132 733}%
\special{pa 3164 740}%
\special{pa 3195 748}%
\special{pa 3225 759}%
\special{pa 3255 771}%
\special{pa 3283 786}%
\special{pa 3310 804}%
\special{pa 3335 824}%
\special{pa 3359 847}%
\special{pa 3380 871}%
\special{pa 3399 897}%
\special{pa 3414 925}%
\special{pa 3427 954}%
\special{pa 3438 985}%
\special{pa 3447 1015}%
\special{pa 3455 1047}%
\special{pa 3463 1078}%
\special{pa 3470 1109}%
\special{pa 3479 1140}%
\special{pa 3488 1171}%
\special{pa 3497 1202}%
\special{pa 3505 1233}%
\special{pa 3514 1264}%
\special{pa 3522 1295}%
\special{pa 3530 1326}%
\special{pa 3537 1357}%
\special{pa 3543 1388}%
\special{pa 3548 1419}%
\special{pa 3553 1450}%
\special{pa 3556 1482}%
\special{pa 3560 1514}%
\special{pa 3563 1546}%
\special{pa 3566 1579}%
\special{pa 3569 1612}%
\special{pa 3572 1644}%
\special{pa 3575 1675}%
\special{pa 3578 1705}%
\special{pa 3580 1736}%
\special{pa 3572 1775}%
\special{pa 3567 1813}%
\special{pa 3585 1827}%
\special{pa 3585 1827}%
\special{sp}%
% SPLINE 2 0 3 0
% 17 2535 2234 2633 1835 2654 1709 2696 1632 2738 1562 2794 1506 2871 1485 2941 1464 3011 1450 3088 1464 3151 1513 3200 1569 3256 1723 3263 1877 3270 1926 3291 2038 3284 2220
% 
\special{pn 8}%
\special{pa 2535 1834}%
\special{pa 2545 1804}%
\special{pa 2555 1773}%
\special{pa 2564 1743}%
\special{pa 2573 1712}%
\special{pa 2582 1681}%
\special{pa 2591 1650}%
\special{pa 2599 1619}%
\special{pa 2607 1588}%
\special{pa 2614 1557}%
\special{pa 2620 1526}%
\special{pa 2625 1494}%
\special{pa 2630 1462}%
\special{pa 2634 1430}%
\special{pa 2636 1397}%
\special{pa 2640 1365}%
\special{pa 2646 1334}%
\special{pa 2656 1304}%
\special{pa 2670 1275}%
\special{pa 2686 1248}%
\special{pa 2703 1220}%
\special{pa 2719 1193}%
\special{pa 2736 1165}%
\special{pa 2755 1139}%
\special{pa 2778 1116}%
\special{pa 2805 1101}%
\special{pa 2836 1092}%
\special{pa 2868 1086}%
\special{pa 2899 1077}%
\special{pa 2930 1068}%
\special{pa 2960 1058}%
\special{pa 2992 1052}%
\special{pa 3024 1050}%
\special{pa 3056 1053}%
\special{pa 3086 1063}%
\special{pa 3114 1080}%
\special{pa 3138 1101}%
\special{pa 3161 1123}%
\special{pa 3183 1147}%
\special{pa 3202 1172}%
\special{pa 3219 1200}%
\special{pa 3232 1229}%
\special{pa 3243 1260}%
\special{pa 3251 1291}%
\special{pa 3256 1323}%
\special{pa 3259 1355}%
\special{pa 3260 1387}%
\special{pa 3260 1419}%
\special{pa 3261 1451}%
\special{pa 3264 1483}%
\special{pa 3268 1515}%
\special{pa 3274 1546}%
\special{pa 3281 1578}%
\special{pa 3287 1609}%
\special{pa 3291 1641}%
\special{pa 3294 1672}%
\special{pa 3294 1704}%
\special{pa 3292 1736}%
\special{pa 3290 1768}%
\special{pa 3286 1800}%
\special{pa 3284 1820}%
\special{sp}%
% SPLINE 2 2 3 0
% 6 1954 2605 1961 2682 1947 2759 1912 2808 1891 2857 1884 2878
% 
\special{pn 8}%
\special{pa 1954 2205}%
\special{pa 1958 2237}%
\special{pa 1961 2268}%
\special{pa 1961 2301}%
\special{pa 1956 2333}%
\special{pa 1945 2363}%
\special{pa 1926 2389}%
\special{pa 1908 2415}%
\special{pa 1895 2445}%
\special{pa 1885 2475}%
\special{pa 1884 2478}%
\special{sp -0.045}%
% STR 2 0 3 0
% 3 1737 2983 1737 3053 2 0
% $X$
\put(17.3700,-26.5300){\makebox(0,0)[lb]{$X$}}%
% SPLINE 2 2 3 0
% 6 2460 1480 2390 1480 2310 1450 2220 1410 2180 1380 2160 1360
% 
\special{pn 8}%
\special{pa 2460 1080}%
\special{pa 2428 1082}%
\special{pa 2396 1081}%
\special{pa 2365 1073}%
\special{pa 2335 1061}%
\special{pa 2306 1048}%
\special{pa 2276 1036}%
\special{pa 2246 1024}%
\special{pa 2218 1009}%
\special{pa 2192 990}%
\special{pa 2168 969}%
\special{pa 2160 960}%
\special{sp -0.045}%
% STR 2 0 3 0
% 3 1930 1250 1930 1350 2 0
% $h^p$
\put(19.3000,-9.5000){\makebox(0,0)[lb]{$h^p$}}%
% SPLINE 2 2 3 0
% 19 4410 2240 4440 2250 4470 2280 4480 2310 4480 2320 4480 2360 4500 2390 4560 2410 4570 2410 4510 2410 4500 2420 4500 2470 4500 2500 4470 2560 4440 2600 4430 2600 4420 2600 4410 2600 4410 2600
% 
\special{pn 8}%
\special{pa 4410 1840}%
\special{pa 4440 1850}%
\special{pa 4464 1871}%
\special{pa 4478 1900}%
\special{pa 4479 1932}%
\special{pa 4481 1964}%
\special{pa 4500 1990}%
\special{pa 4527 2005}%
\special{pa 4561 2010}%
\special{pa 4562 2007}%
\special{pa 4516 2008}%
\special{pa 4497 2032}%
\special{pa 4499 2065}%
\special{pa 4500 2097}%
\special{pa 4490 2125}%
\special{pa 4473 2154}%
\special{pa 4457 2186}%
\special{pa 4432 2200}%
\special{pa 4410 2200}%
\special{sp -0.045}%
% STR 2 0 3 0
% 3 4630 2390 4630 2490 2 0
% [0,1]
\put(46.3000,-20.9000){\makebox(0,0)[lb]{[0,1]}}%
\end{picture}%

\bigbreak\bigbreak\quad\quad
{\bf Figure \ref{proo2}.1:} A handle $h^p$ is attached to $X\x[0,1]$. 
\end{figure}
\end{defn}

\noindent
{\bf Proof of Theorem \ref{Z}.}  
Let $V$ and $V'$ be Seifert hypersurfaces for $L$. 
Recall that $V$ and $V'$ are connected by the definition. 
It suffices to prove that the $\nu$-$\Z[t,t^{-1}]$-Alexander polynomial ($\nu=1,2$) defined by using $V$ is the same as that defined by using $V'$.

By the same manner as 
that in \cite[sections 4 and 5]{Levinesimp},  
and 
that in \cite[Proof of Claim 8.1]{KauffmanOgasaB},  
we have the following: 
There are 
(not necessarily connected) 3-dimensional compact oriented submanifolds 
$V=U_1,U_2,...,U_{u-1},U_u=V'$
$\subset S^4$ $(u\in\N)$
such that $\partial U_*=L$ 
and such that $U_{*+1}$ is obtained from $U_*$ $(2\leq*+1\leq u)$ 
by a surgery by using an embedded 4-dimensional handle. 

If some of $U_\natural$ are not connected, 
use 4-dimensional 1-handles and  then we can suppose that all $U_\natural$ are connected, 
that is, all $U_\natural$ are Seifert hypersurfaces for $L$.  
%There are 3-dimensional connected compact oriented submanifolds 
%$V=U'_1,...,U'_{u'}=V'$ $\subset S^4$ $(u'\in\N)$
%such that $U'_{*+1}$ is obtained from $U'_*$ by a surgery by using an embedded 4-dimensional handle. 

Therefore it suffices to prove the following case: 
$V'$ is obtained from $V$ by a surgery 
by using an embedded 4-dimensional $i$-handle $h^i$ ($i=1,2,3$).  
%Note that the dual handle of a 4-dimensional 3-handle is a 4-dimensional 1-handle. 

Lemmas \ref{Louisiana} and \ref{Maine} imply Theorem \ref{Z}. 

\begin{lem}\label{Louisiana}
 Theorem \ref{Z} holds in the case $i=1,3$. 
\end{lem}

\noindent{\bf Proof of Lemma \ref{Louisiana}.}
 $V'=V\sharp(S^1\x S^2)$ or $V=V'\sharp(S^1\x S^2)$ 
where $\sharp$ is the connected-sum. 
If  $V'=V\sharp(S^1\x S^2)$,  
an Alexander matrix $A(t)$ for $V$ is related to  
an Alexander matrix $A'(t)$ for $V'$  as follows. 
%$A'(1)$　は０か１やけど、言わんで良い 

$$
A(t)=
\arraycolsep5pt
\left(
\begin{array}{@{\,}c|cccc@{\,}}
\text{$t$ or $-t$}&*&\cdot&\cdot&*\\
\hline
0&&&&\\
\cdot&\multicolumn{4}{c}{\raisebox{-10pt}[0pt][0pt]{\Huge $A'(t)$}}\\
\cdot&&&&\\
0&&&&\\
\end{array}
\right), 
$$
\noindent 
where $A(t)$ is an $n\x n$-matrix and $A'(t)$ is an $(n-1)\x (n-1)$-matrix ($n\in\N$). 
If  $V=V'\sharp(S^1\x S^2)$,  
an Alexander matrix $A(t)$ for $V$ is related to  
an Alexander matrix $A'(t)$ for $V'$  as follows. 
%$A'(1)$　は０か１やけど、言わんで良い 

$$
A'(t)=
\arraycolsep5pt
\left(
\begin{array}{@{\,}c|cccc@{\,}}
\text{$t$ or $-t$}&*&\cdot&\cdot&*\\
\hline
0&&&&\\
\cdot&\multicolumn{4}{c}{\raisebox{-10pt}[0pt][0pt]{\Huge $A(t)$}}\\
\cdot&&&&\\
0&&&&\\
\end{array}
\right), 
$$
\noindent 
where $A'(t)$ is an $n\x n$-matrix and $A(t)$ is an $(n-1)\x (n-1)$-matrix.

%{\color{cyan}行か列かどっちがゼロ並ぶほうか区別要るか}

Hence 
det $A(t)$ is $\Z[t,t^{-1}]$-balanced to det $A'(t)$ in the both cases. 

This completes the proof of Lemma \ref{Louisiana}. \qed 

\begin{lem}\label{Maine}
Theorem \ref{Z} holds in the case $i=2$. 
\end{lem}

\noindent{\bf Proof of Lemma \ref{Maine}.}
Let $C\subset V$ be the core of the attaching part of $h^2$. 
Let $N(C)$ be the tubular neighborhood of $C$ in $V$. 
By Definition \ref{subsur}, 
there is $W=(V\times [0,1])\cup h^2$ which is embedded in $S^4$.  
Let $C'\subset V'$ be the core of the attaching part of the dual handle ${\bar h}^2$ 
of $h^2$. Note that ${\bar h}^2$ is a 4-dimensional 2-handle. 
Recall that ${\bar h}^2$ is attached to $V'$ 
and that $W=(V\times [0,1])\cup h^2=(V'\times [0,1])\cup{\bar h}^2$.

There are two cases: 
\smallbreak
\noindent
(1) 
$[C]\in H_1(V;\Z)$ is order finite. 

\smallbreak
\noindent
(2) 
$[C]\in H_1(V;\Z)$ is order infinite. 
\smallbreak

We divide these two cases into four cases. 
\smallbreak
\noindent
(1-1) 
$[C]\in H_1(V;\Z)$ is order finite. 
$[C']\in H_1(V;\Z)$ is order finite. 

\smallbreak
\noindent
(1-2) 
$[C]\in H_1(V;\Z)$ is order finite. 
$[C']\in H_1(V;\Z)$ is order infinite.

\smallbreak
\noindent
(2-1) 
$[C]\in H_1(V;\Z)$ is order infinite. 
For any closed oriented surface $F$ embedded in $V$,  it holds that  
the intersection product of $[F]\in H_2(V, \partial V;\Z)$ and $[C]\in H_1(V;\Z)$ in $V$ 
is zero.

\smallbreak
\noindent
(2-2)   
$[C]\in H_1(V;\Z)$ is order infinite. 
There is a closed oriented surface $F$ embedded in $V$ such that  
the intersection product of $[F]\in H_2(V, \partial V;\Z)$ and $[C]\in H_1(V;\Z)$  in $V$ 
is nonzero.

Lemmas \ref{Maryland}, \ref{Massachusetts}, \ref{Michigan}, and \ref{Minnesota} 
imply Lemma \ref{Maine}. 

\begin{lem}\label{Maryland}
Lemma \ref{Maine} holds in the $(1\text{-}1)$ case.  
\end{lem}

\noindent{\bf Proof of Lemma \ref{Maryland}.}
Proposition \ref{yatto} implies that the (1-1) case does not occur.

\begin{pr}\label{yatto}
If $[C]\in H_1(V;\Z)$ is order finite, 
then $C'\in H_1(V';\Z)$ is order infinite. 
\end{pr}  
 
\noindent{\bf Proof of Proposition \ref{yatto}.}  
Let $[C]$ (resp. $[C']$) be order $q$ (resp. $q'$), where $q, q'\newline\in\N\cup\{0\}$.   
Take a 2-chain $\alpha\subset V$ (resp.  $\alpha'\subset V'$)  
whose boundary is $q\cdot C$ (resp. $q'\cdot C'$).  

Hence $\alpha\cup((q\cdot C)\x[0,1])\cup(q\cdot(\text{the core of $h^2$}))$ 
 (resp. \newline
$\alpha'\cup((q'\cdot C')\x[0,1])\cup(q'\cdot(\text{the core of ${\bar h}^2$}))$ 
is a 2-cycle $\beta$ (resp. $\beta'$) $\subset W$.   
Note that  
$(q\cdot C)\x[0,1]\hookrightarrow V\x[0,1]$
(resp. $(q'\cdot C')\x[0,1]\hookrightarrow V'\x[0,1]$) 
is a level preserving embedding map. 
Recall that $(V\x[0,1])\cup h^2$ is diffeomorphic to 
$(V'\x[0,1])\cup \bar{h}^2$. 
Therefore the intersection product of $\beta$ and $\beta'$ in $W$ is nonzero. 
However, since $W$ is embedeed in $S^4$, this intersection product is zero. 
We arrived at a contradiction. 

This completes the proof of Proposition \ref{yatto}. \qed

This completes the proof of Lemma \ref{Maryland}. \qed

\begin{lem}\label{Massachusetts}
If Lemmas \ref{Michigan} and \ref{Minnesota} hold, 
Lemma \ref{Maine} holds in the $(1\text{-}2)$ case.  
 \end{lem}

\noindent{\bf Proof of Lemma \ref{Massachusetts}.}
Replace $V$ with $V'$, $h^2$ with $\bar{h}^2$, and $C$ with $C'$.  
Therefore the (1-2) case is true if the (2) case is true. 
Note that the (2) case consists of 
the (2-1) case and the (2-2) case. 

This completes the proof of Lemma \ref{Massachusetts}. \qed

\begin{lem}\label{Michigan}
Lemma \ref{Maine} holds in the $(2\text{-}1)$ case.  
 \end{lem}

\noindent{\bf Proof of Lemma \ref{Michigan}.}
We can suppose the following:
There are a positive $p$-Seifert matrix $S_p(V)$ and its related negative $p$-Seifert matrix $N_p(V)$ $(p=1,2)$ associated with $V$ such that 
a square matrix $t\cdot S_p(V)-N_p(V)$ has a row all of whose elements are zero as follows:

$$
\begin{pmatrix}
0&\cdot\cdot\cdot&0\\
*&\cdot\cdot\cdot&*\\
\cdot&\cdot\cdot\cdot&\cdot\\
*&\cdot\cdot\cdot&*\\
\end{pmatrix}.
$$

Hence det$(t\cdot S_p(V)-N_p(V))=0$. 
Hence the $\Q[t,t^{-1}]$-$p$-Alexander polynomial is 
the $\Q[t,t^{-1}]$-balanced class of zero. 
Hence we have the following: 
Let $S_p(V')$ be a positive $p$-Seifert matrix and $N_p(V')$ its related negative $p$-Seifert matrix  $(p=1,2)$ associated with $V'$. 
By Proposition \ref{square} and Notes (2) and (3) to Theorem \ref{Z}, 
it holds that \newline  
det $(t\cdot S_p(V')-N_p(V'))$ 
is $\Q[t,t^{-1}]$-balanced to zero.       
Hence 
det $(t\cdot S_p(V')-N_p(V'))$=0.  

Hence det $(t\cdot S_p(V)-N_p(V))$ and det $(t\cdot S_p(V')-N_p(V'))$ are  
not only $\Q[t,t^{-1}]$-balanced but also $\Z[t,t^{-1}]$-balanced. 

This completes the proof of Lemma \ref{Michigan}. \qed

\begin{lem}\label{Minnesota}
Lemma \ref{Maine} holds in the $(2\text{-}2)$ case.  
 \end{lem}

\noindent{\bf Proof of Lemma \ref{Minnesota}.}
Note that 
$W=(V\x[0,1])\cup h^2$ is diffeomorphic to 
$(V'\x[0,1])\cup {\bar h}^2$. 
Consider the exact sequence by a pair 
$((V\x[0,1])\cup h^2, V)$, where 
we regard $V$ as $V\x\{0\}$: 
$$
\cdot\cdot\cdot\to H_*(V;\Z)\to 
H_*((V\x[0,1])\cup h^2;\Z)\to 
H_*((V\x[0,1])\cup h^2, V;\Z)\to\cdot\cdot\cdot 
$$

and 
the exact sequence by a pair 
$((V'\x[0,1])\cup {\bar h}^2, V')$, where we regard $V'$ as $V'\x\{0\}$: 
$$
\cdot\cdot\cdot\to H_*(V';\Z)\to 
H_*((V'\x[0,1])\cup {\bar h}^2;\Z)\to 
H_*((V'\x[0,1])\cup {\bar h}^2, V';\Z)\to\cdot\cdot\cdot. 
$$

By the existence of $F$, 
$[C']\in H_{2k}(V'; \Z)$ is order finite.  

Let $\xi_1\in H_1(V;\Z)$ be a non-divisible 1-cycle associated with $[C]$. 
Let $\eta_1\in H_2(V;\Z)$ be a non-divisible 2-cycle associated with $[F]$. 
We can suppose the following: 

\smallbreak\noindent(1)
There is a set 
$\{\xi_1, \xi_2,...,\xi_n\}\subset H_1(V;\Z)$, where $n\in\N\cup\{0\}$. 
A set $\{\pi(\xi_1), \pi(\xi_2),...,\pi(\xi_n)\}$
is a basis of $H_1(V;\Z)/{\rm Tor}$,    
where $\pi$ is the natural epimorphism $H_1(V;\Z)\to H_1(V;\Z)/{\rm Tor}$. 

\smallbreak\noindent(2)
We can regard 
$\{\xi_2,...,\xi_n\}\subset H_1(V';\Z)$. 
$\{\pi(\xi_2),...,\pi(\xi_n)\}$  
is a basis of $H_1(V;\Z)/{\rm Tor}$.

\smallbreak\noindent(3)
There is a basis $\{\eta_1, \eta_2,...,\eta_n\}$ of $H_2(V;\Z)$.

\smallbreak\noindent(4)
We can regard $\{\eta_2,...,\eta_n\}$   
as a basis of $H_2(V;\Z)$.

\smallbreak\noindent(5)
Since $H_*(\partial V;\Z)$ is torsion free,  
the intersection product  $\xi_i\cdot\eta_j$ in $V$ (resp. in $V'$) is 
$
\begin{cases}
1& \text{if $i=j=1$}          \\  
\text{$\delta_{ij}$ or zero}  &\text{else.}  
\end{cases}
$

\bigbreak
Hence we have  the following: 
An Alexander matrix $A(t)$ for $V$ is associated with 
 $\{\xi_1,...,\xi_n\}$ and $\{\eta_1,...,\eta_n\}$. 
An Alexander matrix $A'(t)$ for $V'$ is associated with 
 $\{\xi_2,...,\xi_n\}$ and $\{\eta_2,...,\eta_n\}$. 
$A(t)$ is an $n\x n$-matrix. 
$A'(t)$ is an $(n-1)\x (n-1)$-matrix. 
%$A'(1)$　は０か１やけど、言わんで良い 
Seifert pairings s$(\xi_*, \eta_\#)$ ($2\leq*$ and $2\leq\#$) 
are not changed 
when we attach the 4-dimensional 2-handle $h^2$ to $V$.

$$
A(t)=
\arraycolsep5pt
\left(
\begin{array}{@{\,}c|cccc@{\,}}
\text{$t$ or $-t$}&0&\cdot&\cdot&0\\
\hline
*&&&&\\
\cdot&\multicolumn{4}{c}{\raisebox{-10pt}[0pt][0pt]{\Huge $A'(t)$}}\\
\cdot&&&&\\
*&&&&\\
\end{array}
\right), 
$$

Hence 
det $A(t)$ is $\Z[t,t^{-1}]$-balanced to det $A'(t)$. 

This completes the proof of Lemma \ref{Minnesota}. \qed

This completes the proof of Lemma \ref{Maine}.  \qed

This completes the proof of Theorem \ref{Z}. \qed

\bigbreak
\noindent
{\bf Note.}  
It is important that we can suppose that $\xi_1\cdot\eta_1=1.$    
If it does not hold, det $A(t)$ is not $\Z[t,t^{-1}]$-balanced to det $A'(t)$ in general. 
See the example in \cite[\S10]{KauffmanOgasa}.

\bigbreak
\noindent
{\bf Proof of Theorem \ref{two}.} 
In %Proof of Theorem 4.1 %%s 3.2 and 3.3 
\cite[Proof of Theorem 4.1]{Ogasa09}
we proved that 
there is a $\nu$-Alexander matrix $A_{\nu, L_*}(t)$ for  $L_*$ $(*=+,-,0$ and $\nu=1,2)$
such that 
$${\rm det}A_{\nu, L_+}(t)-{\rm det}A_{\nu, L_-}(t)=(t-1)\cdot{\rm det}A_{\nu, L_0}(t).$$ 
This fact and Theorem \ref{Z} imply Theorem \ref{two}. 
%
%２－で、球面でも$t^m$の不定さはのこる  S同値のはったときに$t^m$でる  公式とコンパチも無理
%
%２－は、恒等式と一致するように正規化するのは無理なんや　すぐ上のtt1の例　
%俺の別の論文の結果、引用しとく？　　
%
%Theorem \ref{4k+1}見よ　４k＋１の場合のんは正規化できる。引用しとくか
\qed

\bigbreak
\noindent
{\bf Proof of Proposition \ref{Illinois}.}
Consider the following exact sequence by a pair $(V, \partial V)$  
(Note that $\partial V=S^2\amalg T^2$. Here, $S^2$ denotes $K_1$ and $T^2$ $K_2$.): 
$$
\cdot\cdot\cdot\stackrel{\partial}\to
H_*(S^2\amalg T^2;\Z)\stackrel{\iota}\to 
H_*(V;\Z)\stackrel{\rho}\to 
H_*(V,\partial V;\Z)\stackrel{\partial}\to 
H_{*-1}(S^2\amalg T^2;\Z)\stackrel{\iota}\to\cdot\cdot\cdot. 
$$
We can take sets,  
$\{\sigma_1,...\sigma_n\}$ and  $\{\tau_1,...,\tau_n\}$, 
to satisfy the conditions (1)-(3) in Proposition \ref{Illinois}. 
\qed

\bigbreak
\noindent
{\bf Proof of Theorem \ref{alk}.}   
Take sets,   
$\{\sigma_1,...\sigma_n\}$ and  $\{\tau_1,...,\tau_n\}$, 
as in Proposition \ref{Illinois}. 
Then the 1-Alexander matrix $A(t)$ associated with the ordered sets, 
$\{\sigma_1,...\sigma_n\}$ and  $\{\tau_1,...,\tau_n\}$, 
is written as follows:  
$$
\arraycolsep5pt
\left(
\begin{array}{@{\,}c|cccc@{\,}}
(t-1)\cdot a_{11}&(t-1)\cdot a_{12}&\cdot&\cdot&(t-1)\cdot a_{1n}\\
\hline
(t-1)\cdot a_{21}&&&&\\
\cdot&\multicolumn{4}{c}{\raisebox{-10pt}[0pt][0pt]{\Huge $X(t)$}}\\
\cdot&&&&\\
(t-1)\cdot a_{n1}&&&&\\
\end{array}
\right), 
$$
where we have the following: 
$a_{ij}={\rm s}(\sigma_i,\tau_j)$.  
$X(1)=\delta_{ij}$. 
$|a_{11}|$ is the pseudo-alinking number. 
Hence 

$$
A(t)=(t-1)
\arraycolsep5pt
\left(
\begin{array}{@{\,}c|cccc@{\,}}
a_{11}&  a_{12}&\cdot&\cdot&  a_{1n}\\
\hline
(t-1)\cdot a_{21}&&&&\\
\cdot&\multicolumn{4}{c}{\raisebox{-10pt}[0pt][0pt]{\Huge $X(t)$}}\\
\cdot&&&&\\
(t-1)\cdot a_{n1}&&&&\\
\end{array}
\right). 
$$

Hence 
$$\left.\displaystyle\frac{A(t)}{t-1}\right|_{t=1}= 
\arraycolsep5pt
\left(
\begin{array}{@{\,}c|cccc@{\,}}
a_{11}&  a_{12}&\cdot&\cdot&  a_{1n}\\
\hline
0&&&&\\
\cdot&\multicolumn{4}{c}{\raisebox{-10pt}[0pt][0pt]{\Huge $X(1)$}}\\
\cdot&&&&\\
0&&&&\\
\end{array}
\right). 
$$

Hence 
$
\left|\left.\displaystyle\frac{{\rm det}A(t)}{t-1}\right|_{t=1}\right|
$ 
is 
$ 
|a_{11}| 
$,
and hence is the pseudo-alinking number. 

This completes the proof of Theorem \ref{alk}.  \qed

\bigbreak
\noindent
{\bf Proof of Theorem \ref{Colorado}.} 
Theorem \ref{alk} implies that $(2)\Leftrightarrow(3)$.

We prove that $(1)\Leftrightarrow(2)$.   
By the exact sequence in Proof of Proposition \ref{Illinois}, we have \newline
$H_2(V;\Z)\cong\Z^n$, $H_2(V, \partial V;\Z)\cong\Z^n$,  \newline 
$H_1(V;\Z)\cong\Z^n\oplus T$, $H_2(V, \partial V;\Z)\cong\Z^n\oplus T$,  
where $T$ is the torsion part. 

\noindent 
Note that  
$\rho: 
H_1(V;\Z)\stackrel{\rho}\to 
H_1(V,\partial V;\Z)
$
is not an isomorphism in general. 
See Note \ref{Alaska}.

There is a nonzero element $\alpha\in H_1(S^2\amalg T^2;\Z)$ such that 
$\iota(\alpha)$ is order finite, 
where $\iota$ is the homomorphism in the exact sequence in Proof of Proposition \ref{Illinois}. 
Note that $\alpha$ is represented by an embedded circle $\subset T^2$, 
and let the circle also be called $\alpha$. 
Let $\beta$ be an embedded circle in $T^2$ 
such that $\alpha$ intersects $\beta$ transversely at one point. 
The 1-cycle which represented by $\beta$ is also called $\beta$.

We prove $\iota(\beta)$ is order infinite in $V$. Reason: 
Suppose that $\iota(\beta)$ is order finite. 
Let $P$ (resp. $Q$) be a 2-cycle $\subset V$ whose boundary is $\alpha$ (resp. $\beta$). 
We can suppose that $P$ intersects $Q$ transversely. 
Take $\partial (P\cap Q)$. 
It is a boundary of a 1-cycle $P\cap Q$ and hence it is zero $\in H_0(T^2)$. 
However it is one point by the definition of $\beta$ hence it is not zero $\in H_0(T^2)$. 
We arrived at a contradiction.

Take 
$\{\sigma_1,...\sigma_n\}$ and  $\{\tau_1,...,\tau_n\}$ as in Proposition \ref{Illinois}. 
Since $\beta\in H_1(S^2\amalg T^2;\Z)$, 
$\iota(\beta)\cdot\tau_*=0$ for all $*$. 
Hence 
$\sigma_1$ is a non-divisible 1-cycle associated with $\iota(\beta)$. 
Hence 
$\iota(\beta)=k\cdot\sigma_1$ for a nonzero integer $k$.

Since $S^2\subset\partial V$  (recall that $S^2$ denotes $K_1$), 
the intersection product $\iota(S^2)\cdot \sigma_*=0$ for any $*$. 
Since $H_2(V,\partial V;\Z)$ is torsion-free, 
$\iota(S^2)=\tau_1$. 

Hence 
$|{\rm lk}(\beta, \tau_1)|=|{\rm lk}(\iota(\beta), \tau_1)|=|k\cdot{\rm lk}(\sigma_1, \tau_1)|$. 

Suppose that the alinking number of $L$ is zero. Hence lk$(\beta, \tau_1)=0$. Hence 
the pseudo-alinking number $|{\rm lk}(\sigma_1, \tau_1)|$ is zero.

Suppose that the pseudo-alinking number  $|{\rm lk}(\sigma_1, \tau_1)|$ is zero. 
Hence ${\rm lk}(\beta, \tau_1)=0$. 
We can use $\{\alpha, \beta\}$ as a basis of $H_1(T^2;\Z)$. 
Since $\iota(\alpha)$ is order finite in $V$, lk$(\alpha, \tau_1)=0.$  
Hence lk$(l\cdot\alpha+m\cdot\beta, \tau_1)=0$ for any pair of integers $(l,m)$. 
Hence the alinking number of $L$ is zero. 

Hence  $(1)\Leftrightarrow(2)$.   

This completes of the proof of Theorem \ref{Colorado}. \qed

\begin{note}\label{Alaska} 
If we define $V$ as in (1), 
$\iota:H_1(\partial V;\Z)\to H_1(V;\Z)$ has the property in (2). 

\smallbreak\noindent(1) 
Let $f:S^1\hookrightarrow S^1\x S^2$ be an embedding 
such that \newline
$f_*:H_1(S^1:\Z)\to H_1(S^1\x S^2-{\rm Int}B^3:\Z)$
carries 1 to $p$ ($|p|>1.$ $p\in\Z.$).  
Let $B$ be an embedded 3-ball in $S^1\x S^2$ such that $B\cap N(f(S^1))=\phi$, 
where $N(f(S^1))$ is the tubular neighborhood of $f(S^1)$ in $S^1\x S^2$.    
Let $V$ be $S^1\x S^2-{\rm Int}B^3-{\rm Int}N(f(S^1))$.  
Note that $\partial V=S^2\amalg T^2$. 

\smallbreak\noindent(2) 
There is a non-divisible cycle $\zeta\in H_1(\partial V;\Z)\cong\Z\oplus\Z_p$ 
associated with a cycle $\in\iota(H_1(V;\Z))$ such that $\zeta\notin\iota(H_1(V;\Z))$. 
We have what is written in the fifth line of Proof of Theorem \ref{Colorado}. 
\end{note}

Let $F$ and $G$ be oriented closed connected surface $\subset S^4$. 
An {\it $(F,G)$-link} is  
a 2-dimensional closed oriented submanifold $L=(J,K)\subset S^4$ 
such that $J$ (resp. $K$) is diffeomorphic to $F$ (resp. $G$).   
Let $L=(J,K)$ and $L'=(J', K')$ be $(F,G)$-links in $S^4$. 
We say that $L$ and $L'$ are {\it surface-link-cobordant} 
if there is an embedding map $f:(F\amalg G)\x[0,1]\hookrightarrow S^4\x[0,1]$ 
with the following properties: 

\noindent
For $t=0,1$, 
$f((F\amalg G)\x[0,1])\cap (S^4\x\{t\})$ is $f((F\amalg G)\x\{t\})$. 

\noindent
$f((F\amalg G)\x\{0\})$ in $S^4\x\{0\}$ is $L$. 

\noindent
$f((F\amalg G)\x\{1\})$ in $S^4\x\{1\}$ is $L'$. 

\cite[\S2]{Sato} proved that 
if two surface-links are surface-link-cobordant and 
the alinking number of one of the two is zero, 
 then the alinking number of the other is zero. 
We generalize it and prove the following:

\begin{pr}\label{North Carolina}
Let $L=(J,K)$ and $L'=(J', K')$ be $(F, G)$-links in $S^4$. 
Suppose that  $L$ and $L'$ are surface-link-cobordant. 
Then we have the following: 

\noindent \hskip3cm
${\rm alk}(J\subset L, K\subset L)={\rm alk}(J'\subset L', K'\subset L')$.   

\noindent \hskip3cm
${\rm alk}(K\subset L, J\subset L)={\rm alk}(K'\subset L', J'\subset L')$.   
\end{pr}

\noindent
{\bf Proof of Proposition \ref{North Carolina}.}  
Take a compact oriented 4-manifold $P$ such that 
$\partial P\newline=f(F\x[0,1]) 
\cup \text{(a Seifert hypersurface for $J$)} 
\cup \text{(a Seifert hypersurface for $J'$)}$ \newline
(resp. $f(G\x[0,1]) 
\cup \text{(a Seifert hypersurface for $K$)} 
\cup \text{(a Seifert hypersurface for $K'$)}$) 
Consider 
$P\cap f(G\x[0,1])$ 
(resp. $P\cap f(F\x[0,1])$
and Definition \ref{Arizona}. 
\qed

\bigbreak
\noindent
{\bf Proof of Theorem \ref{Connecticut}.}  
In Proof of Theorem \ref{Colorado}, 
since $H_1(V,\partial V;\Z)$ has a nontrivial torsion in general, 
 $\sigma_1$ is $\iota(\beta)$ or 
a non-divisible 1-cycle associated with $\iota(\beta)$. 
That is, 
$k$ in Proof of Theorem \ref{Colorado} is not $\pm1$ in general. 
Now,  
since $H_1(V,\partial V;\Z)$ is torsion-free, 
 $\iota(\beta)=\pm\sigma_1$ and
$|{\rm lk}(\beta, \tau_1)|=|{\rm lk}(\sigma_1, \tau_1)|$. 
Hence the alinking number is 
$|{\rm lk}(\sigma_1, \tau_1)|$. 
Hence the alinking number is the pseudo-alinking number.

By Theorem \ref{Colorado}, 
Theorem \ref{Connecticut} holds. \qed

\bigbreak
\noindent
{\bf Proof of Corollary \ref{Georgia}.}  
Use the isotopy which changes \cite[Figure 4.4]{Ogasa04}  into \cite[Figure 4.3]{Ogasa04}
and vice versa. 
Use 4-dimensinal 1-handles. 
We obtain a Seifert hypersurface for $L$ whose homology groups are torsion-free.\qed

%$
%H_*(S^2\amalg T^2;\Z)\stackrel{\iota}\to 
%H_*(V;\Z)\stackrel{\rho}\to 
%H_*(V,\partial V;\Z)\stackrel{\partial}\to 
%H_{*-1}(S^2\amalg T^2;\Z)\stackrel{\iota}
%$

\bigbreak
\noindent
{\bf Proof of Theorem \ref{Aa}.}  
There is a (1,2)-pass-move-triple $(L_+, L_-, L_0)$ with the following properties 
(see Figure \ref{proo2}.2):  %{\color{cyan}図} 
\begin{figure}
   \includegraphics[width=12cm]{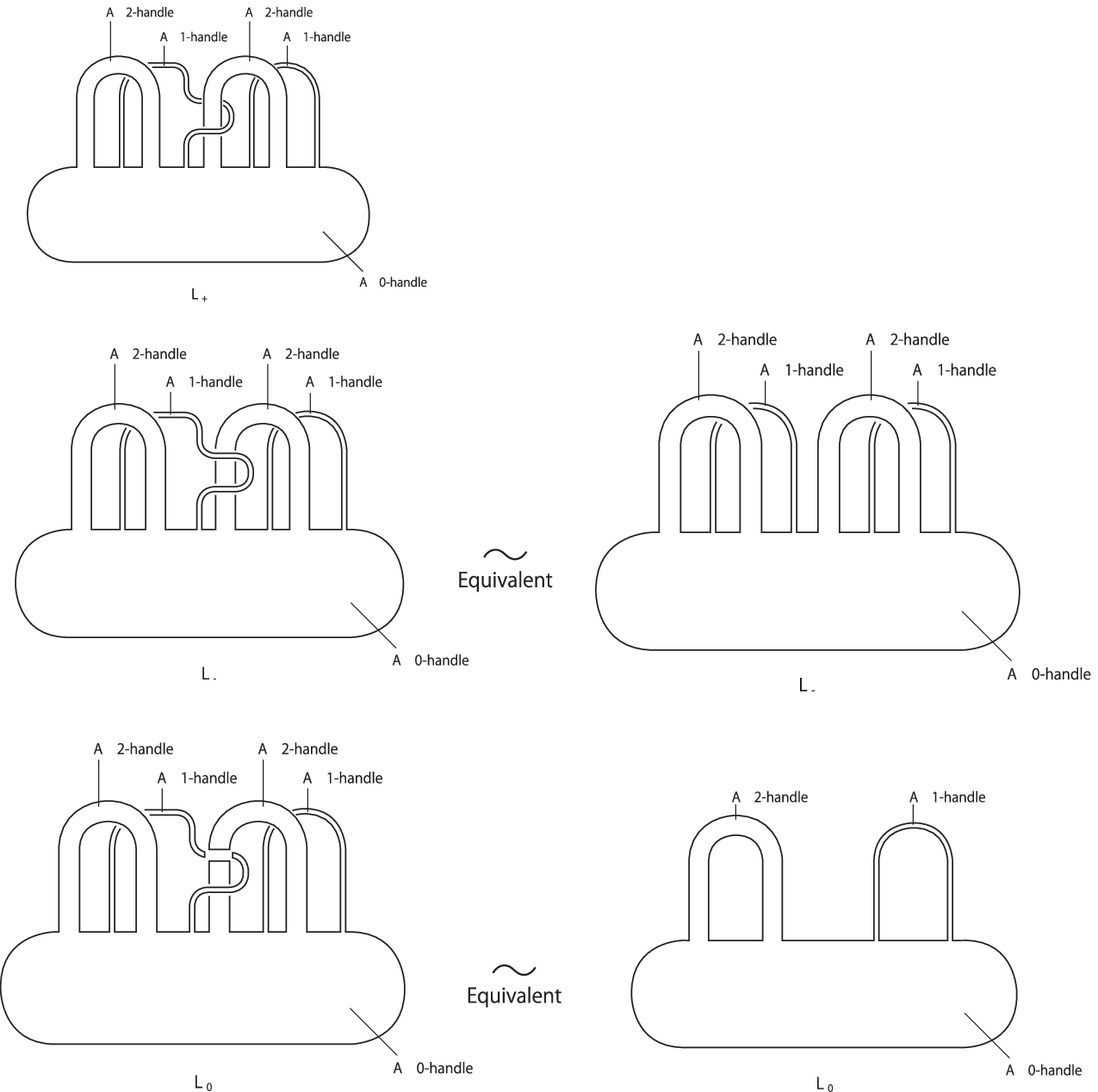}

Figure \ref{proo2}.2: A (1,2)-pass-move-triple $(L_+, L_-, L_0)$
\bigbreak\end{figure}
$L_+$, $L_-$ are diffeomophic to $S^2$. 
$L_0$ is diffeomophic to $S^2\amalg T^2$. 
A seifert hypersurface for $L_+$ (resp, $L_-$)   
is diffeomophic to  \newline 
$(S^2\x S^1)\sharp(S^2\x S^1)-\text{Int}B^3$, 
where $\sharp$ denotes the connected-sum.  
Note that it has a handle decomposition \newline 
(a 3-dimensional 0-handle)$\cup$(two 3-dimensional 1-handles)$\cup$(two 3-dimensional 2-handles). 
A seifert hypersurface for $L_0$ is diffeomophic to %\newline 
$(S^2\x D^1)\natural(D^2\x S^1)$,  
where $\natural$ denotes the boundary-connected-sum.  
Note that it has a handle decomposition \newline 
(a 3-dimensional 0-handle)$\cup$(a 3-dimensional 1-handle)$\cup$(a 3-dimensional 2-handle). 
\newline 
A 1-Alexander matrix for $L_+$ (resp, $L_-$, $L_0$)   
is 
$
\begin{pmatrix}
t&t-1\\
0&t 
\end{pmatrix}
$ 
(resp. 
$
\begin{pmatrix}
t&0\\
0&t
\end{pmatrix}
$, 
$(0)$).

There is a (1,2)-pass-move-triple $(L'_+, L'_-, L'_0)$ 
with the following properties (see Figure \ref{proo2}.3):  
\begin{figure}
   \includegraphics[width=8cm]{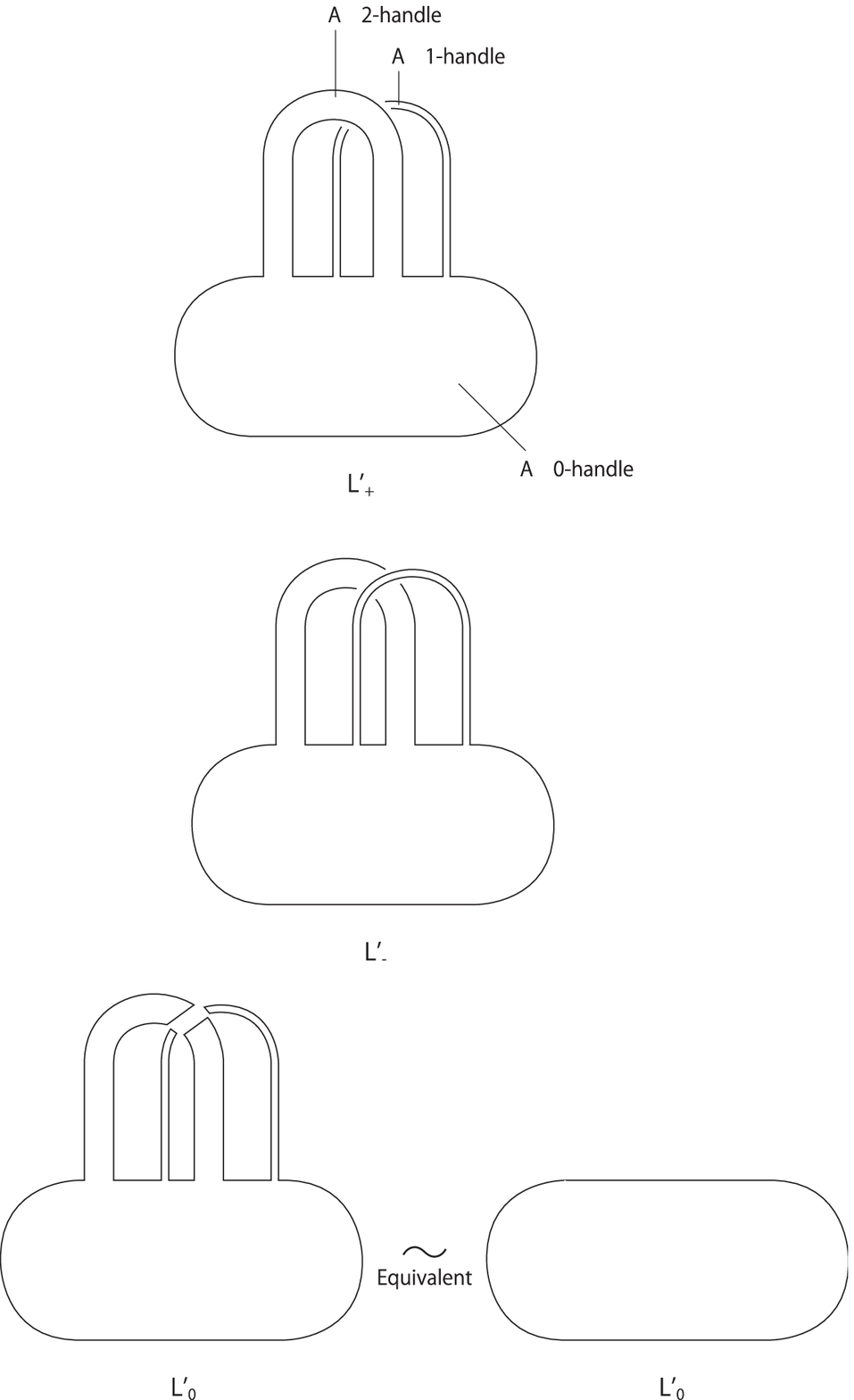}

Figure \ref{proo2}.3: A (1,2)-pass-move-triple $(L'_+, L'_-, L'_0)$ 
\vskip3mm
\end{figure}
$L'_+$, $L'_-$, and $L'_0$ are diffeomophic to $S^2$. 
A seifert hypersurface for $L_+$ (resp, $L_-$)   
is diffeomophic to $(S^2\x S^1)-\text{Int}B^3$.  
Note that it has a handle decomposition \newline 
(a 3-dimensional 0-handle)$\cup$(a 3-dimensional 1-handle)$\cup$(a 3-dimensional 2-handle). \newline
A seifert hypersurface for $L'_0$   
is diffeomophic to a 3-ball, 
which can be regarded a 3-dimensional 0-handle.  Note that $L'_0$ is a trivial 2-knot. 
A 1-Alexander matrix for $L'_+$ (resp, $L'_-$, $L'_0$)   
is 
$(t)$ 
(resp. 
$(1)$, 
$\phi$), 
where $\phi$ denotes the empty matrix.  

By these example,  Theorem \ref{Aa} holds.
\qed

\bigbreak
\noindent
{\bf Proof of Proposition \ref{Idaho}.}   
Propositions \ref{Idaho} follows from Proposition \ref{Mississippi}.   
\qed

\bigbreak
\noindent 
{\bf Proof of Proposition \ref{Mississippi}.}   
There are the following examples (see Figures \ref{proo2}.4 and \ref{proo2}.5).  
\begin{figure}
   \includegraphics[width=9cm]{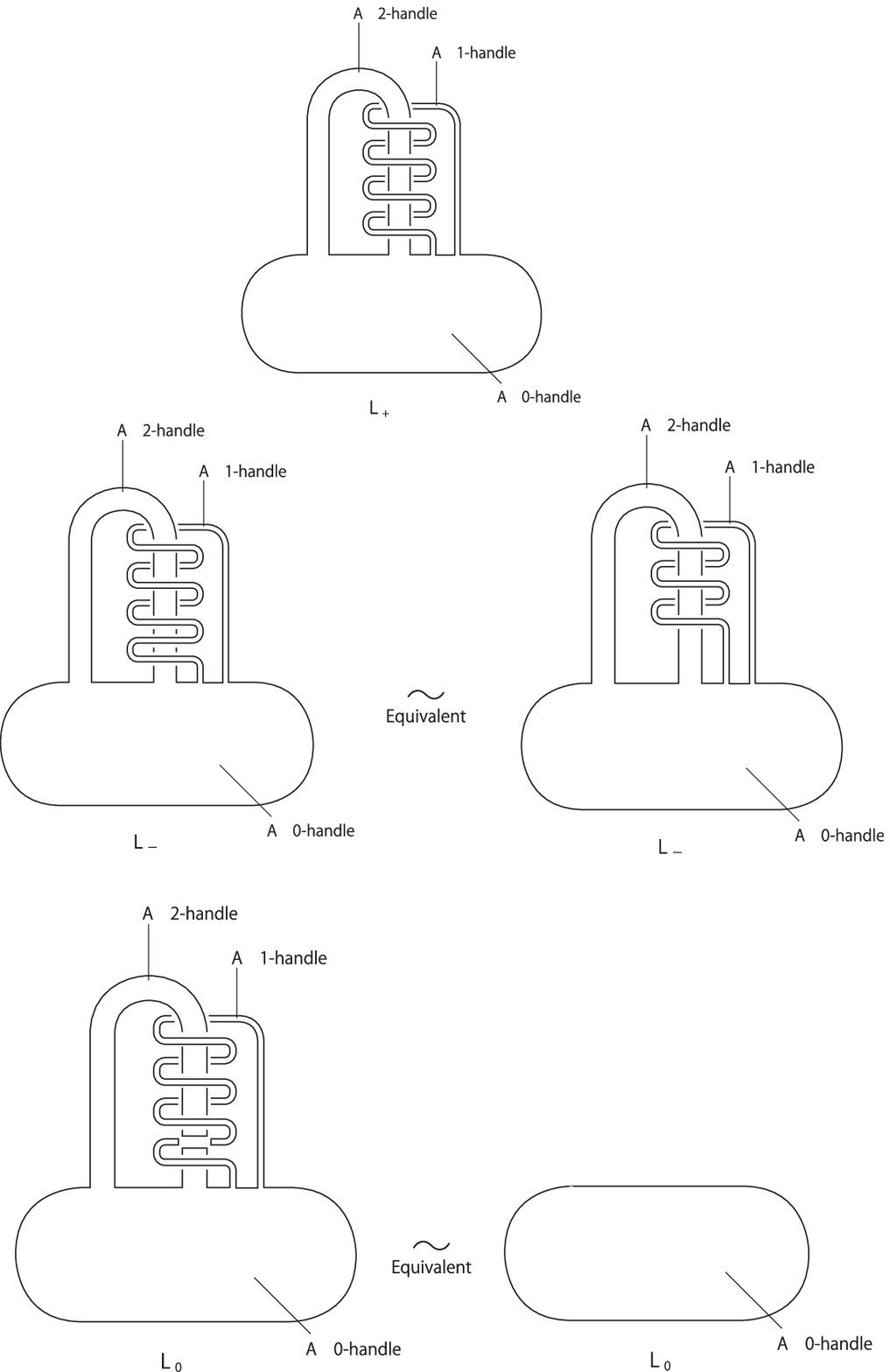}

Figure \ref{proo2}.4: A (1,2)-pass-move-triple $(L_+, L_-, L_0)$
\vskip3mm\end{figure}
\begin{figure} \includegraphics[width=93mm]{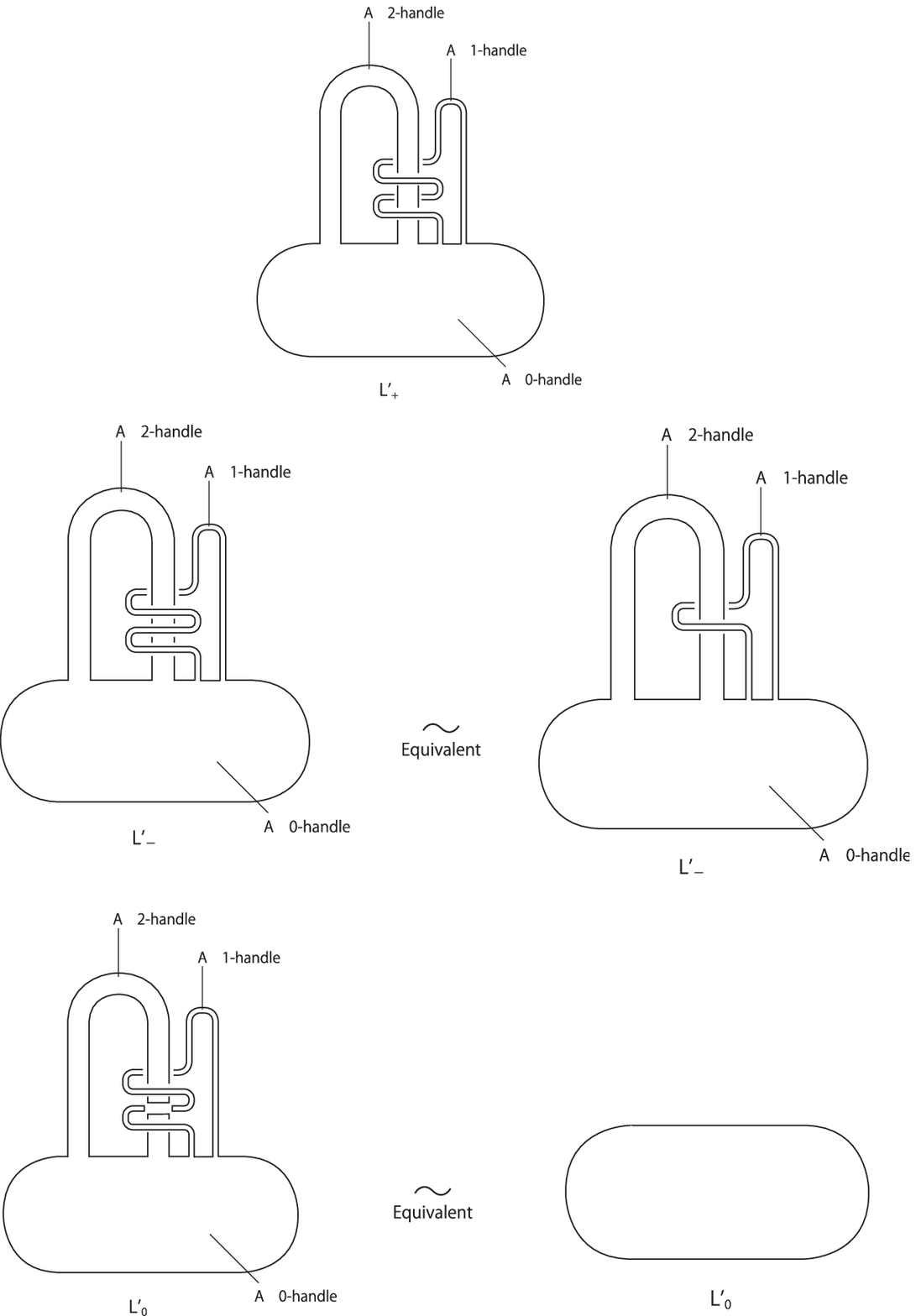}

Figure \ref{proo2}.5: A (1,2)-pass-move-triple $(L'_+, L'_-, L'_0)$  
\vskip3mm\end{figure}
Let $L_+$ (resp. $L_-$, $L'_+$, $L'_-$) be an $(S^2, T^2)$-link 
and bound a Seifert hypersurface which is diffeomorpphic to 
$(S^2\x D^1)\natural(D^2\x S^1)$.    
Note that it has a handle decomposition \newline 
(a 3-dimensional 0-handle)$\cup$(a 3-dimensional 1-handle)$\cup$(a 3-dimensional 2-handle)  \newline   
and that it is not diffeomorphic to $(S^2\x S^1)-\text{Int}B^3$. 
%$(S^1\x D^2)\natural(S^2\x D^1)$, 
%where $\natural$ denotes the boundary-connected-sum. 
We can assume that  a Seifert matrix of 
$L_+$ (resp. $L_-$, $L'_+$, $L'_-$) is a $1\x1$-matrix (4)  (resp. (3), (2), (1)). 
We can suppose that $L_0$ (resp. $L'_0$)  is a trivial 2-knot 
such that $(L_+,L_-,L_0)$ (resp. $(L'_+,L'_-,L'_0)$) is a $(1,2)$-pass-move-triple. 
\qed

%\np
%$\hat\Delta_{K}(t)={\rm{det}}(t^\frac{1}{2}\cdot S_{2k+1}(K)-t^\frac{-1}{2}\cdot N_{2k+1}(K))$

\bigbreak
\section{Proof of results in \S\ref{mhigh}}\label{proohigh}
\noindent
{\bf Proof of Theorem \ref{norm}.}   
Lemmas \ref{Indiana} and \ref{Iowa} imply Theorem \ref{norm}. 
\begin{lem}\label{Indiana}
Theorem \ref{norm} is true if $K$ is homotopy type equivalent to $S^{4k+1}$.  
\end{lem}

\noindent{\bf Proof of Lemma \ref{Indiana}.}
Recall that $S_{2k+1}(V)-N_{2k+1}(V)$ 
is represented by the intersection product on $H_{2k+1}(V;\Z)$. 
By the Poincar\'e duality 
${\rm{det}}(S_{2k+1}(V)-N_{2k+1}(V))=1.$  
Hence 
$\hat\Delta_{K}(1)=1$. 

By the Poincar\'e duality 
a $2k$-Alexander matrix induces an injective homomorphism 
on 
$H_{2k}(\amalg_{-\infty}^{\infty}V\x[-1,1];\Q)
\to
H_{2k}(\amalg_{-\infty}^{\infty}Y;\Q)$. 
Hence the $\Q[t,t^{-1}]$-balanced class of
${\rm{det}}(t\cdot S_{2k+1}(V)-N_{2k+1}(V))$  
is the $(2k+1)$-$\Q[t,t^{-1}]$-Alexander polynomial by \newline Propsoition \ref{square}.
Hence the $\Q[t,t^{-1}]$-balanced class of 
${\rm{det}}(t\cdot S_{2k+1}(V)-N_{2k+1}(V))$  
 is topological invariant of $K$. 
Hence 
the $\Q[t,t^{-1}]$-balanced class of 
${\rm{det}}(t\cdot S_{2k+1}(V)-N_{2k+1}(V))$  
 does not depend on the choice of $V$.

Let $V'$ be a Seifert hypersurface for $K$. 
Let $S_{2k+1}(V')$  be a positive $(2k+1)$-Seifert matrix for $K$ and  
$N_{2k+1}(V')$ its related negative Seifert matrix.    
Let ${\hat\Delta}'_{K}(t) \newline 
={\rm{det}}(t^\frac{1}{2}\cdot S_{2k+1}(V')-t^\frac{-1}{2}\cdot N_{2k+1}(V'))$. 
Recall ${\hat\Delta}_{K}(t)={\rm{det}}(t^\frac{1}{2}\cdot S_{2k+1}(V)-t^\frac{-1}{2}N_{2k+1}(V))$.   
It suffices to prove that ${\hat\Delta}'_{K}(t)={\hat\Delta}_{K}(t)$.

Since $V$ (resp. $V'$) is $(4k+2)$-dimensional 
and $\partial V$ (resp. (resp. $\partial V'$)) is PL homeomorphic to the standard sphere, 
rank$H_{2k+1}(V;\Z)$ (resp. rank$H_{2k+1}(V';\Z)$) is even. 
Therefore 
%By this fact and the previous paragraph %and 
%
%前のん, 
%
there is an integer $n$ such that %\newline
${\hat\Delta}'_{K}(t)=t^n\cdot\hat\Delta_{K}(t)$ $\cdot\cdot\cdot\cdot\cdot(*)$  holds. 
(Note that $\hat\Delta_{K}(1)=\hat\Delta'_{K}(1)=1$.)

By Propositions \ref{tsuika} and \ref{Mars}, 
$
N_{2k+1}(V)=^t\hskip-2mm S_{2k+1}(V). 
$
By the Poincar\'e duality 
the number of the rows of $S_{2k+1}(V)$
and that of the columns of it are the (same) even nonnegative integer. 
Hence we have the following: 
Let $M(t)={\rm{det}}(t^\frac{1}{2}\cdot S_{2k+1}(V)-t^\frac{-1}{2}\cdot N_{2k+1}(V))$. 
Then 
$M(t)= ^t\hskip-2mm M(t^{-1})$. 
Hence 
$\hat\Delta_{K}(t)$ (resp. ${\hat\Delta}'_{K}(t)$) 
has a form 
$$
\displaystyle
\sum^{\rho=l}_{\rho=0}a_\rho\cdot t^{\frac{\rho}{2}}
+
\sum^{\rho=l}_{\rho=0}a_\rho\cdot t^{-\frac{\rho}{2}}. 
$$
By this fact and the above identity $(*)$, we have  
${\hat\Delta}'_{K}(t)=\hat\Delta_{K}(t)$. 
%Hence $\hat\Delta_{K}(t)$ is a topological invariant of $K$ in this case. 

This completes the proof of Lemma \ref{Indiana}.  \qed

\begin{lem}\label{Iowa}
Theorem \ref{norm} is true if $K$ is homotopy type equivalent to $S^{2k+1}\x S^{2k}$. 
\end{lem}

\noindent{\bf Proof of Lemma \ref{Iowa}.}
Let $V$ be a Seifert hypersurface for $K$. 
There are two cases (see Note (2) to Definition \ref{4k+1}):  

\smallbreak\noindent
(I)  Any $2k$-Alexander matrix associated with $V$ induces an injective map \newline\hskip1cm 
on $H_{2k}(\amalg_{-\infty}^{\infty}V\x[-1,1];\Q)\to H_{2k}(\amalg_{-\infty}^{\infty}Y;\Q)$. 
%The homomorphism map $H_{2k}(K;\Q)\to H_{2k}(V;\Q)$ induced by the natural inclusion map is not injective. 

\smallbreak\noindent
(II) 
No $2k$-Alexander matrix associated with $V$ induces an injective map \newline\hskip1cm 
on $H_{2k}(\amalg_{-\infty}^{\infty}V\x[-1,1];\Q)\to H_{2k}(\amalg_{-\infty}^{\infty}Y;\Q)$. 
%The homomorphism map $H_{2k}(K;\Q)\to H_{2k}(V;\Q)$ induced by the natural inclusion map is injective. 

\smallbreak
Lemmas \ref{Kansas} and \ref{Kentucky} imply Lemma \ref{Iowa}. 

\begin{lem}\label{Kansas} Lemma \ref{Iowa} is true in the case \rm{(II)}. \end{lem}

\noindent{\bf Proof of Lemma \ref{Kansas}.}   
Since $V$ satisfies (II), 
the  normalized Alexander polynomial of $K$ defined by using $V$ 
is zero by  Definition \ref{4k+1}. 
By Definition \ref{North Dakota} 
%the definition of the $\Q[t,t^{-1}]$-Alexander polynomial, 
the $\Q[t,t^{-1}]$-Alexander polynomial of $K$ defined by using $V$ 
is the $\Q[t,t^{-1}]$-balanced class of zero. 
 
Let $V'$ be another Seifert hypersurface for $K$. 
By Definition \ref{North Dakota} 
%the definition of the $\Q[t,t^{-1}]$-Alexander polynomial, 
the $\Q[t,t^{-1}]$-Alexander polynomial of $K$ defined by using $V'$  
is the $\Q[t,t^{-1}]$-balanced class of zero 
even if $V'$ satisfies (I) not (II).  
Therefore, by Definitions \ref{North Dakota} and \ref{4k+1} and Proposition \ref{square}  
%the definition of the $\Q[t,t^{-1}]$-(resp. $\Z[t,t^{-1}]$-) Alexander polynomial,  
the normalized Alexander polynomial of $K$ defined by using $V'$ is zero. 

This completes the proof of Lemma \ref{Kansas}. \qed

\begin{lem}\label{Kentucky} Lemma \ref{Iowa} is true in the case \rm{(I)}. \end{lem}

\noindent{\bf Proof of Lemma \ref{Kentucky}.}
In the same manner as  
written in the first part of Proof of Theorem \ref{Z}
it suffices to prove the following case: 
$V$ and $V'$ are Seifert hypersurfaces for $K$. 
$V'$ is obtained from $V$ by a surgery by using an embedded 
$(4k+3)$-dimensional $i$-handle $h^i$ ($1\leq i\leq 4k+2$)
This surgery may change a $(2k+1)$-Alexander matrix associated with $V$ for $K$ 
only if $i=2k+1, 2k+2$.  

Lemma \ref{Kansas} implies that if $V'$ satisfies (II), Lemma \ref{Kentucky} holds. 
Hence it suffices to prove the case where $V'$ satisfies (I). 

Hence both $V$ and $V'$ satisfy (II). 
The dual handle of $h^{2k+2}$ is a $(4k+3)$-dimensional $(2k+1)$-handle. 
Hence it suffices to prove the $i=2k+2$ case. 
Call the core of the attaching part of $h^{2k+2}$, $C$. 

There are two cases: 

\smallbreak\noindent
(1) $[C]\in H_{2k+1}(V,\Z)$ is order finite. 

\smallbreak\noindent
(2) $[C]\in H_{2k+1}(V,\Z)$ is order infinite.

\smallbreak
The case (2) is divided into two cases: 

\smallbreak\noindent
(2-1) $[C]\in H_{2k+1}(V,\Z)$ is order infinite. For all $(2k+1)$-cycle $\alpha$, 
the intersection product $[C]\cdot\alpha=0$. 
 
\smallbreak\noindent
(2-2) $[C]\in H_{2k+1}(V,\Z)$ is order infinite. 
There is a $(2k+1)$-cycle $\alpha$ 
such that the intersection product $[C]\cdot\alpha$ is nonzero.

\bigbreak
Lemmas \ref{Missouri}, \ref{Montana}, and \ref{Nebraska} 
imply Lemma \ref{Kentucky}. 

\begin{lem}\label{Missouri} 
Lemma \ref{Kentucky} holds in  the case $(1).$ 
\end{lem}

\noindent{\bf Proof of Lemma \ref{Missouri}.} 
This surgery does not change a $(2k+1)$-Alexander matrix associated with $V$ for $K$. 
This completes the proof of Lemma \ref{Missouri}.\qed

\begin{lem}\label{Montana} 
Lemma \ref{Kentucky} holds in  the case {\rm $(2$-$1)$}. 
\end{lem}

\noindent{\bf Proof of Lemma \ref{Montana}.} 
There is an Alexander matrix associated with $V$ which has a row (or column) all of whose elements are zero. Hence the $\Q[t,t^{-1}]$-Alexander polynomial of $K$ is 
the $\Q[t,t^{-1}]$-balanced class of zero. 
By Definitions \ref{North Dakota} and \ref{4k+1} 
%By the definition of  the $\Q[t,t^{-1}]$-Alexander polynomial and that of the $\Z[t,t^{-1}]$-Alexander polynomial, 
the  normalized Alexander polynomial is zero. 
%%%{\color{cyan}{これはLemmaにしたほうが良いかな　いや、defが、それぞれ微妙にちゃうから、まとめれん}}  
%By the Poincar\'e duality and the Mayor-Vietoris exact sequence, the case (2-1) is the case (1) after replacing  $V$ and $V'$.  {\color{cyan}{これは間違い｝  

This completes the proof of Lemma \ref{Montana}. \qed

%{\color{cyan}これ、どこ、置こ。  We can suppose that the homomorphism map $H_{2k}(K;\Q)\to H_{2k}(V';\Q)$ induced by the natural inclusion map is not injective. Reason:  If it is injective, replave $V$ and $V'$. By Lemma \ref{Kansas}, Lemma \ref{Kentucky} holds.    } 

\begin{lem}\label{Nebraska} 
Lemma \ref{Kentucky} holds in  the case {\rm$(2$-$2)$}. 
\end{lem}

\noindent{\bf Proof of Lemma \ref{Nebraska}.}  
Let ${\bar h}^{2k+1}$ be the dual handle of  
the $(4k+3)$-dimensional $(2k+2)$-handle $h^{2k+2}$. 
Let $C'$ be the core of the attaching part of ${\bar h}^{2k+1}$. 
Note that \newline 
$(V\x[0,1])\cup h^{2k+2}$ is diffeomorphic to 
$(V'\x[0,1])\cup {\bar h}^{2k+1}$. 

Consider the exact sequence by a pair 
$((V\x[0,1])\cup h^{2k+2}, V)$, where 
we regard $V$ as $V\x\{0\}$: 
$$
\cdot\cdot\cdot\to H_*(V;\Z)\to 
H_*(V\x[0,1])\cup h^{2k+2};\Z)\to 
H_*(V\x[0,1])\cup h^{2k+2}, V;\Z)\to\cdot\cdot\cdot 
$$
and 
the exact sequence by a pair 
$((V'\x[0,1])\cup {\bar h}^{2k+1}, V')$, where 
we regard $V'$ as $V'\x\{0\}$: 
$$
\cdot\cdot\cdot\to H_*(V';\Z)\to 
H_*(V'\x[0,1])\cup {\bar h}^{2k+1};\Z)\to 
H_*(V'\x[0,1])\cup {\bar h}^{2k+1}, V';\Z)\to\cdot\cdot\cdot. 
$$

By the existence of the $(2k+1)$-cycle $\alpha$, 
$[C']\in H_{2k}(V'; \Z)$ is order finite.

Let $\xi\in H_{2k+1}(V;\Z)$ be a non-divisible element associated with $[C]$. 
Since $H_{2k}(\partial V;\Z)$ is torsion-free and the intersection product $\xi\cdot\alpha\neq0$, 
there is a $(2k+1)$-cycle $\eta\in H_{2k+1}(V;\Z)$ such that 
$\eta\cdot\xi=1$. 
We can suppose that $\eta$ is a non-divisible element associated with $\alpha$. 
Therefore we have the following: 
 $A(t)$ (resp. $A'(t)$) 
is ${\rm{det}}(t^\frac{1}{2}\cdot S_{2k+1}(V)-t^\frac{-1}{2}\cdot N_{2k+1}(V))$ 
(resp.  ${\rm{det}}(t^\frac{1}{2}\cdot S_{2k+1}(V')-t^\frac{-1}{2}\cdot N_{2k+1}(V'))$) 
an $(2k+1)$-Alexander matrix associated with $V$ (resp. $V'$).  
We have 

$$A'(t)= 
\arraycolsep5pt
\left(
\begin{array}{@{\,}cc|cccc@{\,}}
0&t^\frac{1}{2}& 0&\cdot&\cdot& 0\\
-t^\frac{-1}{2}& 0&a_{23}&\cdot&\cdot& a_{2n}\\
\hline
0&a_{23}&&&&\\
\cdot&\cdot&&\multicolumn{3}{c}{\raisebox{-10pt}[0pt][0pt]{\Huge \hskip-3mm$A(t)$}}\\
\cdot&\cdot&&&&\\
0&a_{2n}&&&&\\
\end{array}
\right)
$$

or

$$A'(t)= 
\arraycolsep5pt
\left(
\begin{array}{@{\,}cc|cccc@{\,}}
0&-t^\frac{1}{2}& 0&\cdot&\cdot& 0\\
t^\frac{-1}{2}& 0&a_{23}&\cdot&\cdot& a_{2n}\\
\hline
0&a_{23}&&&&\\
\cdot&\cdot&&\multicolumn{3}{c}{\raisebox{-10pt}[0pt][0pt]{\Huge \hskip-3mm$A(t)$}}\\
\cdot&\cdot&&&&\\
0&a_{2n}&&&&\\
\end{array}
\right)
$$

This completes the proof of Lemma \ref{Nebraska}. \qed 

This completes the proof of Lemma \ref{Kentucky}. \qed

This completes the proof of Lemma \ref{Iowa}. \qed

This completes the proof of Theorem \ref{norm}. \qed

\bigbreak
\noindent
{\bf Proof of Theorem \ref{Tokyo}.} 
In \cite[Theorem 3.3]{Ogasa09} we proved the following: 
There is a Seifert hypersurface $V_*$ for $K_*(*=+,-,0)$ 
with an associated $(2k+1)$-Seifert matrix $S_{2k+1}(V_*)$ 
and its related $(2k+1)$-negative Seifert matrix $N_{2k+1}(V_*)$ with the following properties:
\smallbreak
\noindent
(i) 
$${\rm{det}}(t^\frac{1}{2}\cdot S_{2k+1}(V_+)-t^\frac{-1}{2}\cdot N_{2k+1}(V_+))$$
$$-{\rm{det}}(t^\frac{1}{2}\cdot S_{2k+1}(V_-)-t^\frac{-1}{2}\cdot N_{2k+1}(V_-))$$
$$=(t-1)\cdot{\rm{det}}(t^\frac{1}{2}\cdot S_{2k+1}(V_0)-t^\frac{-1}{2}\cdot N_{2k+1}(V_0)).$$
\smallbreak
\noindent
(ii)
$S_{2k+1}(V_+)$ and 
$S_{2k+1}(V_-)$ are $2\nu\x2\nu$-matrices  $(\nu\in\N).$  
$S_{2k+1}(V_+)$ is a$(2\nu-1\x2\nu-1)$-matrix. 
\smallbreak
\noindent
(iii)
The $2k$-Alexander matrix associated with each Seifert hypersurface
defines an  injective map on 
on $H_{2k}(\amalg_{-\infty}^{\infty}V\x[-1,1];\Q)\to H_{2k}(\amalg_{-\infty}^{\infty}Y;\Q)$.

\smallbreak
Hence 
$$\hat\Delta_{K_+}(t)-\hat\Delta_{K_-}(t)
=(t^\frac{1}{2}-t^\frac{-1}{2})\cdot\hat\Delta_{K_0}(t).$$

This completes the proof of Theorem \ref{Tokyo}.  \qed

\bigbreak
\noindent
{\bf Proof of Proposition \ref{hongo}.}
In the same manner as written in the first part of Proof of Theorem \ref{Z} 
it suffices to prove the following case: 
$V$ and $V'$ are Seifert hypersurfaces for $K$. 
$V'$ is obtained from $V$ by a surgery by using an embedded 
$(4k+3)$-dimensional $i$-handle $h^i$ ($1\leqq i\leqq 4k+2$). 
The pseudo-twinkling number may change only if $i\newline=2k+1,2k+2$. 

The dual handle of $h^{2k+2}$ is a $(4k+3)$-dimensional $(2k+1)$-handle. 
Therefore it suffices to prove the following two cases under the condition $i=2k+2$. 

%\smallbreak\noindent (1) 
%Any $2k$-Alexander matrix associated with $V$ induces an injective map 

%\hskip1mm
%on $H_{2k}(\amalg_{-\infty}^{\infty}V\x[-1,1];\Q)\to H_{2k}(\amalg_{-\infty}^{\infty}Y;\Q)$. 

%\smallbreak\noindent (2) 
%No $2k$-Alexander matrix associated with $V$ induces an injective map 

%\hskip1mm 
%on $H_{2k}(\amalg_{-\infty}^{\infty}V\x[-1,1];\Q)\to H_{2k}(\amalg_{-\infty}^{\infty}Y;\Q)$. 

%By the Poincar\'e duality, the Mayor-Vietoris exact sequence and 
%the exact sequence by a pair, we have (1)$\Leftrightarrow(1^\circ)$ %これがさっきの（I）と同じってほんまか？
%and  (2)$\Leftrightarrow(2^\circ)$. %これがさっきの（II）と同じってほんまか？

\smallbreak\noindent
(1) %($1^\circ$) 
There is a non-divisible $(2k+1)$-cycle $\tau\subset V$ such that for any $(2k+1)$-cycle $\alpha\subset V$  
the intersection product $\tau\cdot\alpha$ in $V$ is zero 

\smallbreak\noindent
(2) %($2^\circ$)  
There is not such a cycle as in (1). %($1^\circ$). 

\smallbreak\noindent
%Hence ($1^\circ$) is called (1) and ($2^\circ$) is called (2). 

By Poincar\'e duality and Mayor-Vietoris exact sequence, 
 $\tau$ is a non-divisible cycle in $V$ associated with $*\x S^{2k+1}$ in $K=\partial V$.

The above two cases (1) and (2) are divide into four cases. 

\smallbreak\noindent
(1-1) $V$ satisfies (1). $V'$ satisfies the condition made from (1) 

\hskip4mm
by  replacing $V$ with $V'$ in (1). 

\smallbreak\noindent
(1-2) $V$ satisfies (1). $V'$ satisfies the condition made from (2) 

\hskip4mm
by replacing $V$ with $V'$ in (2).   

\smallbreak\noindent
(2-1) $V$ satisfies (2). $V'$ satisfies the condition made from (1) 

\hskip4mm
by replacing $V$ with $V'$ in (1).   

\smallbreak\noindent
(2-2) $V$ satisfies (2). $V'$ satisfies the condition made from (2) 

\hskip4mm
by replacing $V$ with $V'$ in (2).   

\smallbreak\noindent  
Lemmas \ref{Nevada}, \ref{New Hampshire}, \ref{New Jersey}, and 
\ref{New Mexico} imply Proposition \ref{hongo}.

\begin{lem}\label{Nevada}
Proposition \ref{hongo} holds in the case $(2\text{-}2).$ 
 \end{lem}

\noindent{\bf Proof of Lemma \ref{Nevada}.}
By Definition \ref{twinkling}   
 the pseudo-twinkling number defined by using $V$ (resp. $V'$) is zero. 
\qed

\begin{lem}\label{New Hampshire}
Proposition \ref{hongo} holds in the case $(2\text{-}1).$  
\end{lem}

\noindent{\bf Proof of Lemma \ref{New Hampshire}.}
%By the definition of the pseudo-twinkling number, the pseudo-twinkling number defined by using $V$ is zero. 
When we obtain $V'$ from $V$ by using $h^{2k+2}$, 
there does not appear $\tau'$ in $V'$ as in Definition \ref{4k+1}.
The case (2-1) does not occur. 
%Hence Lemma \ref{Nevada} implies \ref{New Hampshire}. 
\qed

\bigbreak
Lemmas \ref{Nevada} and \ref{New Hampshire} and their proof imply Claim \ref{New York}. 
\begin{cla}\label{New York}
The pseudo-twinkling number is zero in the case $(2).$ 
\end{cla}

\begin{lem}\label{New Jersey}
Proposition \ref{hongo} holds in the case $(1\text{-}2).$ 
\end{lem}

\noindent{\bf Proof of Lemma \ref{New Jersey}.}
% Lemma \ref{New Hampshire} imply Lemma \ref{New Jersey}.  これ間違い
%
By Definition \ref{twinkling}  
 the pseudo-twinkling number defined by using $V'$ is zero. 

%Take $\tau$ in $V$ as in Definition \ref{4k+1}.  
Let $C$ be the attaching part  of $h^{2k+2}$.  
Note that the Seifert pairing $s(C,C)=0$. 

Under the condition (1-2), 
$\tau$ must be a non-divisible $(2k+1)$-cycle associated with $[C]$.  
Hence $s(\tau',\tau')=0$.  
\qed

\begin{lem}\label{New Mexico}
Proposition \ref{hongo} holds in the case $(1\text{-}1).$ 
\end{lem}

\noindent{\bf Proof of Lemma \ref{New Mexico}.}
%Take $\tau$ in $V$ as in Definition \ref{4k+1}.  
The pseudo-twinkling number defined by using $V$ is $s(\tau, \tau)$. 
%Recall that $V'$ is obtained from $V$ by a surgery by using a 
%$(4k+3)$-dimensional $i$-handle  $h^i$ ($i=2k+1, 2k+2$). 
Let $C$ be the core of the attaching part of $h^{2k+2}$. 
%Note that $C$ is a $(i-1)$-sphere. 
There are two cases. 

\smallbreak\noindent(i) 
$\tau$ is a non-divisible cycle in $V$ associated with $C$ . 

\smallbreak\noindent(ii) 
Else. %The other case than (i). 

\smallbreak 
In the case (i),  
$V$ and $V'$ satisfy (1-2) not (1-1). 
 
The cases (ii) 
%the disjoint union of the tubular neighborhood of each cycle which makes $\tau$
%does not change after the surgery by using $h^i$ 
%and $\tau$ is a non-divisible cycle in $V'$ associated with $*\x S^{2k+1}$ in $K=\partial V'$.
%Hence the pseudo-twinkling number defined by using $V'$ is $s(\tau, \tau)$ again. 
follows from Theorem \ref{aletwi} and its proof as written below because 
its proof does not depend of the choice of $V$.

This completes the proof of Lemma \ref{New Mexico}. \qed

This completes the proof of Proposition \ref{hongo}. \qed

\bigbreak
\noindent
{\bf Proof of Theorem \ref{aletwi}.} 
Let $V$ be a Seifert hypersurface for $K$. 

In the case (2) of Proof of Proposition \ref{hongo}.  
By Claim \ref{New York}   the pseudo-twinkling number is zero. 
In this case 
no $2k$-Alexander matrix associated with $V$ induces an injective map.  
By Definition \ref{4k+1} %of the $\Z[t,t^{-1}]$-Alexander polynomial 
the normalized Alexander polynomial  is zero. 
Hence Theorem \ref{aletwi} holds in this case.

In the case (1) of Proof of Proposition \ref{hongo}.  
Let $\{\alpha_1,...,\alpha_\nu\}$ be a basis of $H_{2k+1}(V;\Z)/\text{Tor}$. We can suppose that $\alpha_1=\tau$. 
Then  
$t^{\frac{1}{2}}S-t^{\frac{-1}{2}}\hskip1mm{^t\hskip-1mm S}$
is written as follows, 
where $a_{\*\#}$ is an integer and $a_{11}$ is the pseudo-twinkling number. 

$$
\arraycolsep5pt
\left(
\begin{array}{@{\,}c|cccc@{\,}}
(t^\frac{1}{2}-t^\frac{-1}{2})\cdot a_{11}&(t^\frac{1}{2}-t^\frac{-1}{2})\cdot a_{12}&\cdot&\cdot&(t^\frac{1}{2}-t^\frac{-1}{2})\cdot a_{1*}\\
\hline
(t^\frac{1}{2}-t^\frac{-1}{2})\cdot a_{21}&&&&\\
\cdot&\multicolumn{4}{c}{\raisebox{-10pt}[0pt][0pt]{\Huge $Q(t)$}}\\
\cdot&&&&\\
(t^\frac{1}{2}-t^\frac{-1}{2})\cdot a_{*1}&&&&\\
\end{array}
\right)
$$

$$=(t^\frac{1}{2}-t^\frac{-1}{2})\cdot
\arraycolsep5pt
\left(
\begin{array}{@{\,}c|cccc@{\,}}
a_{11}& a_{12}&\cdot&\cdot& a_{1*}\\
\hline
(t^\frac{1}{2}-t^\frac{-1}{2})\cdot a_{21}&&&&\\
\cdot&\multicolumn{4}{c}{\raisebox{-10pt}[0pt][0pt]{\Huge $Q(t)$}}\\
\cdot&&&&\\
(t^\frac{1}{2}-t^\frac{-1}{2})\cdot a_{*1}&&&&\\
\end{array}
\right)
$$

Hence 
$\left.\displaystyle
\frac{t^{\frac{1}{2}}S-t^{\frac{-1}{2}}\hskip1mm{^t\hskip-1mm S} }
{t^\frac{1}{2}-t^\frac{-1}{2}}
\right|_{t=1}$
is written as follows. 

$$
\arraycolsep5pt
\left(
\begin{array}{@{\,}c|cccc@{\,}}
a_{11}&  a_{12}&\cdot&\cdot& a_{1*}\\
\hline
0&&&&\\
\cdot&\multicolumn{4}{c}{\raisebox{-10pt}[0pt][0pt]{\Huge $Q(1)$}}\\
\cdot&&&&\\
0&&&&\\
\end{array}
\right)
$$

\noindent 
Since $Q(1)$ is a nonsingular matrix and its determinant is $+1$,   
 $\left.\displaystyle\frac{\hat\Delta_{K}(t)}{t^\frac{1}{2}-t^\frac{-1}{2}}\right|_{t=1}=a_{11}$.

\begin{note}\label{Hawaii}
If we define $V$ as in (1), 
$\iota:H_{2k+1}(\partial V;\Z)\to H_{2k+1}(V;\Z)$ has the property in (2).  
 
\smallbreak\noindent
(1)
Take $S^{2k+1}\x S^{2k+1}$.  
Let $f$ be an embedding map 
$S^{2k+1}\hookrightarrow S^{2k+1}\x S^{2k+1}$.  
Let $f(S^{2k+1})\subset S^{2k+1}\x S^{2k+1}$. 
Suppose that the induced map \newline
$f_*:H_{2k+1}(S^{2k+1};\Z)\to H_{2k+1}(S^{2k+1}\x S^{2k+1};\Z)$
is $\Z\to\Z\oplus\Z$ with $1\mapsto (n,0)$, where $|n|\newline>1$.  
Here, we fix a generator of $H_{2k+1}(S^{2k+1};\Z)$ and 
that of 
$H_{2k+1}(S^{2k+1}\x S^{2k+1};\Z)$. 
Let $V$ be $(S^{2k+1}\x S^{2k+1})-{\rm Int}N(f(S^{2k+1}))$, 
where $N(f(S^{2k+1}))$ is the tubular neighborhood of $f(S^{2k+1})$ in  $S^{2k+1}\x S^{2k+1}$.  

\smallbreak\noindent 
(2) 
Let $g$ be a generator of  $H_{2k+1}(\partial V;\Z)\cong\Z$. 
Then $\rho(g)$ is a divisible cycle $\in H_{2k+1}(V;\Z)$. 
\end{note}

\bigbreak
\section{A problem}\label{problem}
\begin{prob}\label{ichiika}
(1) If an invariant of 2-dimensional oriented closed submanifold $\subset S^4$ 
satisfies the identity in Theorem \ref{two} associated with the $(1,2)$-pass-move, 
then is it essentially the $\Z[t,t^{-1}]$-Alexander polynomial? 

\smallbreak\noindent$(2)$
If an invariant of $(4k+1)$-dimensional submanifolds $\subset S^{4k+3}$ 
whose homotopy type is $S^{4k+1}$ or 
$S^{2k}\x S^{2k+1}$ satisfies the identity in Theorem \ref{Tokyo} as written there, 
then is it essentially 
the normalized Alexander polynomial 
(resp. the $\Z[t,t^{-1}]$-Alexander polynomial)? 
\end{prob} 

\noindent{\bf Note.} 
%(1) 
If the answer to 
Problem \ref{ichiika}.(1) is positive, 
it is a new characterization of the $\Z[t,t^{-1}]$-Alexander polynomial of 
2-dimensional closed oriented submanifolds in $S^4$.  
If the answer is negative, we may encounter a new invariant. 

If the answer to Problem \ref{ichiika}.(2) is positive, 
it is a new characterization of 
the normalized Alexander polynomial 
(resp. the $\Z[t,t^{-1}]$-Alexander polynomial) of 
high dimensional knots in the case.   
If the answer is negative, we may encounter a new invariant. 

%
%\smallbreak\noindent$(2)$
There arise similar problems on 
the local-move-identities in \cite{KauffmanOgasa, Ogasa09} to Problem \ref{ichiika}.

\bigbreak\noindent
Eiji Ogasa \quad  Computer Science, Meijigakuin University, Yokohama, Kanagawa, 244-8539, Japan 
\quad pqr100pqr100@yahoo.co.jp  \quad
ogasa@mail1.meijigkakuin.ac.jp

\end{document}